\documentclass[11pt]{article}

\usepackage[utf8]{inputenc}
\usepackage{amsmath,amssymb,epsfig,bbm}
\usepackage{stmaryrd}
\usepackage{comment}
\usepackage[dvips]{color}
\usepackage[T1]{fontenc}

\usepackage{todonotes}
\usepackage{enumitem}
\setlist{nolistsep}

\newtheorem{defi}{Definition}
\newtheorem{prop}[defi]{Proposition}
\newtheorem{theo}[defi]{Theorem}
\newtheorem{conj}[defi]{Conjecture}
\newtheorem{lemm}[defi]{Lemma}
\newtheorem{coro}[defi]{Corollary}
\newtheorem{rema}[defi]{Remark}
\newtheorem{exem}[defi]{Example}
\newtheorem{exems}[defi]{Examples}

\newcommand{\bdefi}{\begin{defi}}
\newcommand{\edefi}{\end{defi}}
\newcommand{\bprop}{\begin{prop}}
\newcommand{\eprop}{\end{prop}}
\newcommand{\btheo}{\begin{theo}}
\newcommand{\etheo}{\end{theo}}
\newcommand{\blemm}{\begin{lemm}}
\newcommand{\brema}{\begin{rema}}
\newcommand{\erema}{\end{rema}}
\newcommand{\bexer}{\begin{exem}}
\newcommand{\eexer}{\end{exem}}
\newcommand{\bexems}{\begin{exems}}
\newcommand{\eexems}{\end{exems}}
\newcommand{\bconj}{\begin{conj}}
\newcommand{\econj}{\end{conj}}
\newcommand{\elemm}{\end{lemm}}
\newcommand{\bcoro}{\begin{coro}}
\newcommand{\ecoro}{\end{coro}}
\newcommand{\dem}{\noindent{\bf Proof. }}
\newcommand{\rem}{\noindent{\bf Remark. }}

\usepackage{mathrsfs}
\renewcommand\mathcal{\mathscr}

\newcommand{\M}{{\cal M}}
\newcommand{\N}{{\cal N}}

\renewcommand{\H}{{\cal H}}
\newcommand{\OOO}{{\cal O}}
\newcommand{\C}{{\cal C}}
\newcommand{\I}{{\cal I}}
\newcommand{\Q}{{\cal Q}}

\newcommand{\maths}[1]{{\mathbb #1}}  

\newcommand{\RR}{\maths{R}}
\newcommand{\NN}{\maths{N}}
\newcommand{\CC}{\maths{C}}
\newcommand{\QQ}{\maths{Q}}

\newcommand{\SSS}{\maths{S}}

\newcommand{\HH}{\maths{H}}

\newcommand{\ZZ}{\maths{Z}}
\newcommand{\PP}{\maths{P}}

\newcommand{\aaa}{{\mathfrak a}}

\newcommand{\ppp}{{\mathfrak p}}

\newcommand{\ccc}{{\mathfrak c}}
\newcommand{\mmm}{{\mathfrak m}}

\newcommand{\weakstar}{\overset{*}\rightharpoonup}
\newcommand{\ra}{\rightarrow}
\newcommand{\bs}{\backslash}

\newcommand{\ov}[1]{{\overline #1}} 
\newcommand{\wt}[1]{{\widetilde{#1}}}
\newcommand{\wh}[1]{{\widehat{#1}}}

\newcommand{\ga}{\gamma}
\newcommand{\Ga}{\Gamma}

\newcommand{\cqfd}{\hfill$\Box$}

\newcommand{\card}{{\operatorname{Card}}}
\renewcommand{\Re}{{\operatorname{Re}}}

\newcommand{\Vol}{\operatorname{Vol}}

\newcommand{\covol}{\operatorname{Covol}}
\newcommand{\id}{\operatorname{id}}
\newcommand{\PSL}{\operatorname{PSL}}
\newcommand{\SL}{\operatorname{SL}}

\newcommand{\arcosh}{\operatorname{argcosh}}

\newcommand{\axis}{\operatorname{Axis}}
\newcommand{\bigO}{\operatorname{O}}
\newcommand{\smallo}{\operatorname{o}}

\renewcommand{\log}{\operatorname{ln}}
\newcommand{\Leb}{\operatorname{Leb}}

\newcommand{\hdr}{{\HH}^2_\RR}
\newcommand{\htr}{{\HH}^3_\RR}
\newcommand{\hcr}{{\HH}^5_\RR}
\newcommand{\hnr}{{\HH}^n_\RR}

\newcommand{\SLO}{\operatorname{SL}_{2}(\OOO)}

\newcommand{\PSLO}{\operatorname{PSL}_{2}(\OOO)}
\newcommand{\SLOK}{\operatorname{SL}_{2}(\OOO_K)}
\newcommand{\PSLOK}{\operatorname{PSL}_{2}(\OOO_K)}
\newcommand{\SLH}{\operatorname{SL}_{2}(\HH)}
\newcommand{\PSLH}{\operatorname{PSL}_{2}(\HH)}
\newcommand{\SLC}{\operatorname{SL}_{2}(\CC)}
\newcommand{\PSLC}{\operatorname{PSL}_{2}(\CC)}
\newcommand{\SLZ}{\operatorname{SL}_{2}(\ZZ)}
\newcommand{\PSLZ}{\operatorname{PSL}_{2}(\ZZ)}

\newcommand{\autom}{\operatorname{SU}_f(\OOO_K)}

\newcommand{\automH}{\operatorname{SU}_f(\OOO)}

\newcommand{\const}{\iota}
\newcommand{\discr}{\operatorname{Disc}}

\newcommand{\SpO}{\operatorname{Sp}_1(\OOO)}

\newcommand{\tr}{\operatorname{\tt tr}}
\newcommand{\n}{\operatorname{\tt n}}
\newcommand{\redtr}{\operatorname{\tt Tr}}
\newcommand{\redn}{\operatorname{\tt N}}

\newcommand{\stab}{\operatorname{Stab}}

\newcommand\normalout{\partial^1_{+}}
\newcommand\normalin{\partial^1_{-}}
\newcommand\normalpm{\partial^1_{\pm}}

\newcounter{fig}

\title{On the arithmetic of crossratios \\ and 
generalised Mertens' formulas}
\author{Jouni Parkkonen \and Fr\'ed\'eric Paulin} 

\begin{document}
\bibliographystyle{../alphanum}
\maketitle

\begin{abstract}
  We develop the relation between hyperbolic geometry and arithmetic
  equidistribution problems that arises from the action of arithmetic
  groups on real hyperbolic spaces, especially in dimension $\leq 5$.
  We prove generalisations of Mertens' formula for quadratic imaginary
  number fields and definite quaternion algebras over $\QQ$, counting
  results of quadratic irrationals with respect to two different
  natural complexities, and counting results of representations of
  (algebraic) integers by binary quadratic, Hermitian and Hamiltonian
  forms with error bounds. For each such statement, we prove an
  equidistribution result of the corresponding arithmetically defined
  points.  Furthermore, we study the asymptotic properties of
  crossratios of such points, and expand Pollicott's recent results on
  the Schottky-Klein prime functions.
\end{abstract}

\section{Introduction}

The aim of this paper is to prove various equidistributions results
(and their related asymptotic counting results) of arithmetically
defined points in tori (quotient of the real or complex line $\RR$ or
$\CC$, or of Hamilton's quaternion algebra $\HH$, by arithmetic
lattices), organised using appropriate complexities (see for instance,
amongst the numerous litterature, \cite{Oh10, Harcos10, Serre12,
  BenQui12} and their references for other types of results).

We denote by $\Delta_x$ the unit Dirac mass at a point $x$, by
$\weakstar$ the weak-star convergence of measures, and by $\Leb_E$ the
standard Lebesgue measure on a Euclidean space $E$. Given an imaginary
quadratic number field $K$, we will denote by $\OOO_K$ its ring of
integers, by $D_K$ its discriminant, by $\zeta_K$ its zeta function
and by $\n$ its norm.

\medskip
The first two statements are equidistribution results of rational
points (satisfying congruence properties) in the complex field and
Hamilton's quaternion algebra, analogous to the equidistribution
result of Farey fractions in the real field, where the complexity is
the norm of the denominator.

\btheo\label{theo:appliHermintro} Let $\mmm$ be a (nonzero) fractional
ideal of the ring of integers of an imaginary quadratic number field
$K$, with norm $\n(\mmm)$. As $s\to+\infty$,
$$
\frac{|D_K|\,\zeta_K(2)}{2\,\pi\,  s^2}
\sum_{\substack{(u,v)\in\mmm\times\mmm\;\\
\;\n(\mmm)^{-1}\n(v)\leq s,
\ \OOO_K u+\OOO_K v=\mmm}} \Delta_{\frac uv}\;\weakstar\;\Leb_\CC\,.
$$
\etheo

Let $\HH$ be Hamilton's quaternion algebra over $\RR$, with reduced
norm $\redn: x\mapsto x\overline{x}$. Let $A$ be a quaternion algebra
over $\QQ$, which is definite ($A\otimes_\QQ\RR=\HH$), with reduced
discriminant $D_A$. Let $\OOO$ be a maximal order in $A$, and let
$\mmm$ be a (nonzero) left ideal of $\OOO$, with reduced norm
$\redn(\mmm)$ (see \cite{Vigneras80} for definitions). We refer to
Equation \eqref{eq:equidistribtwoideals} for a more general result.

\btheo\label{theo:appliHamintro} As $s$ tends to $+\infty$, we have
$$
\frac{\pi^2\,\zeta(3)\,\prod_{p|D_A}(p^3-1)}{360\,D_A\; s^4}
\sum_{\substack{(u,v)\in\mmm\times\mmm\;\\
\;\redn(\mmm)^{-1}\redn(v)\leq s,
\ \OOO u+\OOO v=\mmm}}\Delta_{uv^{-1}}\;\weakstar\;\Leb_\HH\,.
$$
\etheo

The next three statements are equidistribution results on $\RR$ or
$\CC$ of points arithmetically constructed by using quadratic
irrationals.

We first prove that the set of traces $\tr \alpha$ of the real
quadratic irrationals $\alpha$ over $\QQ$ at distance at least $s>0$
from their Galois conjugate $\alpha^\sigma$, in a given orbit $G \cdot
\alpha_0$ by homographies of a finite index subgroup $G$ of the
modular group $\PSLZ$, equidistributes to the Lebesgue measure on
$\RR$, as $s\ra 0$, using as complexity (the inverse of) the distance
to their Galois conjugate (see \cite{ParPau11MZ,ParPau12JMD} for
background on this height). Given a real integral quadratic irrational
$\alpha$ over $\QQ$, we denote by $R_\alpha$ the regulator of
$\ZZ+\alpha\ZZ$, and by $Q_\alpha(t)= t^2-(\tr \alpha)\,t+\n(\alpha)$
its associated monic quadratic polynomial (the integrality assumption
is only used in this introduction, to simplify the notation).  We
denote by $H_{x}$ the stabiliser of $x \in \PP_{1}(\RR)=
\RR\cup\{\infty\}$ in a subgroup $H$ of $\PSL_2(\RR)$.

\btheo\label{theo:equidtraceintro} Let $\alpha_0$ be a real integral
quadratic irrational over $\QQ$ and let $G$ be a finite index subgroup
of $\PSLZ$.  Then, as $\epsilon\ra 0$,
$$
\frac{\pi^2\;[\PSLZ:G]\;\epsilon}
{6\;[\PSLZ_{\alpha_0}:G_{\alpha_0}]\;R_{\alpha_0} }\;\;
\sum_{\alpha\in G \cdot \alpha_0
\;:\; |\alpha-\alpha^\sigma|\geq \epsilon}\Delta_{\tr\, \alpha}\;\;
\stackrel{*}{\rightharpoonup}\;\; \Leb_\RR\,.
$$
\etheo 

We refer to Theorems \ref{theo:appliarithdim3} and
\ref{theo:appliarithdim5} for extensions of the above result to
quadratic irrationals over an imaginary quadratic number field (using
relative traces) or over a rational quaternion algebra, and we only
quote in this introduction the following special case of Theorem
\ref{theo:appliarithdim3}.

\bcoro \label{coro:orbianchiintro} Let $\phi=\frac{1+\sqrt{5}}2$ be
the Golden Ratio, let $K$ be an imaginary quadratic number field with
$D_K\neq -4$, let $\ccc$ be a nonzero ideal in $\OOO_K$, and let
$\Ga_{0}(\ccc)$ be the Hecke congruence subgroup $\Big\{ \pm
\Big(\!\begin{array}{cc}a & b \\ c & d \end{array} \!\Big) \in \PSLOK
\;:\; c\in\ccc\Big\}$. Then, as $\epsilon\ra 0$,
$$
\frac{|D_K|^{\frac{3}{2}}\; \zeta_{K}(2)\,
  \n(\ccc)\prod_{\ppp |\ccc} \big(1+\frac{1}{\n(\ppp)}\big)\;\epsilon^2}
{4\,\pi^2\;k_\ccc\;\log\phi} \; 
\sum_{\alpha\in \Ga_{0}(\ccc) \cdot \phi,\,
|\alpha-\alpha^\sigma|\geq \epsilon}\Delta_{\tr\,\alpha}\;
\weakstar\;\Leb_\CC\,.
$$
where $k_\ccc$ is the smallest $k\in\NN-\{0\}$ such that the
$2k$-th term of Fibonacci's sequence belongs to  $\ccc$.
\ecoro

\medskip Given two real quadratic irrationals $\alpha,\beta$ over
$\QQ$, we introduce the {\it relative heigth} $h_\alpha(\beta)$ of
$\beta$ with respect to $\alpha$, measuring how close the (unordered)
pair $\{\beta, \beta^\sigma\}$ is to the pair
$\{\alpha,\alpha^\sigma\}$, by
$$
h_{\alpha}(\beta)=
\frac{ \min\{|\beta-\alpha|\,|\beta^\sigma-\alpha^\sigma|,\;
|\beta-\alpha^\sigma|\,|\beta^\sigma-\alpha|\}}{|\beta-\beta^\sigma|}\,.
$$ See Equation \eqref{eq:relheiwithcrossratio} for an expression of
$h_{\alpha}(\beta)$ using crossratios of $\alpha,\beta,\alpha^\sigma,
\beta^\sigma$. Consider the points
\begin{equation}\label{eq:defixpmalphbet}
x^\pm_{\alpha}(\beta)=\frac{\n(\beta)-\n(\alpha)\pm
\big((\n(\alpha)-\n(\beta))^2+ (\tr\,\beta-\tr\,\alpha)
(\tr\,\beta\n(\alpha)-\tr\,\alpha\n(\beta))\big)^{\frac{1}{2}}}
{\tr\,\beta-\tr\,\alpha}\,.
\end{equation}
We will prove that when $\beta$ varies in an orbit of the modular
group $\PSLZ$, the relative height $h_{\alpha}(\beta)$ is a well
defined complexity modulo the stabiliser of $\alpha$, and, except
finitely many orbits under this stabiliser, the points
$x^\pm_{\alpha}(\beta)$ are well defined and real. The following
results says that these points, when $\beta$ has relative height at
most $s$ tending to $\infty$, equidistribute to the (unique up to
scalar) measure on $\RR-\{\alpha_0,\alpha_0^\sigma\}$ which is
absolutely continuous with respect to the Lebesgue measure and
invariant under the stabiliser of $\{\alpha_0,\alpha_0^\sigma\}$ under
$\PSL_2(\RR)$.

\btheo\label{theo:relcomplcountintro} Let $\alpha,\beta$ be real
integral quadratic irrationals and let $G$ be a finite index subgroup
of $\PSLZ$. Then, as $s\ra +\infty$,
$$
\frac{\pi^2\,[\PSLZ:G]}
{24\,[\PSLZ_{\beta}:G_{\beta}]\,R_{\beta}\,s}\;
\sum_{\beta'\in G\cdot\beta,\;h_{\alpha}(\beta')\le s} 
\Delta_{x^-_{\alpha}(\beta')}+\Delta_{x^+_{\alpha}(\beta')}
\;\weakstar\;\frac{dt}{|Q_\alpha(t)|}\,.
$$
\etheo

This result, whose proof uses asymptotic properties of crossratios,
implies that 
$$
\card \{\beta'\in \operatorname{SO}_{Q_\alpha}(\ZZ)\bs
\PSLZ\cdot\beta,\; h_{\alpha}(\beta')\le s\} \sim
\frac{48\,R_{\alpha}\,R_{\beta}}
{\pi^2\,|\alpha-\alpha^\sigma|}\; s\,.
$$
We refer to Theorem \ref{theo:relcomplcount} for a version with
congruences, and to Theorem \ref{theo:relcomplcountBianchi} for an
extension of this result to quadratic irrationals over an imaginary
quadratic extension of $\QQ$.

\medskip

Towards the same measure, we also have the following equidistribution
result of integral representations of quadratic norm forms (see
Section \ref{subsec:normform} for generalisations).

\btheo\label{theo:normformintro} If $\alpha$ is a real quadratic
irrational over $\QQ$, then, as $s\ra+\infty$,
$$
\frac{\pi^2}{12\;s}\;
\sum_{(u,v)\in \ZZ^2,\;(u,v)=1,\;|u^2-\tr \alpha \,uv +\n(\alpha)\,v^2|\leq s} 
\Delta_{\frac uv}
\;\weakstar\;\frac{dt}{|Q_{\alpha}(t)|}\,,
$$
\etheo

\medskip Our final equidistribution result for this introduction is
the following equidistribution of coefficients of binary Hermitian
forms in an orbit under the Bianchi group $\SLOK$, using as complexity
their first coefficient (see Subsection
\ref{subsect:errortermintrepbybinary} for extensions, in particular to
binary Hamiltonian forms). Given an imaginary quadratic number field
$K$, let $f(u,v)=a\,|u|^2+2\,\Re(b\,\ov u\,v) +c\,|v|^2$ be a binary
Hermitian form, which is integral over $\OOO_K$ (that is
$a=a(f),c=c(f)\in\ZZ$ and $b=b(f)\in\OOO_K$), and indefinite (that is
$\discr(f)=|b|^2-ac>0$). We denote by $\cdot$ the action of $\SLC$ on
the set of binary Hermitian forms by precomposition, and by $\autom$
the stabiliser of $f$ in $\SLOK$.

\btheo \label{theo:mainHermerrorintro} Let $G$ be a finite index
subgroup of $\SLOK$. As $s\ra+\infty$,
$$
\frac{[\SLOK:G]\,|D_K|^{\frac{3}2}\,\zeta_K(2)\,\discr(f)}
{\iota_G\,\pi^2\,\covol(\autom\cap G)\,s^2}\; \sum_{f'\in
  \,G \cdot f,\; \,0\,<\,|a(f')|\,\le\, s} \Delta_{\frac{b(f')}{a(f')}}
\;\weakstar\;\Leb_\CC\,.
$$
\etheo

Identifying $\OOO_K$ with the upper triangular unipotent sugbroup of
$\SLOK$, this theorem implies the following counting result of
integral binary Hermitian forms (ordered by their first coefficient)
in an orbit of the Bianchi group $\SLOK$: as $s\ra+\infty$,
$$
\card\{f'\in\, \OOO_K\bs \SLOK\cdot f\;:\;0<|a(f')|\leq s\}
\sim \frac{\pi^2\,\iota_G\,\covol(\autom)}
{2\,|D_K|^{\frac{1}2}\,\zeta_K(2)\,\discr(f)}
\;\;
s^2\;.
$$

All these theorems come with an error term, and will be (not
immediate) consequences of the geometric equidistribution results of
\cite{ParPau14}. The particular case of the result of loc.~cit.~that
we will use, stated in Section \ref{sect:geometry}, is the following
one (see also \cite{OhSha13, OhShaCircles} for related counting and
equidistribution results in real hyperbolic spaces and \cite{Kim13,
  ParPauHeis} in negatively curved symmetric spaces) : given a lattice
$\Ga$ in the isometry group of the real $n$-dimensional hyperbolic
space $\HH^n_\RR$, and two horoballs or totally geodesic subspaces
$D^-,D^+$ whose stabilisers in $\Ga$ have cofinite volume, the initial
tangent vectors of the common perpendiculars between $D_-$ and the
images under $\Ga$ of $D^+$ equidistribute in the unit normal bundle
of $D^-$.

The paper is organised according to the arithmetic applications: In
Section \ref{sect:mertens}, we apply Section \ref{sect:geometry} with
both $D^-$ and $D^+$ horoballs to prove generalisations of the
classical Mertens formula, describing the asymptotic behaviour of the
average order of Euler's function, for the rings of integers of
quadratic imaginary number fields and maximal orders in definite
quaternion algebras over $\QQ$.  In Section \ref{sect:irrationals}, we
consider counting and equidistribution of quadratic irrationals in
terms of two complexities, the height introduced in \cite{ParPau11MZ}
and the relative height mentionned above. In Subsection
\ref{subsect:errortermquadirrat}, the geometric result is applied
when $D^-$ is a horoball and $D^+$ is a geodesic line, and in
Subsection \ref{subsec:relcompquadirrat} and
\ref{sect:countingcrossratios} (where we extend Pollicott's result on
the asymptotic of crossratios to prove the convergence of generalised
Schottky-Klein functions) when both $D^-$ and $D^+$ are geodesic
lines.  In the final section, we consider representations of integers
by binary quadratic, Hermitian and Hamiltonian forms, applying the
geometric result when $D^-$ is a horoball and $D^+$ is either a
geodesic line or a totally geodesic hyperplane or (when considering
positive definite forms) a point.

\medskip\noindent{\small {\it Acknowledgement: } The first author
  thanks the Université Paris-Sud (Orsay) for a visit of a month
  and a half which allowed most of the writing of this paper, under
  the financial support of the ERC grant GADA 208091.}

\section{Geometric counting and equidistribution}
\label{sect:geometry}

In this section, we briefly review a simplified version of the
geometric counting and equidistribution results proved in
\cite{ParPau14}, whose arithmetic applications will be considered in
the other parts of this paper (see also \cite{ParPauRev} for related
references and \cite{ParPauHeis} for the case of locally symmetric
spaces).

Let $n\geq 2$, let $\Ga$ be a discrete nonelementary group of
isometries of the $n$-dimensional real hyperbolic space $\hnr$, and
let $M=\Ga\bs\hnr$ be the quotient orbifold. Let $D^-$ and $D^+$ be
nonempty proper closed convex subsets in $\hnr$, such that the
families $(\ga D^-)_{\ga\in\Ga}$ and $(\ga D^+)_{\ga\in\Ga}$ are
locally finite in $\hnr$ (see \cite[\S 3.3]{ParPau14} for more general
families in any simply connected complete Riemannian manifold with
pinched negative curvature).

\medskip We denote by $\partial_\infty \hnr$ the boundary at infinity
of $\hnr$, by $\Lambda \Ga$ the limit set of $\Ga$ and by
$(\xi,x,y)\mapsto \beta_\xi(x,y)$ the Busemann cocycle on
$\partial_\infty \hnr\times \hnr\times \hnr$ defined by
$$
(x,y,\xi)\mapsto \beta_{\xi}(x,y)=
\lim_{t\to+\infty}d(\rho_t,x)-d(\rho_t,y)\;,
$$
where $\rho:t\mapsto \rho_t$ is any geodesic ray with point at
infinity $\xi$ and $d$ is the hyperbolic distance.  
%%See for instance \cite[\S II.8]{BriHae99}. 
%In the upper halfspace model, the boundary at infinity of the real
%hyperbolic space $\hnr$ is $\RR^{n-1}\cup\{\infty\}$, and the 
%Busemann cocycle of $\hnr$ is
%$$
%\beta_\infty(x,y)=\log\frac{y_n}{x_n}
%$$
%for all $x,y\in\hnr$ and
%$$
%\beta_\xi(x,y)=
%\ln(\frac {y_{n}}{x_{n}}\frac{\|x-\xi\|^2}{\|y-\xi\|^2})
%$$
%for all $x,y\in\hnr$ and all $\xi\in\RR^{n-1}$, where $\|z\|$ is 
%the Euclidean norm and $z_n$ the last coordinate of $z\in \RR^{n}$.

For every $v\in T^1\hnr$, let $\pi(v)\in \hnr$ be its origin, and let
$v_-, v_+$ be the points at infinity of the geodesic line in $\hnr$
whose tangent vector at time $t=0$ is $v$. We denote by $\normalpm
D^\mp$ the {\it outer/inner unit normal bundle} of $\partial D^\mp$,
that is, the set of $v\in T^1\hnr$ such that $\pi(v)\in \partial
D^\mp$ and the closest point projection on $D^\mp$ of $v_\pm
\in \partial_\infty \wt M-\partial_\infty D^\mp$ is $\pi(x)$. For
every $\ga,\ga'$ in $\Ga$ such that $\ga D^-$ and $\ga' D^+$ have a
common perpendicular (that is, if the closures $\overline{\ga D^-}$
and $\overline{\ga' D^+}$ in $\hnr\cup\partial_\infty \hnr$ are
disjoint), we denote by $\alpha_{\ga,\,\ga'}$ this common
perpendicular (starting from $\ga D^-$ at time $t=0$), by
$\ell(\alpha_{\ga,\,\ga'})$ its length, by $v^-_{\ga,\,\ga'} \in \ga
\normalout D^-$ its initial tangent vector and by $v^+_{\ga,\,\ga'}
\in \ga'\normalin D^+$ its terminal tangent vector.  Let $\Ga_{D^-}$
and $\Ga_{D^+}$ be the stabilisers in $\Ga$ of the subsets $D^-$ and
$D^+$, respectively.  The {\em multiplicity} of $\alpha_{\ga,\ga'}$ is
$$
m_{\ga,\ga'}=\frac 1{\card(\ga\Ga_{D^-}\ga^{-1}\cap\ga'\Ga_{D^+}{\ga'}^{-1})}\,,
$$
which equals $1$ when $\Ga$ acts freely on $T^1\hnr$ (for instance
when $\Ga$ is torsion-free). Let
$$
\N_{D^-,\,D^+}(t)=
\sum_{\substack{
(\ga,\,\ga')\in \,\Ga\bs((\Ga/\Ga_{D^-})\times (\Ga/\Ga_{D^+}))\\
\phantom{\big|}\overline{\ga D^-}\,\cap \,\overline{\ga' D^+}\,=\emptyset,\; 
\ell(\alpha_{\ga,\, \ga'})\leq t}} m_{\ga,\ga'}=
\sum_{\substack{
[\ga]\in\, \Ga_{D^-}\bs\Ga/\Ga_{D^+}\\
\phantom{\big|}\overline{D^-}\,\cap \,\overline{\ga D^+}\,=\emptyset,\; 
\ell(\alpha_{e,\, \ga})\leq t}} m_{e,\ga}
\;,
$$
where $\Ga$ acts diagonally on $\Ga\times\Ga$. When $\Ga$ has no
torsion, $\N_{D^-,\,D^+}(t)$ is the number (with multiplicities coming
from the fact that $\Ga_{D^\pm}\bs D^\pm$ is not assumed to be
embedded in $M$) of the common perpendiculars of length at most $t$
between the images of $D^-$ and $D^+$ in $M$. We refer to \cite[\S
4]{ParPau14} for the use of H\"older potentials on $T^1\hnr$ to modify
this counting function by adding weights, 
%\todo{We should remember to check this idea!}  
which could be useful for some further arithmetic applications.

Recall the following notions (see for instance \cite{Roblin03}). The
{\em critical exponent} of $\Ga$ is
$$
\delta_{\Ga}=\limsup_{N\to+\infty}\frac 1N \ln\card\{\ga\in\Ga:
d(x_{0},\ga x_{0})\leq N\}\,,
$$
which is positive, finite and independent of the base point $x_{0}\in
\hnr$. Let $(\mu_{x})_{x\in \hnr}$ be a {\em Patterson density} for
$\Ga$, that is, a family $(\mu_{x})_{x\in \hnr}$ of nonzero finite
measures on $\partial_{\infty}\hnr$ whose support is $\Lambda\Ga$,
such that $\ga_*\mu_x=\mu_{\ga x}$ and
$$
\frac{d\mu_{x}}{d\mu_{y}}(\xi)=e^{-\delta_{\Ga}\beta_{\xi}(x,y)}
$$
for all $\ga\in\Ga$, $x,y\in\hnr$ and $\xi\in\partial_{\infty}\hnr$.
The
{\em Bowen-Margulis measure} $\wt m_{\rm BM}$ for $\Ga$ on $T^1\hnr$
is defined, using Hopf's parametrisation $v\mapsto
(v_-,v_+,\beta_{v_+}(x_0,\pi(v))\,)$ of $T^1\hnr$, by
$$
d\wt m_{\rm BM}(v)=e^{-\delta_{\Ga}(\beta_{v_{-}}(\pi(v),\,x_{0})+
\beta_{v_{+}}(\pi(v),\,x_{0}))}\;
d\mu_{x_{0}}(v_{-})\,d\mu_{x_{0}}(v_{+})\,dt\,.  
$$
The measure $\wt m_{\rm BM}$ is nonzero and independent of
$x_{0}\in\hnr$.  It is invariant under the geodesic flow, the
antipodal map $v\mapsto -v$ and the action of $\Ga$, and thus defines
a nonzero measure $m_{\rm BM}$ on $T^1M$, called the {\em
  Bowen-Margulis measure} on $M=\Ga\bs \hnr$, which is invariant under
the geodesic flow of $M$ and the antipodal map.  When $m_{\rm BM}$ is
finite, denoting the total mass of a measure $m$ by $\|m\|$, the
probability measure $\frac{m_{\rm BM}} {\|m_{\rm BM}\|}$ is then
uniquely defined, and it is the unique probability measure of maximal
entropy for the geodesic flow (see \cite{OtaPei04}). This holds for
instance when $M$ has finite volume or when $\Ga$ is geometrically
finite.

\medskip Using the endpoint homeomorphism $v\mapsto v_+$ from
$\normalout {D^-}$ to $\partial_{\infty}\hnr-\partial_{\infty}D^-$,
the {\em skinning measure}
$\wt\sigma_{D^-}$ of $\Ga$ on $\normalout{D^-}$ is defined by
$$
d\wt\sigma_{D^-}(v) =  
e^{-\delta\,\beta_{v_{+}}(\pi(v),\,x_{0})}\,d\mu_{x_{0}}(v_{+})\,,
$$
see \cite[\S 1.2]{OhSha13} when $D^-$ is a horoball or a totally
geodesic subspace in $\HH^n_\RR$ and \cite{ParPau13ETDS},
\cite{ParPau14} for the general case of convex subsets in variable
curvature and with a potential.

The measure $\wt\sigma_{D^-}$ is independent of $x_{0}\in\hnr$, it is
nonzero if $\Lambda\Ga$ is not contained in $\partial_{\infty}D^-$,
and satisfies $\wt\sigma_{\ga D^-} =\ga_{*}\wt\sigma_{D^-}$ for every
$\ga\in\Ga$. Since the family $(\ga D^-)_{\ga\in\Ga}$ is locally
finite in $\hnr$, the measure $\sum_{\ga \in \Ga/\Ga_{D^-}}
\;\ga_*\wt\sigma_{D^-}$ is a well defined $\Ga$-invariant locally
finite (Borel nonnegative) measure on $T^1\hnr$, hence induces a
locally finite measure $\sigma_{D^-}$ on $T^1M=\Ga\bs T^1\hnr$, called
the {\em skinning measure} of $D^-$ in $T^1M$. We refer to \cite[\S
5]{OhShaCircles} and \cite[Theo.~9]{ParPau13ETDS} for finiteness
criteria of the skinning measure $\sigma_{D^-}$, in particular
satisfied when $M$ has finite volume and if either $D^-$ is a horoball
centred at a parabolic fixed point of $\Ga$ or if $D^-$ is a totally
geodesic subspace.

\medskip The following result on the asymptotic behaviour of the
counting function $\N_{D^-,\,D^+}$ in real hyperbolic space is a
special case of much more general results \cite[Coro.~20, 21,
Theo.~28]{ParPau14}. Furthermore, the equidistribution result of the
initial and terminal tangent vectors of the common perpendiculars
holds simultaneously in the outer and inner tangent bundles of $D^-$
and $D^+$. We refer to \cite{ParPauRev} for a survey of the particular
cases known before \cite{ParPau14} due to Huber, Margulis, Herrmann,
Cosentino, Roblin, Oh-Shah, Martin-McKee-Wambach, Pollicott, and the
authors for instance.

For every $t\geq 0$, let 
$$
m_t(x)=\sum_{\ga\in \Ga/\Ga_{D^+}\;:\;
  \overline{D^-}\,\cap \,\overline{\ga D^+}\,= \emptyset,\;
\alpha_{e,\, \ga}(0)=x,\; \ell(\alpha_{e,\, \ga})\leq t} m_{e,\ga}
$$ 
be the multiplicity of a point $x \in \partial D^-$ as the origin of
common perpendiculars with length at most $t$ from $D^-$ to the
elements of the $\Ga$-orbit of $D^+$.  We denote by $\Delta_x$ the
unit Dirac mass at a point $x$.

\btheo\label{theo:mainequicountdown} Let $\Ga,D^-,D^+$ be as
above. Assume that the measures $m_{\rm BM},\sigma_{D^-},\sigma_{D^+}$
are nonzero and finite. Then
$$
\N_{D^-,\,D^+}(t)\;\sim\;
\frac{\|\sigma_{D^-}\|\;\|\sigma_{D^+}\|}{\delta_\Ga\;\|m_{\rm
    BM}\|}\; e^{\delta_\Ga \,t}\;,
$$
as $t\ra+\infty$. If $\Ga$ is arithmetic or if $M$ is compact, then
the error term is $\bigO(e^{(\delta_\Ga-\kappa) t})$ for some $\kappa
>0$.  Furthermore, the origins of the common perpendiculars
equidistribute in the boundary of $D^-$: 
\begin{equation}\label{eq:equidistribdown}
\lim_{t\ra+\infty}\; 
\frac{\delta_\Ga\;\|m_{\rm BM}\|}{\|\sigma_{D^-}\|\;\|\sigma_{D^+}\|}
\;e^{-\delta_\Ga\, t}\;\sum_{x\in\partial D^-} m_{t}(x) \;
\Delta_{x}\;=\; \frac{\pi_*\wt\sigma_{D^-}}{\|\sigma_{D^-}\|}\;
\end{equation}
for the weak-star convergence of measures on the locally compact space
$T^1\hnr$. 
\cqfd
\etheo

For smooth functions $\psi$ with compact support on $\partial D^-$,
there is an error term in the above equidistribution claim, of the form
$\bigO(s^{\delta_\Ga-\kappa}\|\psi\|_\ell)$ where $\kappa>0$ and
$\|\psi\|_\ell$ is the Sobolev norm of $\psi$ for some $\ell\in\NN$,
as proved in \cite[Theo.~28]{ParPau14}.

\medskip When $M$ has finite volume, we have $\delta_\Ga=n-1$, the
Bowen-Margulis measure $m_{\rm BM}$ coincides up to a multiplicative
constant with the Liouville measure on $T^1M$, and the skinning
measures of points, horoballs and totally geodesic subspaces $D^\pm$
coincide again up to a multiplicative constant with the (homogeneous)
Riemannian measures on $\normalpm D^\mp$ induced by the Riemannian
metric of $T^1\hnr$.  These proportionality constants were computed in
\cite[\S 7]{ParPauRev} and \cite[Prop.~29]{ParPau14}:
$$
\|m_{\rm BM}\|=2^{n-1}\,(n-1)\,\Vol(\SSS^{n-1})\,\Vol(M)\,,
$$
if $D^-$ is a horoball, then $\|\sigma_{D^-}\| =2^{n-1}(n-1)
\Vol(\Ga_{D^-}\bs D^-)$, and if $D^-$ is a totally geodesic
submanifold in $\hnr$ of dimension $k^-\in\{0,\dots,n-1\}$ with
pointwise stabiliser of order $m^-$, then $\wt\sigma_{D^-}=
\Vol_{\normalout D^-}$, so that $\|\sigma_{D^-}\|
=\frac{\Vol(\SSS^{n-k^--1})}{m^-}\Vol(\Ga_{D^-}\bs D^-)$.  See
\cite[\S 3]{ParPauHeis} for the computation of the proportionality
constants in the complex hyperbolic case.

Using these explicit expressions, we now reformulate Theorem
\ref{theo:mainequicountdown} using these explicit expressions,
considering the following cases.

\smallskip
 (1) If $D^-$ and $D^+$ are totally geodesic submanifolds in
$\hnr$ of dimensions $k^-$ and $k^+$ in $\{1,\dots,n-1\}$,
respectively, such that $\Vol(\Ga_{D^-}\bs D^-)$ and
$\Vol(\Ga_{D^+}\bs D^+)$ are finite, let
$$
c(D^-,\,D^+)=\frac{\Vol(\SSS^{n-k^--1})\Vol(\SSS^{n-k^+-1})}
{2^{n-1}(n-1)\Vol(\SSS^{n-1})} 
\frac{\Vol(\Ga_{D^-}\bs D^-)\Vol(\Ga_{D^+}\bs D^+)}{\Vol(M)}\,.
$$ 

\smallskip (2) If $D^-$ and $D^+$ are horoballs in $\hnr$ centred at
parabolic fixed points of $\Ga$, let
$$
c(D^-,D^+)=\frac{2^{n-1}(n-1)} {\Vol(\SSS^{n-1})}
\frac{\Vol(\Ga_{D^-}\bs D^-)\Vol(\Ga_{D^+}\bs D^+)} {\Vol(M)}\,.
$$

\smallskip (3) If $D^-$ is a horoball in $\hnr$ centred at a parabolic
fixed point of $\Ga$ and $D^+$ is a totally geodesic subspace in
$\hnr$ of dimension $k^+\in \{1,\dots,n-1\}$ such that
$\Vol(\Ga_{D^+}\bs D^+)$ is finite, let
$$
c(D^-,D^+)=\frac{\Vol(\SSS^{n-k^+-1})}{\Vol(\SSS^{n-1})}
\frac{\Vol(\Ga_{D^-}\bs D^-)\Vol(\Ga_{D^+}\bs D^+)}
{\Vol(M)}\,.
$$

\smallskip (4) If $D^-$ is a horoball in $\hnr$ centred at a parabolic
fixed point of $\Ga$ and $D^+$ is a point of $\hnr$, let
$$
c(D^-,D^+)=
\frac{\Vol(\Ga_{D^-}\bs D^-)}
{\Vol(M)}\,.
$$

\bcoro\label{coro:constcurvcount} Let $\Ga$ be a discrete group of
isometries of $\hnr$ such that the orbifold $M=\Ga\bs\hnr$ has finite
volume. In each of the cases (1) to (4) above, if $m^\pm$ is the
cardinality of the pointwise fixator of $D^\pm$, then
$$
\N_{D^-,\,D^+}(t)\sim \frac{c(D^-,D^+)}{m^-\,m^+}\;e^{(n-1)t}\,.
$$
If $\Ga$ is arithmetic or if $M$ is compact, then there exists 
$\kappa>0$ such that, as $t\ra+\infty$,
$$
\N_{D^-,\,D^+}(t)=\frac{c(D^-,D^+)}{m^-\,m^+}\;
e^{(n-1)t}\big(1+\operatorname{O}(e^{-\kappa t})\big)\;.
$$
Furthermore, the origins of the common perpendiculars equidistribute
in $\partial D^-$ to the induced Riemannian measure: 
if $D^-$ is a horoball centred at a parabolic fixed point of $\Ga$, then
\begin{equation}\label{eq:distribhorob}
\frac{m^+\,(n-1)\Vol(\Ga_{D^-}\bs D^-)}{c(D^-,D^+)}
\;e^{-(n-1)t}\;\sum_{x\in\partial D^-} m_{t}(x) \;
\Delta_{x}\;\weakstar\; \Vol_{\partial D^-}\,,
\end{equation}
as $t\to+\infty$, 
and if $D^-$ is a totally geodesic subspace in $\hnr$ of dimension
$k^-\in \{1,\dots,n-1\}$, such that $\Vol(\Ga_{D^-}\bs D^-)$ is
finite, then as $t\to+\infty$, 
\begin{equation}\label{eq:distribtotgeod}
\frac{m^-\,m^+\,\Vol(\Ga_{D^-}\bs D^-)}{c(D^-,D^+)}
\;e^{-(n-1) t}\;\sum_{x\in D^-} m_{t}(x) \;
\Delta_{x}\;\weakstar\; \Vol_{D^-}\,.
\;\;\;\Box
\end{equation}
\ecoro

Again, for smooth functions $\psi$ with compact support on $\partial
D^-$, there is an error term in the above two equidistribution claims, of
the form $\bigO(s^{n-1-\kappa}\|\psi\|_\ell)$ where $\kappa>0$ and
$\|\psi\|_\ell$ is the Sobolev norm of $\psi$ for some $\ell\in\NN$.

\medskip 
We end this section by recalling some terminology concerning the
isometries of $\hnr$.  The {\em translation length} of an isometry
$\ga$ of $\hnr$ is
$$
\ell(\ga)=\inf_{x\in\hnr}d(x,\ga x)\,.
$$
If $\ell(\ga)>0$, then $\ga$ is {\em loxodromic}.  Each loxodromic
isometry $\ga$ stabilises a unique geodesic line in $\hnr$, called the
{\em translation axis} of $\ga$ and denoted by $\axis\ga$, on which it
acts as a translation by $\ell(\ga)$.  The points at infinity of
$\axis\ga$ are the two fixed points of $\ga$ in $\partial_\infty\hnr$,
and we denote by $\ga^-$ and $\ga^+$ the attracting and repelling
fixed points of $\ga$ respectively.  Any loxodromic element of $\Ga$
is contained in a maximal cyclic subgroup.

An element $\ga\in\Ga$ is {\em primitive} if the cyclic subgroup
$\ga^\ZZ$ it generates is a maximal cyclic subgroup of $\Ga$.  A
loxodromic element $\ga\in\Ga$ is $\Ga${\em-reciprocal} if there is an
element in $\Ga$ that switches the two fixed points of $\ga$.  If
$\ga$ is $\Ga$-reciprocal, then let $\iota_\Ga(\ga)=2$, otherwise, we set
$\iota_\Ga(\ga)=1$.  If $\ga$ is a primitive loxodromic element of
$\Ga$, then the stabiliser of $\axis \ga$ is generated by $\ga$, an
elliptic element that switches the two points at infinity of the axis
of $\ga$ if $\ga$ is reciprocal, and a (possibly trivial) group of
finite order $m_\Ga(\ga)$, which is the pointwise stabiliser of
$\axis\ga$, so that 
\begin{equation}\label{eq:miotacardstab}
m_\Ga(\ga)=\frac{1}{\iota_\Ga(\ga)}\,
[\stab_\Ga (\axis \ga):\ga^\ZZ]\,.
\end{equation}

\section{Generalised Mertens' formulas}
\label{sect:mertens}

A classical result, known as {\it Mertens' formula} (see for example
\cite[Thm.~330]{HarWri08}, a better error term is due to Walfisz
\cite{Walfisz63}), describes the asymptotic behaviour of the average
order $\Phi:\NN-\{0\}\ra \NN-\{0\}$ of Euler's function $\varphi$,
where $\varphi(n)=\card((\ZZ/n\ZZ)^\times)$~:
$$
\Phi(n)=\sum_{k=1}^n\varphi(k)=
\frac 3{\pi^2}\,n^2+\operatorname{O}(n\log n)\,.
$$

Corollary \ref{coro:constcurvcount} provides a geometric proof of
Mertens' formula with a less explicit error term, as follows.  The
modular group $\Ga_\QQ= \PSLZ$ isometrically acts on the upper half
plane model of the real hyperbolic plane $\hdr$ via M\"obius
transformations. Let $\H_\infty$ be the horoball in $\hdr$ that
consists of all points with vertical coordinate at least $1$, and let
$\Ga_\infty$ be its stabiliser in $\Ga_\QQ$.  The images of
$\H_\infty$ under $\Ga_\QQ$ different from $\H_\infty$ are the disks
in $\hdr$ tangent to the real line at the rational points $\frac pq$
(where $q>0$ and $(p,q)=1$) with Euclidean diameter $\frac 1{q^2}$.
The common perpendicular from $\H_\infty$ to this disc exists if and
only if $q>1$, its length is $\log q^2$, and its multiplicity is
$1$. For every $n\in\NN-\{0\}$, the cardinality of the set of rational
points $\frac pq\in \,]0,1]$ with $q\leq n$ is $\Phi(n)$. Thus,
$\Phi(n)= \N_{\H_\infty,\H_\infty}(\log
n^2)+1$. Now$\Vol(\Ga_\infty\bs\H_\infty) =1$ and
$\Vol(\Ga_\QQ\backslash\hdr)=\frac\pi 3$, and we can apply Corollary
\ref{coro:constcurvcount} to conclude that for some $\kappa>0$, as
$n\ra+\infty$,
\begin{equation}\label{eq:Mertens}
\Phi(n)=\frac 3{\pi^2}\,n^2+\operatorname{O}(n^{2-\kappa})\,.
\end{equation}

Let us reformulate Mertens' formula in a way that is easily extendable
in more general contexts (see also \cite{ParPauHeis}).  Let the
additive group $\ZZ$ act on $\ZZ\times\ZZ$ by horizontal shears
(transvections): $k\cdot (u,v)=(u+kv,v)$, and define
$$
\psi(s)=\card\;\ZZ\,\bs
\big\{(u,v)\in\ZZ\times\ZZ\;: (u,v)=1,\  |v|\le s \big\}\,.
$$
We easily have $\psi(n)=2\,\Phi(n)+2$ for every $n\in\NN-\{0\}$, so
that Mertens' formula \eqref{eq:Mertens} is equivalent with
$\psi(n)=\frac 6{\pi^2}\,s+\operatorname{O}(s^{1-\kappa})$.
Furthermore, a straightforward application of Corollary
\ref{coro:constcurvcount} shows that as $s\to+\infty$, we have
$$
\frac{\,\pi^2}{3\,s}\sum_{(u,\,v)=1,\;1\le v\le s}\Delta_{(\frac uv,\,1)}
\weakstar\Vol_{\partial\H_\infty}\,.
$$
Observing that the pushforward of the measures on $\partial \H_\infty$
by the map $(x,1)\mapsto x$ from $\partial\H_\infty$ to $\RR$ is
linear, is continuous for the weak-star topology, maps the unit Dirac
mass at a point $p$ to the unit Dirac mass at $f(p)$ and the
volume measure of $\partial \H_\infty$ to the Lebesgue measure, we get
the well-known result on the equidistribution of the Farey fractions:
as $s\to+\infty$, we have
$$
\frac{\,\pi^2}{6\,s}\;\;\sum_{(u,\,v)=1,\;|v|\le s}\Delta_{\frac uv}\;
\weakstar\;\Leb_\RR\,.
$$

In subsections \ref{subsect:mertensbianchi} and
\ref{subsect:mertenshamilton}, we generalise the above results to
quadratic imaginary number fields and definite quaternion algebras
over $\QQ$.

\subsection{A Mertens' formula for the rings of integers 
of imaginary quadratic number fields}
\label{subsect:mertensbianchi}

Let $K$ be an imaginary quadratic number field, with ring of integers
$\OOO_K$, discriminant $D_K$, zeta function $\zeta_K$ and norm $\n$.
Let $\omega_K=\card({\OOO_K}^\times)$ (with $\omega_K=2$ if $D_K\neq
-3,-4$ for future use).  Let $\mmm$ be a (nonzero) fractional ideal of
$\OOO_K$, with norm $\n(\mmm)$. Note that the action of the additive
group $\OOO_K$ on $\CC\times \CC$ by the horizontal shears
$k\cdot(u,v)=(u+kv,v)$ preserves $\mmm\times\mmm$.

We consider the counting function $\psi_{\mmm}:[0,+\infty[\;\ra \NN$
defined by
$$
\psi_{\mmm}(s)=\card\;\OOO_K\,\bs
\big\{(u,v)\in\mmm\times\mmm\;:\;\n(\mmm)^{-1}\n(v)\leq s,
\ \OOO_K u+\OOO_K v=\mmm\big\}\,.
$$
Note that $\psi_\mmm$ depends only on the ideal class of $\mmm$ and
thus we can assume in the computations that $\mmm$ is integral.  The
{\em Euler function} $\varphi_K$ of $K$ is defined on the set of
(nonzero) integral ideals $\aaa$ of $\OOO_K$ by $\varphi_K(\aaa)=
\card ((\OOO_K/\aaa)^\times)$, see for example \cite[\S 4A]{Cohn78}.
Thus,
$$
\psi_{\OOO_K}(s)= \sum_{v\in\OOO_K,\;0<\n(v)\le s} \varphi_K(v\OOO_K)\,,
$$ 
and the first claim of the following result, in the special case
$\mmm=\OOO_K$, is an analog of Mertens' formula, due to \cite[Satz
2]{Grotz79} (with a better error term), see also \cite[\S
4.3]{Cosentino99}. The equidistribution part of Theorem
\ref{theo:appliHerm} (stated as Theorem \ref{theo:appliHermintro} in
the introduction) generalises \cite[Th.~4]{Cosentino99} which covers
the case $\mmm=\OOO_K$ without an explicit proportionality constant
but with an explicit estimate on the speed of equidistribution.  Note
that the ideal class group of $\OOO_K$ is in general nontrivial, hence
the extension to general $\mmm$ is interesting.

\btheo\label{theo:appliHerm} There exists $\kappa>0$ such that, as
$s\to+\infty$,
$$
\psi_{\mmm}(s)= \frac{\pi}
{\zeta_K(2)\,\sqrt{|D_K|}}\;\;s^2+\operatorname{O}(s^{2-\kappa})\,.
$$
Furthermore, as $s\to+\infty$,
$$
\frac{|D_K|\,\zeta_K(2)}{2\,\pi\,  s^2}
\sum_{\substack{(u,v)\in\mmm\times\mmm\;\\
\;\n(\mmm)^{-1}\n(v)\leq s,
\ \OOO_K u+\OOO_K v=\mmm}}\Delta_{\frac uv}\;\weakstar\;\Leb_\CC\,,
$$
with error term $\bigO(s^{2-\kappa}\|\psi\|_\ell)$
when evaluated on $\C^\ell$-smooth functions $\psi$ with compact
support on $\CC$, for $\ell$ big enough. 
\etheo
 
\dem The projective action of $\PSLC$ on the Riemann sphere
$\PP^1(\CC)=\CC\cup\{\infty\}$ identifies by the Poincar\'e extension
with the group of orientation preserving isometries of the upper
halfspace model of the real hyperbolic space $\htr$.  We denote the
image in $\PSLC$ of any subgroup $G$ of $\SLC$ by $\ov G$ and of any
element $g$ of $\SLC$ again by $g$.  The {\em Bianchi group}
$\Ga_K=\SLOK$ is a (nonuniform) arithmetic lattice in $\SLC$, whose
covolume is given by Humbert's formula (see \cite[\S 8.8 and
9.6]{ElsGruMen98})
\begin{equation}\label{eq:covolBianchi}
\Vol(\,\overline{\Ga_K}\,\bs\htr)=
\frac{|D_K|^{\frac 32}\zeta_K (2)}{4\pi^2}\,.
\end{equation}
 
Let $\Ga_{x,\,y}$ be the stabiliser of $(x,y)\in\OOO_K\times\OOO_K$
under the linear action of $\Ga_K$ on $\OOO_K\times\OOO_K$. In
particular, $ \Ga_{1,\,0}$ is the upper triangular unipotent subgroup of
$\Ga_K$, acting on $ \OOO_K\times\OOO_K$ via horizontal shears.
Furthermore, $x\OOO_K+y\OOO_K=u\OOO_K+v\OOO_K$ if and only if
$(u,v)\in\Ga_K(x,y)$, and the ideal class group of $K$ corresponds
bijectively to the set $\overline{\Ga_K}\,\bs\PP^1(K)$ of cusps of the
orbifold $\overline{\Ga_K}\,\bs\htr$ by the map induced by
$x\OOO_K+y\OOO_K\mapsto \rho=\frac xy\in K\cup\{\infty\}$, see for
example Section 7.2 of \cite{ElsGruMen98} for details. 

To prove the first claim (it follows by integration from the second
one, but the essential ingredients of the proof of either statement are the same), 
note that by the above, if
$\mmm=x\OOO_K+y\OOO_K$, then
$$
\psi_\mmm(s)=
\card\ \Ga_{1,\,0}\bs\big\{(u,v)\in\Ga_K(x,y):|v|^2\le \n(\mmm)s\big\}\,.
$$

If $y= 0$, let $\ga_\rho=\id$ and we may assume that $x=1$ since
$\psi_\mmm$ depends only on the ideal class of $\mmm$. If $y\neq 0$,
let
$$
\ga_\rho=\begin{pmatrix}\rho& -1\\1& \;\,0\end{pmatrix}\in\SLC\;.
$$
Let $\tau\in\;]0,1]$. Let $\H_\infty$ be the horoball in $\htr$ that
consists of all points with Euclidean height at least $1/\tau$ and let
$\H_\rho=\ga_\rho\H_\infty$. For any $\rho'\in K\cup \{\infty\}$, let
$\overline{\Ga_{\H_{\rho'}}}$ be the stabiliser in $\overline{\Ga_K}$
of the horoball $\H_{\rho'}$.  As in \cite{ParPau11BLMS} after its
Equation (4), we fix $\tau$ small enough such that $\H_\infty$ and
$\H_\rho$ are precisely invariant under $\overline{\Ga_K}$ (that is,
it meets one of its image by an element of $\overline{\Ga_K}$ only if
it coincides with it).

For any $g\in\SLC$ such that $g\H_\infty$ and $\H_\infty$ are
disjoint, it is easy to check using the explicit expression of the
Poincar\'e extension (see for example \cite[Eq.~4.1.4.]{Beardon83})
that the length of the common perpendicular of $\H_\infty$ and
$g\H_\infty$ is $|\ln(\tau^{-2}|c|^2)|$, where $c$ is the $(2,1)$-entry
of $g$.  If $y\ne 0$, then $(u,v)=g(x,y)$ if and only if $(\frac
uy,\frac vy)=g\ga_\rho(1,0)$, and the $(2,1)$-entry of $g\ga_\rho$ is
$\frac vy$.  Thus, the length $\ell(\delta_g)$ of the common
perpendicular $\delta_g$ between $\H_\infty$ and $g\H_\rho= g\ga_\rho
\H_\infty$, if these two horoballs are disjoint, is $\big|\!\ln(\frac{|v|^2}
{\tau^2|y|^2})\big|$ if $y\ne 0$ and $|\ln(\tau^{-2}|v|^2)|$ otherwise.

For the rest of the proof, we concentrate on the case $y\ne 0$. The
case $y=0$ is treated similarly. By discreteness, there are only
finitely many double classes $[g]\in \Ga_{1,\,0}\bs\Ga_K/\Ga_{x,\,y}$
such that $\H_\infty$ and $g\H_\rho$ are not disjoint or such that
$|v|\leq |y|$ or such that the multiplicity of $\delta_g$ is different
from $1$. Since the stabilisers $\Ga_{1,\,0}$ and
$\Ga_{x,\,y}$ do not contain $-\id$, we have
\begin{align*}
\psi_\mmm(s)&= \card\ \big\{[g]\in\Ga_{1,\,0}\bs\Ga_K/\Ga_{x,\,y}\;:\;
\ell(\delta_g)\le \ln\frac{\n(\mmm)s}{\tau^2|y|^2}\;\big\} 
+\operatorname{O}(1)\\ & =2\, \card\  
\big\{[g]\in\overline{\Ga_{1,\,0}}\bs\overline{\Ga_K}/\overline{\Ga_{x,\,y}}
\;:\;\ell(\delta_g)\le
\ln\frac{\n(\mmm)s}{\tau^2|y|^2}\;\big\} +\operatorname{O}(1)\,.
\end{align*}
We use \cite[Lem.~7]{ParPau11BLMS} with $C=\overline{\Ga_K}$,
$A=\overline{\Ga_{1,\,0}}$, $A'=\overline{\Ga_{\H_\infty}}$,
$B=\overline{\Ga_{x,\,y}}$, $B'=\overline{\Ga_{\H_\rho}}$, since there
are only finitely many $[g]\in A'\bs C/ B'$ such that $g^{-1}A'g\cap
B'\neq \{1\}$. Since $[\overline{\Ga_{\H_\infty}}:\overline{\Ga_{1,\,0}}]
=[\overline{\Ga_{\H_\rho}}:\overline{\Ga_{x,\,y}}]=\frac{\omega_K}{2}$,
we then have
$$
\psi_\mmm(s)= \frac{{\omega_K}^2}{2}\;\card\big\{[g]\in
\overline{\Ga_{\H_\infty}}\,\bs\overline{\Ga_K}/
\overline{\Ga_{\H_\rho}}\;:\;\ell(\delta_g)\le
\ln\frac{\n(\mmm)s}{\tau^2|y|^2}\;\big\}+\operatorname{O}(1)
\,.
$$

By \cite[Lem.~6]{ParPau11BLMS} and the last equality of the proof of
Theorem 4 on page 1055 of loc.~cit., we have
$\Vol(\,\overline{\Ga_{\H_\infty}}\,\bs\H_\infty)=
\frac{\tau^2\sqrt{|D_K|}}{2\,{\omega_K}}$ and
$\Vol(\,\overline{\Ga_{\H_\rho}}\,\bs\H_\rho)=
\frac{\tau^2\sqrt{|D_K|}}{2\,{\omega_K}}\frac{|y|^4}{\n(\mmm)^2}$. Thus,
Corollary \ref{coro:constcurvcount}, applied to $n=3$, $\Ga=
\overline{\Ga_K}$, $D^-=\H_\infty$ and $D^+=\H_\rho$, gives, since the
pointwise stabilisers of $\H_\infty$ and $\H_\rho$ are trivial,
$$
\psi_\mmm(s) = 
\frac{{\omega_K}^2}{2}\;\frac{2^2\,2\,\tau^4\,|D_K|\,|y|^4\,4\pi^2}
{4\,{\omega_K}^2\,\n(\mmm)^2\;4\pi\,|D_K|^{\frac 32}\,\zeta_K(2)}
\frac{\n(\mmm)^2}{\tau^4|y|^4}\,s^2(1+\bigO(s^{-\kappa}))\,,
$$
which, after simplification, proves  the first claim.

To prove the second claim when $y\ne 0$, note that Equation
\eqref{eq:distribhorob} in Corollary \ref{coro:constcurvcount} implies
that as $t\to+\infty$, with the appropriate error term,
$$
\frac{|D_K|\,\zeta_K(2)\,\omega_K}{2\,\pi}
\frac{\n(\mmm)^2}{\tau^2\,|y|^4}\,e^{-2t}
\sum_{z\in\partial\H_\infty}m_t(z)\,\Delta_z\weakstar\Vol_{\partial\H_\infty}\,.
$$ Using the change of variable $t=\ln\frac{\n(\mmm)s}{\tau^2|y|^2}$,
the fact that the geodesic lines containing the common perpendiculars
from $\H_\infty$ to the images of $\H_\rho$ by the elements of
$\overline{\Ga_K}$ are exactly the geodesic lines from $\infty$ to an
element of $\overline{\Ga_K}\cdot\frac xy$, and the above facts on
disjointness, lengths and multiplicities, this gives, since the map
$(u,v)\mapsto u/v$ is $\omega_K$-to-$1$,
$$
\frac{|D_K|\,\zeta_K(2)\,\omega_K\,\tau^2}{2\,\pi\,s^2}\;\frac{1}{\omega_K}
\sum_{\substack{(u,v)\in\mmm\times\mmm\;\\
\;\n(\mmm)^{-1}|v|^2\leq s,
\ \OOO_K u+\OOO_K v=\mmm}}\Delta_{(\frac uv,\,\frac{1}{\tau})}
\;\weakstar\;\Vol_{\partial\H_\infty}\,,
$$
as $s\to+\infty$.  The claim now follows from the observation that the
pushforward of $\Vol_{\partial\H_\infty}$ by the endpoint map
$(x,\frac{1}{\tau})\mapsto x$ is $\tau^2\Leb_\CC$. The proof of the
case $y=0$ is similar.  \cqfd

\subsection{A Mertens' formula for maximal orders of 
rational quaternion algebras}
\label{subsect:mertenshamilton}

Let $\HH$ be Hamilton's quaternion algebra over $\RR$, with $x\mapsto
\overline{x}$ its conjugation, $\redn: x\mapsto x\overline{x}$ its
reduced norm and $\redtr: x\mapsto x+\overline{x}$ its reduced
trace. We endow $\HH$ with its standard Euclidean structure (making
its standard basis orthonormal). Let $A$ be a quaternion algebra over
$\QQ$, which is definite ($A\otimes_\QQ\RR= \HH$), with reduced
discriminant $D_A$. Let $\OOO$ be a maximal order in $A$, and let
$\mmm$ be a (nonzero) left ideal of $\OOO$, with reduced norm
$\redn(\mmm)$ (see \cite{Vigneras80} for definitions).

The additive group $\OOO$ acts on the left on $\HH\times
\HH$ by the horizontal shears (transvections) $o\cdot(u,v) =(u+ov,v)$.
Let us consider the following counting function of the generating
pairs of elements of $\mmm$, defined for every $s>0$, by
$$
\psi_{\mmm}(s)=\card\;\OOO\,\bs
\big\{(u,v)\in\mmm\times\mmm\;:\;\redn(\mmm)^{-1}\redn(v)\leq s,
\;\;\;\OOO u+\OOO v=\mmm\big\}\,,
$$

\btheo\label{theo:appliHam} There exists $\kappa>0$ such that, as $s$
tends to $+\infty$,
$$
\psi_{\mmm}(s)= 
\frac{90\, D_A^2}{\pi^2\,\zeta(3)\,\prod_{p|D_A}(p^3-1)}
\;\;s^4+\operatorname{O}(s^{4-\kappa}) \;,
$$
with $p$ ranging over positive rational primes.  Furthermore,
$$
\frac{\pi^2\,\zeta(3)\,\prod_{p|D_A}(p^3-1)}{360\,D_A\; s^4}
\sum_{\substack{(u,v)\in\mmm\times\mmm\;\\
\;\redn(\mmm)^{-1}\redn(v)\leq s,
\ \OOO u+\OOO v=\mmm}}\Delta_{uv^{-1}}\;\weakstar\;\Leb_\HH
$$
as $s\to+\infty$, with error term $\bigO(s^{4-\kappa}\|\psi\|_\ell)$
when evaluated on $\C^\ell$-smooth functions $\psi$ with compact
support on $\HH$, for $\ell$ big enough.  \etheo

\dem We denote by $[\mmm]$ the ideal class of a left
fractional ideal $\mmm$ of $\OOO$, and by $_\OOO\!\I$ the set of left
ideal classes of $\OOO$, and by $\OOO_r(I)=\{x\in A\;:\; Ix \subset
I\}$ the right order of a $\ZZ$-lattice $I$ in $A$ (see
\cite{Vigneras80} for definitions).  Note that $\psi_{\mmm}$ depends
only on $[\mmm]$. For every $(u,v)$ in $A\times A-\{(0,0)\}$, consider
the two left fractional ideals of $\OOO$
$$
I_{u,\,v}=\OOO u+\OOO v\;,\;\;K_{u,\,v}=\Big\{\begin{array}{cl}
\OOO u\cap\OOO v &{\rm if~} uv\neq 0\;,\\ 
\OOO&{\rm otherwise.}\end{array}
$$
Given $\mmm,\mmm'$ two left fractional ideals of $\OOO$ and $s>0$, let
$$
\psi_{\mmm,\,\mmm'}(s)=\card\;\;\OOO\,\bs
\big\{(u,v)\in\mmm\times\mmm\;:\;\redn(v)\leq \redn(\mmm)s,
\;\;\;I_{u,\,v}=\mmm,\;[K_{u,\,v}]=[\mmm']\big\}\,,
$$
so that
\begin{equation}\label{eq:somdeunadeuidealfrac}
\psi_{\mmm}=\sum_{[\mmm']\;\in\;_\OOO\!\I}\;\psi_{\mmm,\,\mmm'}\;.
\end{equation}
We will give in Equation \eqref{eq:psimmprime} an asymptotic to
$\psi_{\mmm,\,\mmm'}(s)$ as $s\ra+\infty$, and in Equation
\eqref{eq:equidistribtwoideals} the related equidistribution result,
for each $[\mmm']\;\in\;_\OOO\!\I$ (interesting in themselves), and
then easily infer Theorem \ref{theo:appliHam}.

We refer for instance to \cite{Kellerhals03}, \cite[\S 3]{ParPau13ANT}
for the following properties of quaternionic homographies. Let $\SLH$
be the group of $2\times 2$ matrices with coefficients in $\HH$ and
Dieudonn\'e determinant $1$. Recall that the Dieudonné determinant of
$\begin{pmatrix} a & b \\c & d\end{pmatrix}$ is
\begin{equation}\label{eq:detdieud}
\redn(ad) +\redn(bc)-\redtr(a\,\overline{c}\,d\,\overline{b}\,)\,.
\end{equation}
The group $\SLH$ acts linearly on the left on the right $\HH$-module
$\HH\times\HH$, hence projectively on the left on the right projective
line $\PP^1_r(\HH)= (\HH\times \HH-\{0\})/\HH^\times$. The group
$\PSLH=\SLH/\{\pm {\rm Id}\}$ identifies by the Poincar\'e extension
procedure with the group of orientation preserving isometries of the
upper halfspace model $\HH\times\,]0,+\infty[$ of the real hyperbolic
space $\hcr$ of dimension $5$, with Riemannian metric 
$$ds^2(x)=
\frac{ds^2_\HH(z) +dr^2}{r^2}$$ at the point $x=(z,r)$.  For any
subgroup $G$ of $\SLH$, we denote by $\overline G$ its image in
$\PSLH$.

Let $\Ga_\OOO=\SLO=\SLH \cap\M_2(\OOO)$ be the {\it Hamilton-Bianchi
  group} of $\OOO$, which is a (nonuniform) arithmetic lattice in the
connected real Lie group $\SLH$ (see for instance \cite[page
1104]{ParPau10GT}).  Given $(x,y)\in \OOO\times \OOO$, let
$\Ga_{x,\,y}$ be the stabiliser of $(x,y)$ for the left linear action
of $\Ga_\OOO$.  By \cite[Rem.~7]{ParPau13ANT}, the map from the set
$\overline{\Ga_\OOO}\,\bs\PP^1_r(\OOO)$ of cusps of
$\overline{\Ga_\OOO}\,\bs\hcr$ into $_\OOO\!\I\times \;_\OOO\!\I$
which associates, to the orbit of $[u:v]$ in $\PP^1_r(\OOO)$ under
$\overline{\Ga_\OOO}$, the pair of ideal classes $([I_{u,\,v}],
[K_{u,\,v}])$ is a bijection. We hence fix a (nonzero) element
$(x,y)\in\OOO\times\OOO$ such that $[I_{x,\,y}]=[\mmm]$ and
$[K_{x,\,y}] = [\mmm']$, assuming that $x=1$ if $y=0$.

Let us prove the first claim.  Since the reduced norm of an invertible
element of $\OOO$ is $1$, the index in $\Ga_{1,\,0}$ of its unipotent
upper triangular subgroup is $|\OOO^\times|$, and we have
$$
\psi_{\mmm,\,\mmm'}(s)=|\OOO^\times|\;
\card\;\;\Ga_{1,\,0}\bs
\big\{(u,v)\in\Ga_\OOO(x,y)\;:\;\redn(v)\leq \redn(\mmm)s\big\}\,.
$$

Let $\tau\in\;]0,1]$. Let $\H_\infty$ be the horoball in $\hcr$
consisting of the points of Euclidean height at least $1/\tau$. By
\cite[Lem.~6.7]{ParPau10GT}, if $c$ is the $(2,1)$-entry of a matrix
$g\in \SLH$ such that $\H_\infty$ and $g\H_\infty$ are disjoint, then
the hyperbolic distance between $\H_\infty$ and $g\H_\infty$ is $|\ln
(\tau^{-2}\redn(c))|$.

Let $\rho=xy^{-1}\in A\cup\{\infty\}$. If $y= 0$, let $\ga_\rho=\id$,
otherwise let
$$
\ga_\rho=\Big(\!\begin{array}{cc}\rho& -1\\1&
  \;\,0\end{array}\!\Big)\in\SLH\;.
$$
Let $\H_\rho=\ga_\rho\H_\infty$, which is a horoball centered at
$\rho$.  For every $g\in \Ga_\OOO$ such that $\H_\infty$ and
$g\H_\rho$ are disjoint, let $\delta_g$ be the common perpendicular
from $\H_\infty$ to $g\H_\rho= g\ga_\rho\H_\infty$, with length
$\ell(\delta_g) >0$.  Hence if $\H_\infty$ and $g\H_\rho$ are
disjoint, since $(u,v)=g(x,y)$ if and only if $(uy^{-1},vy^{-1})=
g\ga_\rho(1,0)$ when $y\neq 0$, we have $\ell(\delta_g)=\big|\ln
\frac{\redn(v)}{\tau^2\redn(y)}\big|$ if $y\neq 0$ and
$\ell(\delta_g)=|\ln (\tau^{-2}\redn(v))|$ otherwise. By
discreteness, there are only finitely many double classes $[g]\in
\Ga_{1,\,0}\bs\Ga_\OOO/\Ga_{x,\,y}$ such that $\H_\infty$ and
$g\H_\rho$ are not disjoint or such that $\redn(v)\leq \redn(y)$
or such that the multiplicity of $\delta_g$ is different from
$1$. Assume that $y\neq 0$ (the case $y=0$ is treated
similarly). Then, as $s\ra+\infty$,
\begin{align*}
\psi_{\mmm,\,\mmm'}(s)&=|\OOO^\times|\; \card
\big\{[g]\in \Ga_{1,\,0}\bs\Ga_\OOO/\Ga_{x,\,y}\;:\;\ell(\delta_g)\leq 
\ln \frac{\redn(\mmm)s}{\tau^2\redn(y)}\big\}+\operatorname{O}(1)\,.
\\&=2\,|\OOO^\times|\; \card \big\{[g]\in
\overline{\Ga_{1,\,0}}\;\bs\overline{\Ga_\OOO}/\;\overline{\Ga_{x,\,y}}
\;:\; \ell(\delta_g)\leq
\ln \frac{\redn(\mmm)s}{\tau^2\redn(y)}\big\}+\operatorname{O}(1)\,.
\end{align*}

Let $\overline{\Ga_{\H_\rho}}$ be the stabiliser in
$\overline{\Ga_\OOO }$ of the horoball $\H_\rho$. By
\cite[Lem.~7]{ParPau11BLMS}, the number $ \psi_{\mmm,\,\mmm'}(s)$ is
equal to
$$
2\,|\OOO^\times|\, 
[\overline{\Ga_{\H_\infty}}:\overline{\Ga_{1,\,0}}]
[\overline{\Ga_{\H_\rho}}:\overline{\Ga_{x,\,y}}]\; 
\card \big\{[g]\in \overline{\Ga_{\H_\infty}}\;\bs
\overline{\Ga_\OOO}/\;\overline{\Ga_{\H_\rho}}
\;:\; \ell(\delta_g)\leq \ln 
\frac{\redn(\mmm)s}{\tau^2\redn(y)}\big\}+\operatorname{O}(1)\,.
$$

Note that $[\overline{\Ga_{\H_\infty}}:\overline{\Ga_{1,\,0}}]=
\frac{|\OOO^\times|}{2}$.  If $\tau$ is small enough, then $\H_\infty$
and $\H_\rho$ are precisely invariant under $\overline{\Ga_\OOO}$.
Hence, using Corollary \ref{coro:constcurvcount}, applied to $n=5$,
$\Ga= \overline{\Ga_\OOO}$, $D^-=\H_\infty$ and $D^+=\H_\rho$, we
have, since the pointwise stabilisers of $\H_\infty$ and $\H_\rho$ are
trivial,
\begin{align*}
\psi_{\mmm,\,\mmm'}(s)& = |\OOO^\times|^2
[\,\overline{\Ga_{\H_\rho}}:\overline{\Ga_{x,\,y}}\,]\;
\N_{D^-,\,D^+}\big(\ln \frac{\redn(\mmm)s}{\tau^2\redn(y)}\big)
+\operatorname{O}(1)\\ & =|\OOO^\times|^2
[\,\overline{\Ga_{\H_\rho}}:\overline{\Ga_{x,\,y}}\,]\;c(D^-,D^+)\;
\big(\frac{\redn(\mmm)s}{\tau^2\redn(y)}\big)^4
\big(1+\operatorname{O}(e^{-\kappa t})\big)\,.
\end{align*}
By the Remark before Lemma 15 of \cite{ParPau13ANT}, we have
$\Vol(\,\overline{\Ga_{\H_\infty}}\,\bs\H_\infty)=
\frac{D_A\,\tau^4}{8\,|\OOO^\times|^2}$, and by Lemma 15, Equation
(33) and the centred equation three lines before Remark 19 in
\cite{ParPau13ANT},
\begin{equation}\label{eq:xycuspvol}
\Vol(\,\overline{\Ga_{\H_\rho}} \bs \H_\rho)=
\frac{D_A\,\redn(y)^4\,\tau^4}{16\,|\OOO_r(K_{x,\,y})^\times|\,  
[\,\overline{\Ga_{\H_\rho}}:\overline{\Ga_{x,\,y}}\,]\,\redn(\mmm)^4}\,.
\end{equation}
By Theorem 3 of \cite{ParPau13ANT} (due to Emery, see the appendix of
\cite{ParPau13ANT}), we have
\begin{equation}\label{eq:covolHamilton}
\Vol(M)=\frac{\zeta(3)\prod_{p|D_A}(p^3-1)(p-1)}{11520}\,.
\end{equation}
Recall that $\Vol(\SSS^4)=\frac{8\pi^2}{3}$. Thus,
$\psi_{\mmm,\,\mmm'}(s)$ equals
\begin{align}\label{eq:psimmprime}
&\frac{|\OOO^\times|^2
[\,\overline{\Ga_{\H_\rho}}:\overline{\Ga_{x,\,y}}\,]
\;2^4\cdot 4\cdot 3\cdot 11520\, 
D_A^2\, \tau^8\,\redn(y)^4\;(\redn(\mmm)s)^4
\,(1+\operatorname{O}(s^{-\kappa}))}{8\pi^2\,\zeta(3)\,
{\displaystyle\prod_{p|D_A}}(p^3-1)(p-1)\;8\,|\OOO^\times|^2\,16\,
|\OOO_r(K_{x,\,y})^\times|\,[\overline{\Ga_{\H_\rho}}:\overline{\Ga_{x,\,y}}]
\;\redn(\mmm)^4\,\redn(y)^4\,\tau^8 } \nonumber
  \\
 & = \;\frac{2160\,D_A^2}{\pi^2\,\zeta(3)\,|\OOO_r(K_{x,\,y})^\times|\,
\prod_{p|D_A}(p^3-1)(p-1)}\;s^4\big(1+\operatorname{O}(s^{-\kappa})\big)\,.
\end{align}
By page 134 of \cite{Deuring68} (see Equation (8) of
\cite{ParPau13ANT}), we get
\begin{equation}\label{eq:deuring}
\sum_{[\mmm']\in\;_\OOO\!\I}\frac 1{|\OOO_r(\mmm')^\times|}=
\frac 1{24}\prod_{p|D_A}(p-1)\,.
\end{equation}
Thus, Equation \eqref{eq:somdeunadeuidealfrac} and the above
computations give the first claim of Theorem \ref{theo:appliHam}.

To prove the second claim (when $y\ne 0$, the case $y=0$ being
similar), Corollary \ref{coro:constcurvcount} implies that
$$
\frac{\pi^2\,\zeta(3)\,\prod_{p|D_A}(p^3-1)(p-1)\;
|\OOO_r(K_{x,\,y})^\times|\,
[\,\overline{\Ga_{\H_\rho}}:\overline{\Ga_{x,\,y}}\,]\;\redn(\mmm)^4}
{4320\,D_A\,\redn(y)^4\;\tau^4}
\;e^{-4t}\sum_{z\in\partial\H_\infty}m_t(z)\,\Delta_z
$$
weak-star converges to $\Vol_{\partial\H_\infty}$ as $t\to+\infty$,
with the appropriate error term.  Using the change of variable
$t=\ln\frac{\redn(\mmm)s}{\tau^2\redn(y)}$, the fact that the
geodesic lines containing the common perpendicular from $\infty$ to an
image of $\H_\rho$ by an element of $\overline{\Ga_\OOO}$ are exactly
the geodesic lines from $\infty$ to an element of
$\overline{\Ga_\OOO}\cdot xy^{-1}$, and the above facts on
disjointness, lengths and multiplicities, this gives, since the map
$(u,v)\mapsto uv^{-1}$ is $|\OOO^\times|$-to-$1$,
$$
\frac{\pi^2\,\zeta(3)\,\prod_{p|D_A}(p^3-1)(p-1)\;
|\OOO_r(K_{x,\,y})^\times|\,
[\,\overline{\Ga_{\H_\rho}}:\overline{\Ga_{x,\,y}}\,]}
{4320\,D_A\;|\OOO^\times|\;s^4\;\tau^{-4}}
\sum_{\substack{(u,\,v)\in\mmm\times\mmm,\;\redn(\mmm)^{-1}\redn(v)\leq s\\
I_{u,\,v}=\mmm,\;[K_{u,\,v}]=[\mmm']}} \Delta_{(uv^{-1},\,\frac 1\tau)}
$$
weak-star converges to $\Vol_{\partial\H_\infty}$ as $s\to+\infty$.
Since $[\,\overline{\Ga_{\H_\rho}}:\overline{\Ga_{x,\,y}}\,]=
\frac{|\OOO^\times|}{2}$ and the pushforward of
$\Vol_{\partial\H_\infty}$ by the endpoint map $(x,\frac 1\tau)\mapsto
x$ is $\tau^4\Leb_\HH$, we have
\begin{equation}\label{eq:equidistribtwoideals}
\frac{\pi^2\,\zeta(3)\,\prod_{p|D_A}(p^3-1)(p-1)\;
|\OOO_r(K_{x,\,y})^\times|}{8640\,D_A\;s^4}
\sum_{\substack{(u,\,v)\in\mmm\times\mmm,\;\redn(\mmm)^{-1}\redn(v)\leq s\\
I_{u,\,v}=\mmm,\;[K_{u,\,v}]=[\mmm']}}\Delta_{uv^{-1}}
\;\weakstar\;\Leb_\HH\,,
\end{equation}
The second claim of Theorem \ref{theo:appliHam} now follows by
dividing both sides by $|\OOO_r(K_{x,\,y})^\times|$ and summing over
$[m']$ as above. \cqfd

\section{Counting and equidistribution of quadratic irrationals 
and  crossratios}
\label{sect:irrationals}

Let $K$ be either $\QQ$ or an imaginary quadratic number field, and,
respectively, $\wh K=\RR$ or $\wh K=\CC$. Let $\OOO_K$ be the ring of
integers, $D_K$ the discriminant, and $\zeta_K$ the zeta function of
$K$. In this section, we denote by $\cdot$ the action by homographies
of the group $\Ga_K=\PSLOK$ on $\PP^1(\wh K)= \wh K\cup\{\infty\}$.
For every finite index subgroup $G$ of $\Ga_K$, and every $x\in \wh
K\cup\{\infty\}$, let $G_x$ be the stabiliser of $x$ in $G$, with
$\Ga_x=(\Ga_K)_x$ to simplify the notation.

For every $\alpha\in \wh K$ which is quadratic irrational over $K$,

$\bullet$~ denote by $\alpha^\sigma$ its Galois conjugate over $K$, by
$\tr \alpha=\alpha+\alpha^\sigma$ and $\n(\alpha)= \alpha
\alpha^\sigma$ its relative trace and relative norm, and by
$$
Q_\alpha(X)=X^2-\tr\alpha\; X +\n(\alpha)
$$ 
the standard monic quadratic polynomial with roots $\alpha$ and $\alpha^\sigma$;

$\bullet$~ denote by $h(\alpha)=\frac{2}{|\alpha-\alpha^\sigma|}$ the
natural complexity of $\alpha$ in an orbit of $\Ga_K$, modulo
translations by $\OOO_K$, introduced in \cite{ParPau11MZ} and
motivated in \cite[\S 4.1]{ParPau12JMD} (when $\alpha$ is integral,
$\frac{4}{h(\alpha)^2}=\discr Q_\alpha$ is the discriminant of the minimal
polynomial of $\alpha$);

$\bullet$~ if $K=\QQ$, let $q_\alpha$ be the least common multiple of
the denominators of the rationals $\tr\, \alpha$ and $\n(\alpha)$, let
$D_{\alpha}={q_\alpha}^2((\tr\, \alpha)^2 -4\n(\alpha))$ (when $\alpha$
is integral, we have $q_\alpha=1$ and $D_\alpha$ is the discriminant
of the order $\ZZ+\alpha\ZZ$), let $(t_\alpha,u_\alpha)$ be the
fundamental solution of the Pell-Fermat equation
$t^2-D_{\alpha}u^2=4$, and let $R_\alpha= \arcosh\frac{t_\alpha}{2}$
(which is the regulator of $\ZZ+\alpha\ZZ$ when $\alpha$ is integral);

$\bullet$~ for every finite index subgroup $G$ of $\Ga_K$, we define
the $G$-{\em reciprocity index} $\iota_G(\alpha)$ of $\alpha$ as
follows: we set $\iota_G(\alpha)=2$ if $\alpha$ is {\it
  $G$-reciprocal}, that is, if some element of $G$ maps $\alpha$ to
$\alpha^\sigma$, see \cite{Sarnak07} and \cite[Prop.~4.3]{ParPau12JMD}
for characterisations (when $K=\QQ$, $G=\Ga_K$ and $\alpha$ is
integral, this is equivalent to saying that the order $\ZZ+\alpha\ZZ$
contains a unit of norm $-1$), and we set $\iota_G(\alpha)=1$
otherwise.

To conclude with a geometric remark, for every $\alpha$ as above,
there exists a unique primitive loxodromic element
$\wh{\alpha}\in\Ga_K$ such that the repelling fixed point
$\wh{\alpha}^-$ of $\wh{\alpha}$ is equal to $\alpha$ (see for
instance \cite[Lem.~6.2]{ParPau11MZ}). Its attractive fixed point
$\wh{\alpha}^+$ is then $\alpha^\sigma$. With the notation at the end
of Section \ref{sect:geometry}, we have $\iota_G(\alpha)=
\iota_G(\wh{\alpha})$ for every finite index subgroup $G$ of $\Ga_K$.
Finally, for all $\ga\in\Ga_K$, we have $\wh{\ga\alpha}=
\ga\wh{\alpha}\ga^{-1}$.

\subsection{Equidistribution and error terms in counting 
functions of quadratic irrationals}
\label{subsect:errortermquadirrat}

We give in this subsection an error term to the counting asymptotics
of \cite{ParPau12JMD} of the number of quadratic irrationals
$\alpha\in \wh K$ over $K$ with complexity $h(\alpha)$ at most $s$, in
an orbit of a finite index subgroup of $\Ga_K$, and we prove an
equidistribution result for the set of traces $\tr\,\alpha$ of these
elements, as $s\ra+\infty$.

\btheo\label{theo:appliarithneg2} Let $\alpha_0$ be a real quadratic
irrational over $\QQ$ and let $G$ be a finite index subgroup of
$\Ga_\QQ$.  There exists $\kappa>0$ such that, as $s\ra+\infty$,
\begin{multline*}
\card\{\alpha\in G \cdot \alpha_0
+[\Ga_\infty:G_\infty]\ZZ\;:\; h(\alpha)\leq s\} \\= \frac{6
\;[\Ga_\infty:G_\infty]\,[\Ga_{\alpha_0}:G_{\alpha_0}]
\;R_{\alpha_0}}{\pi^2\;[\Ga_\QQ:G]}\;
\;s +\operatorname{O}(s^{1-\kappa})\,.
\end{multline*}
Furthermore, for the weak-star convergence of measures, 
we have
$$
\lim_{s\ra +\infty}\;\;\frac{\pi^2\;[\Ga_\QQ:G]}
{3\;[\Ga_{\alpha_0}:G_{\alpha_0}]\;R_{\alpha_0}\;s}\;\;
\sum_{\alpha\in G \cdot \alpha_0
\;:\; h(\alpha)\leq s}\Delta_{\tr\,\alpha}\;=\; \Leb_\RR\,.
$$
\etheo

Note that the indices (except maybe $[\Ga_{\alpha_0}:G_{\alpha_0}]$)
appearing in the above statement are for instance known when $G$ is
the principal congruence subgroup modulo a prime $p$, or the Hecke
congruence subgroup modulo $p$. When $G=\Ga_\QQ$, the first claim is
known (with a better error term, see for instance \cite[page
164]{Cohn62}), but it is new in particular when $G$ is not a
congruence subgroup of $\Ga_\QQ$.  For smooth functions $\psi$ with
compact support on $\RR$, there is an error term in the
equidistribution claim, of the form $\bigO(s^{1-\kappa}\|\psi\|_\ell)$
where $\kappa>0$ and $\|\psi\|_\ell$ is the Sobolev norm of $\psi$ for
some $\ell\in\NN$, see the comment following Theorem
\ref{theo:mainequicountdown}. Theorem \ref{theo:equidtraceintro} in
the Introduction follows from this result.

\medskip \dem The proof of the counting result with an error term is
similar to that of \cite[Th\'eo.~4.4]{ParPau12JMD} upon replacing
\cite[Coro.~3.10]{ParPau12JMD} by the above Corollary
\ref{coro:constcurvcount} (the factor $\frac{2}{\iota_G(\alpha_0)}$
which appears in \cite[Th\'eo.~4.4]{ParPau12JMD} comes from the fact
that we are counting $\alpha\in G \cdot \{\alpha_0,\alpha_0^\sigma\}$
therein). Noting that $\tr (\alpha+m)=\tr \alpha+2m$ for all
$m\in\ZZ$, it also follows by applying the equidistribution result
with error term on the quotient of $\RR$ by $2[\Ga_\infty:G_\infty]
\,\ZZ$ to the constant function $1$.

To prove the second claim, we apply Equation \eqref{eq:distribhorob}
with $n=2$, $\Ga=G$, $D^-$ the horoball centred at $\infty$ consisting
of the points with vertical coordinates at least $1$ in the upper
halfplane model of $\hdr$, and $D^+$ the geodesic line with points at
infinity $\alpha_0$ and $\alpha_0^\sigma$, whose stabiliser in $G$ we
denote by $G_{D^+}$. The image of $D^+$ by an element of $G$ is the
geodesic line $]\alpha, \alpha^\sigma[$ with points at infinity
$\alpha$, $\alpha^\sigma$ for some $\alpha\in G \cdot \alpha_0$. Such
an $\alpha$ is uniquely determined if $\alpha_0$ is not
$G$-reciprocal, otherwise there are exactly two choices, $\alpha$ and
$\alpha^\sigma$.  The origin of the common perpendicular from $D^-$ to
$]\alpha, \alpha^\sigma[$ (which exists except for finitely many
$G_\infty$-orbits of $\alpha\in G \cdot \alpha_0$) is the point
$x_\alpha=(\frac{\tr \,\alpha}{2},1) \in\partial D^-$, its length is
$\ln h(\alpha)$, and its multiplicity is $1$, since $\PSLZ$ acts
freely on $T^1\hdr$. The induced Riemannian measure on $\partial D^-$
is the Lebesgue measure. Therefore, by Equation
\eqref{eq:distribhorob} and the value of the constant $c(D^-,D^+)$
given above Corollary \ref{coro:constcurvcount}, we have
\begin{equation}\label{eq:equidquadirratQ}
\lim_{s\ra +\infty}\;\;\frac{\Vol(\SSS^{1})\;\Vol(G\bs \hdr)}
{\iota_{G}(\alpha_0)\;\Vol(\SSS^{0})\;\Vol(G_{D^+}\bs D^+)\;s}\;\;
\sum_{\alpha\in G \cdot \alpha_0
\;:\; h(\alpha)\leq s}\Delta_{x_\alpha}\;=\; 
\Leb_{\partial D^-}\;.
\end{equation}
We have 
\begin{equation}\label{eq:volumehdrmodG}
\Vol(G\bs\hdr)=[\Ga_\QQ:G]\,\Vol(\Ga_\QQ\bs\hdr)=[\Ga_\QQ:G]\,\frac\pi 3\,.
\end{equation}
Recall that $\wh{\alpha_0}\in\Ga_\QQ$ is a primitive loxodromic
element with fixed points $\alpha_0$ and $\alpha_0^\sigma$. By for
instance \cite[p.~173] {Beardon83}, the translation length
$\ell(\wh{\alpha_0})$ of $\ga_0$ satisfies $\cosh
\frac{\ell(\wh{\alpha_0})}{2}= \frac{|\operatorname{tr}
  {\wh{\alpha_0}}\,|} {2}$. By for instance the proof of
\cite[Prop.~4.1]{ParPau12JMD}, we have $|\operatorname{tr}
{\wh{\alpha_0}}\,| = t_{\alpha_0}$. Hence
\begin{equation}\label{eq:translengthregulator}
\ell(\wh{\alpha_0})= 2\,R_{\alpha_0}\,,
\end{equation}
so that
\begin{equation}\label{eq:longgeodmodG}
  \Vol(G_{D^+}\bs D^+)=
  \frac{[\Ga_{\alpha_0}:G_{\alpha_0}]}{\iota_G(\alpha_0)}\,
  \Vol(\Ga_{\alpha_0}\bs D^+)=
  \frac{[\Ga_{\alpha_0}:G_{\alpha_0}]\,\ell(\wh{\alpha_0})} 
  {\iota_G(\alpha_0)}=
  \frac{2\,[\Ga_{\alpha_0}:G_{\alpha_0}]\,R_{\alpha_0}}
  {\iota_G(\alpha_0)}\,,
\end{equation}
The result now follows from Equations \eqref {eq:equidquadirratQ},
\eqref{eq:volumehdrmodG} and \eqref{eq:longgeodmodG} by applying the
pushforwards of measures by the map $f:(x,1)\mapsto 2x$ from $\partial
D^-$ to $\RR^{n-1}$ (which sends the Lebesgue measure to
$\frac{1}{2^{n-1}}$ times the Lebesgue measure).  \cqfd

\btheo\label{theo:appliarithdim3} Let $K$ be an imaginary quadratic
number field, let $\alpha_0\in\CC$ be a quadratic irrational over $K$,
let $G$ be a finite index subgroup of $\Ga_K$, and let $\Lambda$ be
the lattice of $\lambda\in\OOO_K$ such that $\pm\Big(\begin{array}{cc}
  1&\lambda \\0&1\end{array} \Big)\in G$.  There exists $\kappa>0$
such that, as $s\ra+\infty$,
\begin{multline*}
\card\{\alpha\in G \cdot
\alpha_0+\Lambda\;:\; h(\alpha)\leq s\}
\\=\frac{\pi^2\;[\OOO_K:\Lambda]\,
[\Ga_{\alpha_0}:G_{\alpha_0}]
  \;\big|\log\big|\frac{\operatorname{tr}\,\wh{\alpha_0}+
    \sqrt{(\operatorname{tr}\,\wh{\alpha_0})^2-4}}{2}\big|\;\big|}
{m_{\Ga_K}(\alpha_0)\,[\Ga_K:G]\;|D_K|\;\zeta_{K}(2)}\; \;s^2 
+\operatorname{O}(s^{2-\kappa})\;,
\end{multline*}
where $\wh{\alpha_0}$ is a primitive element of $\Ga_K$ fixing
$\alpha_0$ with absolute values of eigenvalues different from $1$, and
$m_G(\alpha_0)$ the number of elements in $G$ fixing $\alpha_0$ with
absolute values of eigenvalues equal to $1$.  Furthermore, for the
weak-star convergence of measures, we have
$$
\lim_{s\ra +\infty}\;\;
\frac{2\,|D_K|^{\frac 32}\;\zeta_{K}(2)\;[\Ga_K:G]\;m_{\Ga_K}(\alpha_0)}
{\pi^2\;[\Ga_{\alpha_0}:G_{\alpha_0}]\;
\big|\log\big|\frac{\operatorname{tr}\,\wh{\alpha_0}+
\sqrt{(\operatorname{tr}\,\wh{\alpha_0})^2-4}}{2}\big|\;\big|\;\;s^2}\;\;
\sum_{\alpha\in G \cdot \alpha_0
\;:\; h(\alpha)\leq s}\Delta_{\tr\, \alpha}
\;=\; \Leb_\CC\,.
$$
\etheo

For smooth functions $\psi$ with compact support on $\CC$, there is an
error term in the above equidistribution claim, of the form
$\bigO(s^{2-\kappa} \|\psi\|_\ell)$ where $\kappa>0$ and
$\|\psi\|_\ell$ is the Sobolev norm of $\psi$ for some $\ell\in\NN$.

\medskip
\dem 
The proof of the counting claim with an error term is similar to that
of \cite[Th\'eo.~4.6]{ParPau12JMD} upon replacing
\cite[Coro.~3.10]{ParPau12JMD} by the above Corollary
\ref{coro:constcurvcount}, and using the simplification 
$$
[\Ga_\infty:\Lambda]=[\Ga_\infty:\OOO_K][\OOO_K:\Lambda]=
\frac{| {\OOO_K}^\times|}{2}\,[\OOO_K:\Lambda]\;.
$$

To prove the second claim, we apply Equation \eqref{eq:distribhorob}
with $n=3$, $\Ga=G$, $D^-$ the horoball centred at $\infty$ consisting
of the points with vertical coordinates at least $1$ in the upper
halfspace model of $\htr$, and $D^+$ the geodesic line with points at
infinity $\alpha_0$ and $\alpha_0^\sigma$, whose stabiliser in $G$ we
denote by $G_{D^+}$. The cardinality of the pointwise stabiliser in
$G$ of $D^+$ is $m_G(\alpha_0)$ as defined in the statement of Theorem
\ref{theo:appliarithdim3} (we have $m_G(\alpha_0)=m_G(\wh{\alpha_0})$
with the notation at the end of Section \ref{sect:geometry}).

The set of points in $\partial D^-$ fixed by a nontrivial element of
the stabiliser of $\infty$ in $\Ga_ K$ is discrete since, for
instance, the subgroup of translations in a discrete group of
isometries of $\RR^{n-1}$ has finite index by Bieberbach's theorem.
Hence the multiplicity of the common perpendicular from $D^-$ to an
image of $D^+$ by an element of $G$ is different from $1$ only
when its origin belongs to a discrete subset $S$ of $\partial D^-$.
Furthermore, since $\Ga_{D^+}\bs D^+$ is compact, there exists
$\epsilon >0$ (depending only on $\Ga_K$ and $\alpha_0$) such that a
geodesic arc, leaving perpendicularly from an image of $D^+$ and
arriving perpendicularly at another image of $D^+$, has length at
least $\epsilon$. Hence the number of common perpendiculars from a
given point in $S$ to an image of $D^+$ grows at most linearly in the
length.  Therefore these common perpendiculars do not contribute
asymptotically to the equidistributing sum, and as in the previous
proof, we have
\begin{equation}\label{eq:equidquadirratK}
\lim_{s\ra +\infty}\;\;
\frac{2\,m_G(\alpha_0)\,\Vol(\SSS^{2})\;\Vol(G\bs \htr)}
{\iota_{G}(\alpha_0)\;\Vol(\SSS^{1})\;\Vol(G_{D^+}\bs D^+)\;s^2}\;\;
\sum_{\alpha\in G \cdot \alpha_0
\;:\; h(\alpha)\leq s} \Delta_{(\frac{\tr \,\alpha}{2},1)}\;=\; 
\Leb_{\partial D^-}\;.
\end{equation}
Recall that the translation length $\ell(\ga_0)$ of an element $\ga_0$
of $\PSLC$ is
\begin{equation}\label{eq:formultranslength}
\ell(\ga_0)=2\big|\log\big|\frac{\operatorname{tr}\,\ga_0+
    \sqrt{{\operatorname{tr}\ga_0}^2-4}}{2}\big|\;\big|
\end{equation}
(which is independent of the choice of the square root of the complex
number $\operatorname{tr}\ga_0^2-4$ and of the choice of the lift of
$\ga_0$ in $\SLC$ modulo $\{\pm\id\}$). With $r$ the smallest positive
integer such that ${\wh{\alpha_0}}^r\in G$, we have
$[\Ga_{\alpha_0}:G_{\alpha_0}]=\frac{m_{\Ga_K}(\alpha_0)}{m_G(\alpha_0)}\,r$,
so that
\begin{equation}\label{eq:calcvollonggeod}
\Vol(G_{D^+}\bs D^+)=
\frac{r}{\iota_G(\alpha_0)}\,\ell(\wh{\alpha_0})=
\frac{m_G(\alpha_0)\,[\Ga_{\alpha_0}:G_{\alpha_0}]}
{m_{\Ga_K}(\alpha_0)\,\iota_G(\alpha_0)}\, \ell(\wh{\alpha_0})\,.
\end{equation}
As in the previous proof, the result now follows from Equations
\eqref{eq:equidquadirratK}, \eqref{eq:formultranslength} and
\eqref{eq:covolBianchi} by applying the pushforwards of measures by
the map $(x,1)\mapsto 2x$ from $\partial D^-$ to $\CC$. 
\cqfd

\bcoro \label{coro:orbianchi} Let $\phi=\frac{1+\sqrt{5}}2$ be the
Golden Ratio, let $K$ be an imaginary quadratic number field with
$D_K\neq -4$, let $\ccc$ be a nonzero ideal in $\OOO_K$, and let
$\Ga_{0}(\ccc)$ be the Hecke congruence subgroup
$\Big\{ \pm \Big(\!\begin{array}{cc}a & b \\ c & d 
\end{array} \!\Big)\in \Ga_K \;:\; c\in\ccc\Big\}$. Then there
exists $\kappa>0$ such that, as $s\ra+\infty$, $$
\card\{\alpha\in \Ga_{0}(\ccc) \cdot \phi + \OOO_K\,:\; 
h(\alpha)\leq s\}=
\frac{2\pi^2\;k_\ccc\;\log\phi} {|D_K|\;  \zeta_{K}(2)\,
  \n(\ccc)\prod_{\ppp |\ccc} \big(1+\frac{1}{\n(\ppp)}\big)}\; \;s^2+
\bigO(s^{2-\kappa})\;,
$$
where $k_\ccc$ is the smallest $k\in\NN-\{0\}$ such that the $2k$-th
term of the standard Fibonacci sequence belongs to $\ccc$, and the
product is over the prime ideals in $\OOO_K$ dividing $\ccc$.  
\ecoro

\dem The proof of this corollary, as well as the one of Corollary
\ref{coro:orbianchiintro} in the introduction, are similar to the one
of \cite[Coro.~4.7]{ParPau12JMD}, using the fact proven therein that
$m_{\Ga_0(\ccc)}(\phi)=1$, so that (the element $\wh{\phi}^{-1}$ is denoted
by $\ga_1=\begin{pmatrix} 2 & 1 \\ 1 & 1\end{pmatrix}$
in loc. cit., and its translation length is $2\log \phi$)
$$
\frac{[\Ga_{\phi}:(\Ga_{0}(\ccc))_{\phi}]}{m_{\Ga_K}(\phi)}
=\frac{[\wh{\phi}\,^\ZZ:\wh{\phi}\,^\ZZ\cap\Ga_{0}(\ccc)]}
{m_{\Ga_0(\ccc)}(\phi)}= k_\ccc\,.
\;\;\;\Box
$$

\bigskip 
Let $\HH$ be Hamilton's quaternion algebra over $\RR$, let $A$ be a
definite quaternion algebra over $\QQ$, and let $\OOO$ be a maximal
order in $A$ (we refer to Subsection \ref{subsect:mertenshamilton} for
background, as well as for the notation $\overline{\stackrel{\;\;}
  {\;}}$, $\redn$, $\redtr$, $D_A$, $\PSLH$, $\PSLO$, $\cdot$ ). For
every $x\in \HH\cup\{\infty \}$, we denote by $G_x$ the stabiliser of
$x $ for the action by homographies of a subgroup $G$ of $\PSLH$. For
every $\ga=\begin{pmatrix}a& b\\c& d\end{pmatrix} \in \SLH$, we denote
by $X_\ga$ the largest real root of
$$
2x^3 -c_1\, x^2 + 2(c_2c_3 -1\big)x+(c_1-{c_2}^2-{c_3}^2)\,,
$$
where $c_1=\redn(a+d)+\redtr(ad-bc)$, $c_2=\frac{1}{2}\redtr(a+d)$,
$c_3=\frac{1}{2}\redtr\big((ad-bc)\overline{a}+(da-cb)
\overline{d}\,\big)$. We will say that an element $\alpha\in\HH$ is a
{\it loxodromic quadratic irrational} over $\OOO$ if there exists $\ga
\in \SLO$ such that $|X_\ga|\neq 1$ fixing $\alpha$, in which case we
denote by $\alpha^\sigma$ its other fixed point (which exists, is
unique, and is independent of such a $\ga$, see below for proofs and
a justification of the terminology). Note that by the noncommutativity of
$\HH$, there are several types of quadratic equations over $\OOO$.

The following result is a counting and equidistribution result in
$\HH$ of loxodromic quadratic irrationals over $\OOO$ in a given
homographic orbit under $\PSLO$.

\btheo\label{theo:appliarithdim5} Let $\OOO$ be a maximal order in a
definite quaternion algebra over $\QQ$, let $\alpha_0\in\HH$ be a
loxodromic quadratic irrational over $\OOO$, let $G$ be a finite index
subgroup of $\PSLO$, and let $\Lambda$ be the lattice of
$\lambda\in\OOO$ such that $\pm\Big(\begin{array}{cc} 1&\lambda
  \\0&1\end{array} \Big)\in G$.  There exists $\kappa>0$ such that, as
$\epsilon\ra 0$,
\begin{multline*}
\card\{\alpha\in G \cdot
\alpha_0+\Lambda\;:\; \redn(\alpha-\alpha^\sigma)\geq \epsilon\}=
\\\frac{8640\;D_A\;[\OOO:\Lambda]\,[\PSLO_{\alpha_0}:G_{\alpha_0}]
  \;\big|\ln|X_{\ga_0} +({X_{\ga_0}}^2-1)^{1/2}|\;\big|}
{\zeta(3)\,m_{\PSLO}(\alpha_0)\;[\PSLO:G]\;
\prod_{p|D_A}(p^3-1)(p-1)}\; \;\epsilon^{-2} +\bigO(\epsilon^{\kappa-2})\;,
\end{multline*}
where $\ga_0$ is a primitive element of $\SLO$ fixing $\alpha_0$ with
$|X_{\ga_0}|\neq 1$, $m_G(\alpha_0)$ is the number of elements $\ga$
in $G$ fixing $\alpha_0$ with $|X_{\ga}|= 1$, and the product is over
the primes $p$ dividing $D_A$.  Furthermore, the Lebesgue measure
$Leb_\HH$ on $\HH$ is the weak-star limit, as $\epsilon\ra 0$, of the
measures
$$
\frac{\zeta(3)\;m_{\Ga_K}(\alpha_0)\;[\PSLO:G]\;
{\prod_{p|D_A}}(p^3-1)(p-1)\;\;\epsilon^2}
{2160\;[\PSLO_{\alpha_0}:G_{\alpha_0}]\;
\big|\ln|X_{\ga_0} +({X_{\ga_0}}^2-1)^{1/2}|\,\big|}\;
\sum_{\alpha\in G \cdot \alpha_0
\;:\; \redn(\alpha-\alpha^\sigma)\geq \epsilon}\Delta_{\alpha+\alpha^\sigma}\,.
$$
\etheo

For smooth functions $\psi$ with compact support on $\HH$, there is an
error term in the equidistribution claim, of the form $\bigO
(\epsilon^{\kappa-2} \|\psi\|_\ell)$ where $\kappa>0$ and $\|\psi\|_\ell$ is
the Sobolev norm of $\psi$ for some $\ell\in\NN$.

\medskip \dem First recall (see \cite[\S 3]{ParSho09}) that for every
$\ga \in \SLH$, we have $|X_\ga|\geq 1$, and that $\ga$ is loxodromic
for its isometric action on the upper halfspace model of $\hcr$ by
Poincaré's extension (see for instance Equation (14) in
\cite{ParPau13ANT}) if and only if $|X_\ga|\neq 1$. If this holds,
then (see \cite[\S 3]{ParSho09}) $\ga$ is conjugated to
$\begin{pmatrix}t_+& 0\\0& t_- \end{pmatrix}$ with $\redn(t_\pm)=
X_\ga\pm(X_\ga^2-1)^{1/2}$.  Since this diagonal matrix acts on the
vertical axis $\{0\} \times\; ]0,+\infty[$ in $\hcr$ by $(0,r)\mapsto
(0,\redn(t_+)r)$, the translation length $\ell(\ga)$ of $\ga$ is
therefore
\begin{equation}\label{eq:translengthhamilton}
\ell(\ga)=\big|\ln|X_\ga +(X_\ga^2-1)^{1/2}|\;\big|\;.
\end{equation}
As a side remark, if $\ga \in \SLO$ satisfies $|X_\ga|> 1$, in
particular, it has exactly two fixed points, which are the only two
solutions of some quadratic equation $ax+b=xcx+xd$ where $a,b,c,d
\in\OOO$ and $\redn(ad)+
\redn(bc)-\redtr(a\,\overline{c}\,d\,\overline{b}\,) =1$.

\medskip To prove the equidistribution claim of Theorem
\ref{theo:appliarithdim5}, we apply Equation \eqref{eq:distribhorob}
with $n=5$, $\Ga=G$, $D^-$ the horoball centred at $\infty$ consisting
of the points with vertical coordinate at least $1$ in $\hcr$, and
$D^+$ the geodesic line with points at infinity $\alpha_0$ and
$\alpha_0^\sigma$, whose stabiliser in $G$ we denote by $G_{D^+}$. We
define $\iota_{G}(\alpha_0)=2$ if there exists an element in $G$
sending $\alpha_0$ to ${\alpha_0}^\sigma$ and $\iota_{G}(\alpha_0)=1$
otherwise.  Except for finitely many $G_\infty$-orbits of $\alpha\in
G\cdot \alpha_0$, the common perpendicular from $D^-$ to the geodesic
line with points at infinity $\alpha$, $\alpha^\sigma$ exists and has
hyperbolic length $\ln \frac{2}{\redn(\alpha-\alpha^\sigma)^{1/2}}$, and
its origin is $(\frac{\alpha+\alpha^\sigma}{2},1)\in\partial D^-$. As
in the proof of Theorem \ref{theo:appliarithdim3}, we have
\begin{equation}\label{eq:equidquadirratH}
  \lim_{t\ra +\infty}\;\;\frac{4\;m_{G}(\alpha_0)\,
\Vol(\SSS^{4})\;\Vol(G\bs \hcr)}
  {\iota_{G}(\alpha_0)\,\Vol(\SSS^{3})\;\Vol(G_{D^+}\bs D^+)\;e^{4t}}\;\;
  \sum_{\alpha\in G \cdot \alpha_0\;:\; 
    \redn(\alpha-\alpha^\sigma)\geq 4\,e^{-2t}} 
\Delta_{(\frac{\alpha+\alpha^\sigma}{2},\,1)}\;=\; 
  \operatorname{Leb}_{\partial D^-}\;.
\end{equation}
We have $\Vol(\SSS^{4})=\frac{8\pi^2}{3}$, $\Vol(\SSS^{3})=2\pi^2$ and,
as in the proof of Theorem \ref{theo:appliarithdim3},
$$
\Vol(G_{D^+}\bs D^+)= \frac{[\PSLO_{\alpha_0}:G_{\alpha_0}]\,
  m_{G}(\alpha_0)} {\iota_G(\alpha_0)\,m_{\PSLO}(\alpha_0)}\,\ell(\ga_0)\,.
$$
Taking $\epsilon=4\,e^{-2\,t}$, the result follows, as in the proof of
Theorem \ref{theo:appliarithdim3}, from Equations
\eqref{eq:equidquadirratH}, \eqref{eq:translengthhamilton} and
\eqref{eq:covolHamilton} by applying the pushforwards of measures by
the map $(x,1)\mapsto 2x$ from $\partial D^-$ to $\HH$.

\medskip 
To prove the first claim of Theorem \ref{theo:appliarithdim5}, note
that by \cite[Prop.~5.5]{KraOse90} for instance, we have
$$
\Vol((2\Lambda)\bs \HH)=2^{4}\,[\OOO:\Lambda]\,\Vol(\OOO\bs \HH)=
2^{4}\,[\OOO:\Lambda]\,\frac{D_A}{4}\,.
$$
The result then follows by considering the measures induced on the
compact quotient $(2\Lambda)\bs \HH$, and by applying the equidistribution
result with the error term to the constant function $1$.  
\cqfd

\subsection{Relative complexity of loxodromic elements} 
\label{subsect:relatcomploxo}

This subsection is a geometric one, paving the way to the next one,
which is arithmetic. It explains in particular why the crossratios
play an important geometric role in this paper.

Let $n\geq 2$ and $\Ga$ be as in the beginning
of Section \ref{sect:geometry}, and let $\ga_0\in\Ga$ be a fixed
primitive loxodromic element.

We define the {\em relative height} of a primitive loxodromic
element $\ga\in\Ga$ with respect to $\gamma_0$ to be the length of the
common perpendicular of the translation axes of $\ga_0$ and $\ga$ if
they are disjoint and $0$ otherwise. Note that
$$
h_{\gamma_0^\epsilon}(\ga^{\epsilon'})=h_{\gamma_0}(\ga)=
h_{g\gamma_0g^{-1}}(g\ga g^{-1})=h_{\gamma_0}(g'\ga g'^{-1})
$$ 
for all $\epsilon,\epsilon'\in\{\pm 1\}$, $g\in\Ga$ and
$g'\in\stab_\Ga(\axis\ga_0)$.  Furthermore, it follows from the
properness of the action of $\Ga$ as in \cite[Lemma 3.1] {ParPau11MZ}
that the set
$$
\{[\ga]\in \ga_0^\ZZ\bs\Ga/\ga_1^\ZZ\;:\;h_{\ga_0}(\ga\ga_1\ga^{-1})\le t\}
$$
is finite for all primitive loxodromic elements $\ga_0,\ga_1\in\Ga$
and $t\geq 0$.  This shows that the relative height with respect to a
given primitive loxodromic element is a reasonable complexity within a
given conjugacy class of primitive loxodromic elements of
$\Ga$. Theorem \ref{theo:mainequicountdown} easily implies the
following counting result of primitive loxodromic elements using their
relative heights.

\bcoro\label{coro:complexitycount} If the Bowen-Margulis measure
$m_{\rm BM}$ of $\Ga\bs\hnr$ is finite, then for all primitive
loxodromic elements $\ga_0,\ga_1\in\Ga$, as $t\ra+\infty$,
$$
\sum_{[\ga]\,\in\, \ga_0^\ZZ\bs\Ga/\ga_1^\ZZ\,,\;h_{\ga_0}(\ga\ga_1\ga^{-1})
\,\le\, t} m_{e,\ga} \;\sim\; C_{\ga_0,\ga_1}\,e^{\delta_\Ga t}\,,
$$
where $C_{\ga_0,\ga_1}=\iota_\Ga(\ga_0)\,\iota_\Ga(\ga_1)\,
m_\Ga(\ga_0) \,m_\Ga(\ga_1)\,\|\sigma_{\axis(\ga_0)}\|\,
\|\sigma_{\axis(\ga_1)}\| \,\delta_\Ga^{-1}\,\|m_{\rm BM}\|^{-1}$,
with an error term $\bigO(e^{(n-1)\, t})$ if furthermore $\Ga\bs\hnr$
is compact or arithmetic. If $\Ga\bs\hnr$ has finite volume, then
$$
C_{\ga_0,\ga_1}=\frac{\Vol(\SSS^{n-2})^2\,\ell(\ga_0)\,\ell(\ga_1)}
{2^{n-1}\,(n-1)\,\Vol(\SSS^{n-1})\,\Vol(\Ga\bs \hnr)}\,.
$$
\ecoro

In particular, when $\Ga$ is torsion free, $$\card\{[\ga]\in
\ga_0^\ZZ\bs\Ga/\ga_1^\ZZ:h_{\ga_0}(\ga\ga_1\ga^{-1}) \le t\}\sim
C_{\ga_0,\ga_1}\,e^{\delta_\Ga t}$$ as $t\ra+\infty$.

\medskip
\dem By \cite[Lem.~7]{ParPau11BLMS}, except above finitely many double
classes, the canonical map from $\ga_0^\ZZ\bs\Ga/\ga_1^\ZZ$ to
$\stab_\Ga(\axis(\ga_0))\bs\Ga/\stab_\Ga(\axis(\ga_1))$ is a
$k$-to-$1$ map, where
$$
k=[\stab_\Ga(\axis(\ga_0)):\ga_0^\ZZ]\,[\stab_\Ga(\axis(\ga_1)):\ga_1^\ZZ]
=\iota_\Ga(\ga_0)\,\iota_\Ga(\ga_1)\,m_\Ga(\ga_0)\,m_\Ga(\ga_1)
$$
(using Equation \eqref{eq:miotacardstab}). When $[\ga]$ varies in a
fiber of this map, the quantities $m_{e,\ga}$ (defined in Section
\ref{sect:geometry}) and $h_{\ga_0}(\ga\ga_1\ga^{-1})$ are
constant. As seen in the proof of Theorem \ref{theo:appliarithdim3},
if $\Ga\bs\hnr$ has finite volume, then
$\Vol(\stab_\Ga(\axis(\ga_i))\bs \axis(\ga_i))=
\frac{\ell(\ga_i)}{\iota_G(\ga_i)}$ for $i=0,1$.  The
result then immediately follows from Theorem
\ref{theo:mainequicountdown} and Corollary \ref{coro:constcurvcount}.
\cqfd

\bigskip In the end of this subsection, we give an asymptotic formula
relating the relative height of two primitive loxodromic elements with
the crossratio of its fixed points. 

Recall that the {\em crossratio}
of four pairwise distinct points $a,b,c,d$ in
$\PP_1(\RR)=\RR\cup\{\infty\}$ or in $\PP_1(\CC)= \CC\cup\{\infty\}$
is
$$
[a,b,c,d]=\frac{(c-a)\,(d-b)}{(c-b)\,(d-a)}\,,
$$
with the standard conventions when one of the points is $\infty$.
Using Ahlfors's terminology,
%\todo{his notation is $|a,b,c,d|$} 
the {\em absolute crossratio} of four pairwise distinct points
$a,b,c,d$ in the one-point compactification $\RR^{n-1}\cup\{\infty\}$ of
the standard $n-1$-dimensional Euclidean space is
$$
\llbracket a,b,c,d\rrbracket=
\frac{\|c-a\|\,\|d-b\|}{\|c-b\|\,\|d-a\|}\,,
$$
where $\|\cdot\|$ denotes the standard Euclidean norm and with
conventions analogous to the definition of crossratio when one of the
points is $\infty$.

We will denote by $]a_-,a_+[$ the (oriented) geodesic line in $\hnr$
whose pair of endpoints is a given pair $(a_-,a_+)$ of distinct points
in $\RR^{n-1}\cup\{\infty\}$. In particular, for every loxodromic element
$\ga\in\Ga$, we  have $\axis(\ga)=\;]\ga^-,\ga^+[$.

\blemm\label{lem:complexitycrossratio} (1) Two geodesic lines
$]a_-,a_+[$ and $]b_-,b_+[$ in $\hdr$ or $\htr$ are orthogonal if and
only if $[a_-,a_+,b_-,b_+]=-1$.

\smallskip (2) For all primitive loxodromic elements
$\gamma_0,\ga\in\Ga$, as $h_{\ga_0}(\ga)\to+\infty$, we have
$$
e^{h_{\ga_0}(\ga)}=
\frac4{\llbracket\ga_0^-,\ga^-,\ga_0^+,\ga^+\rrbracket}+\bigO(1)\,.
$$

\elemm

\dem (1) This is a classical fact, see for example pages 15 and 31 of
\cite{Fenchel89}.

\medskip (2) This is a corollary of \cite[Lem.~2.2]{ParPau11MZ}. Note
that the crossratio in \cite{ParPau11MZ} is the logarithm of the
absolute crossratio in this paper.  \cqfd

\subsection{Relative complexity of  quadratic irrationals} 
\label{subsec:relcompquadirrat}

Let $K$ (and its related terminology) be as in the beginning of Section
\ref{sect:irrationals}.  Let $\alpha$ and $\beta$ be quadratic
irrationals in $\wh K$ over $K$.  We define the {\em relative
  height} of $\beta$ with respect to $\alpha$ as
\begin{equation}\label{eq:relhei}
h_{\alpha}(\beta)=
\min\Big\{\frac{|\beta-\alpha|\,|\beta^\sigma-\alpha^\sigma|}
{|\beta-\beta^\sigma|},\;\frac{
|\beta-\alpha^\sigma|\,|\beta^\sigma-\alpha|}
{|\beta-\beta^\sigma|}\Big\}\,.
\end{equation}
Another expression of this complexity using the absolute crossratios
is, if $\beta\notin\{\alpha,\alpha^\sigma\}$,
\begin{equation}\label{eq:relheiwithcrossratio}
h_{\alpha}(\beta)=\frac{|\alpha-\alpha^\sigma|}
{\max\{\llbracket\alpha,\beta,\alpha^\sigma,\beta^\sigma \rrbracket,\;
\llbracket\alpha,\beta^\sigma,\alpha^\sigma,\beta \rrbracket \}}\,.
\end{equation}

Here are some elementary properties of the relative heights.

\blemm\label{lem:proprirelatheight} 
Let $\alpha$ and $\beta$ in $\wh K$ be quadratic irrationals
over $K$. Then

\smallskip\noindent(1) $h_{\alpha^\rho}(\beta^\tau)=
h_{\alpha}(\beta)$ for all $\rho,\tau\in\{\id,\sigma\}$.

\smallskip\noindent(2) $h_{\alpha}(\beta)=0$ if and only if $\beta
\in\{\alpha,\alpha^\sigma\}$.

\smallskip\noindent(3) $h_{\alpha}(\ga \beta)= h_{\alpha}(\beta)$ for
every $\ga\in\stab_{\Ga_K}(\{\alpha,\alpha^\sigma\})$.  
\elemm

\dem We will prove (3), the first two claims being immediate
consequences of the definition.  We denote by $\dot{\;\;}$ the derivative.
Assume first that $\ga\in \Ga_K$ fixes $\alpha$ (and consequently also
$\alpha^\sigma$).  Then, using a well-known identity for linear
fractional transformations (see for example \cite[p.~19]{Ahlfors81})
and the well-known fact that $\dot{\ga}(\alpha)\,
\dot{\ga}(\alpha^\sigma) =1$, we have
\begin{align*}
\frac{|\ga\beta-\alpha|\,|\ga\beta^\sigma-\alpha^\sigma|}
{|\ga\beta-\ga\beta^\sigma|}&=
\frac{|\ga\beta-\ga\alpha|\,|\ga\beta^\sigma-\ga\alpha^\sigma|}
{|\ga\beta-\ga\beta^\sigma|}
\\ & =\sqrt{\frac{|\dot\ga\beta|\,|\dot\ga\alpha|\,
|\dot\ga\beta^\sigma|\,|\dot\ga\alpha^\sigma|}{|\dot\ga\beta|
\,|\dot\ga\beta^\sigma|}}
\, \frac{|\beta-\alpha|\,|\beta^\sigma-\alpha^\sigma|}
{|\beta-\beta^\sigma|}
=\frac{|\beta-\alpha|\,|\beta^\sigma-\alpha^\sigma|}
{|\beta-\beta^\sigma|}\,.
\end{align*}
The invariance by $\ga$ of the second term (seen as a function of
$\beta$) in the definition of the relative height is checked by a
similar computation, replacing $\alpha$ by $\alpha^\sigma$.  If now
$\ga\in \Ga_K$ exchanges $\alpha$ and $\alpha^\sigma$, then a similar
computation gives that $\ga$ exchanges the two terms of the definition
of the relative height. The result follows. \cqfd

\blemm \label{lem:finitrelatheight} For all $\alpha,\beta$ in $\wh
K$ which are quadratic irrationals over $K$ and $s\geq 0$, the set
\begin{equation}\label{eq:finiterelcomplexity}
\{\beta'\in \Ga_{\alpha}\bs\Ga_K\cdot \beta\;:\;h_{\alpha}(\beta')\le s\}
\end{equation}
is finite.
\elemm

\dem Note that this set is well defined by Lemma
\ref{lem:proprirelatheight} (3). By Equation
\eqref{eq:relheiwithcrossratio}, by Lemma
\ref{lem:complexitycrossratio} (2), since $h_{\wh\alpha}(\wh\beta)
=h_{\wh\alpha}(\wh\beta^{-1})$ and since the fixed points of
$\wh\alpha$ and $\wh\beta$ are $\alpha,\alpha^\sigma$ and
$\beta,\beta^\sigma$ respectively, we have
\begin{equation}\label{eq:complexities}
h_{\alpha}(\beta)= \frac{|\alpha-\alpha^\sigma|}4\, 
e^{h_{\wh\alpha}(\wh\beta)}+ \bigO(|\alpha-\alpha^\sigma|)\,.
\end{equation}
Since the set $\{[\ga]\in \wh\alpha\,^\ZZ\bs\Ga/\wh\beta\,^\ZZ\;:\;
h_{\wh\alpha}(\ga\wh\beta\ga^{-1})\le t\}$ is finite for all $t\geq 0$
(see Subsection \ref{subsect:relatcomploxo}), since $\wh\alpha\,^\ZZ$
and $\wh\beta\,^\ZZ$ have finite index in $\Ga_{\alpha}$ and
$\Ga_{\beta}$ respectively, and since
$\ga\wh\beta\ga^{-1}=\wh{\ga\beta}$ for all $\ga\in\Ga_K$, this proves
the result. \cqfd

\medskip
This shows that the relative height with respect to a given
quadratic irrational $\alpha_0$ is, indeed, a reasonable complexity
for quadratic irrationals in a given orbit of the modular group
$\Ga_K$ (modulo the stabiliser of $\alpha_0$).

Recall that for all quadratic irrationals $\alpha$ and $\beta$ in $\wh
K$ over $K$, we have given the expression of two points
$x_{\alpha}^\pm (\beta)$ in $\wh K$ in Equation
\eqref{eq:defixpmalphbet} of the introduction. We will see in the
proof (Lemma \ref{lem:calcboutscomonperp}) of Theorems
\ref{theo:relcomplcount} and \ref{theo:relcomplcountBianchi} that
these points are well defined when $\beta$ varies in a given orbit of
$\Ga_K$, except for finitely many $\Ga_\alpha$-classes. We are now
going to prove the equidistribution of these points, when counted
according to their relative heights. We state separately the cases
$K=\QQ$ and $K$ a quadratic imaginary number field, but it is natural to give a common proof for both statements.

\btheo\label{theo:relcomplcount} Let $K=\QQ$, let $G$ be a subgroup of
finite index in $\Ga_K = \PSLZ$ and let $\alpha_0,\beta_0$ be real
quadratic irrationals over $\QQ$.  Then there exists $\kappa>0$ such
that, as $s\ra+\infty$,
\begin{multline*}
\card\{\beta\in G_{\alpha_0}\bs G\cdot\beta_0\;:\;
h_{\alpha_0}(\beta)\le s\}\\
=\frac{24\,[\Ga_{\alpha_0}:G_{\alpha_0}]\,[\Ga_{\beta_0}:G_{\beta_0}]\,
R_{\alpha_0}\,R_{\beta_0}\,h(\alpha_0)}{\pi^2\, [\Ga_\QQ:G]}\,
\,s\,(1+\bigO(s^{-\kappa}))\,.
\end{multline*}
Furthermore, as $s\ra+\infty$, we have the following convergence of
measures on $\RR-\{\alpha_0,\alpha_0^\sigma\}$:
\begin{equation}\label{eq:appldistribtotgeod}
\frac{\pi^2\,[\Ga_\QQ:G]}
{24\,[\Ga_{\beta_0}:G_{\beta_0}]\,R_{\beta_0}\,s}\;
\sum_{\beta\in G\cdot\beta_0,\;h_{\alpha_0}(\beta)\le s} 
\Delta_{x^-_{\alpha_0}(\beta)}+\Delta_{x^+_{\alpha_0}(\beta)}
\;\weakstar\; 
\frac{d\Leb_\RR(t)}{|Q_{\alpha_0}(t)|}\,.
\end{equation} 
\etheo

\btheo\label{theo:relcomplcountBianchi} Let $K$ be a quadratic
imaginary number field, let $G$ be a subgroup of finite index in
$\Ga_K=\PSLOK$ and let $\alpha_0,\beta_0\in\CC$ be quadratic
irrationals over $K$.  Then there exists $\kappa>0$ such that, as
$s\ra+\infty$,
\begin{multline*}
\card\{\beta\in G_{\alpha_0}\bs G\cdot\beta_0\;:\;
h_{\alpha_0}(\beta)\le s\}\\=
\frac{8\pi^3\,[\Ga_{\alpha_0}:G_{\alpha_0}]\,
[\Ga_{\beta_0}:G_{\beta_0}]\,
\big|\log\big|\frac{\operatorname{tr}\,\wh{\alpha_0}+
    \sqrt{{\operatorname{tr}\wh{\alpha_0}}^2-4}}{2}\big|\big|
\;\big|\log\big|\frac{\operatorname{tr}\,\wh{\beta_0}+
    \sqrt{{\operatorname{tr}\,\wh{\beta_0}\,}^2-4}}{2}\big|\big|}
{m_{\Ga_K}(\alpha_0)\,m_{\Ga_K}(\beta_0)\,h(\alpha_0)^{-2}\,[\Ga_K:G]\,
|D_K|^{\frac 32}\,\zeta_K (2)}\;
\,s^2\,(1+\bigO(s^{-\kappa}))\,.
\end{multline*}
Furthermore, as $s\ra+\infty$, we have the following convergence of
measures on $\CC-\{\alpha_0,\alpha_0^\sigma\}$:
\begin{equation}\label{eq:appldistribtotgeodBianchi}
\frac{m_{\Ga_K}(\beta_0)\,[\Ga_K:G]\,|D_K|^{\frac 32}\,\zeta_K (2)}
{16 \,\pi^2\,[\Ga_{\beta_0}:G_{\beta_0}]\,
\big|\log\big|\frac{\operatorname{tr}\,\wh{\beta_0}+
\sqrt{{\operatorname{tr}\,\wh{\beta_0}\,}^2-4}}{2}\big|\big|\;s^2}\;
\sum_{\beta\in G\cdot\beta_0,\;h_{\alpha_0}(\beta)\le s} 
\Delta_{x^-_{\alpha_0}(\beta)}+\Delta_{x^+_{\alpha_0}(\beta)}
\;\weakstar\; 
\frac{d\Leb_\CC(z)}{|Q_{\alpha_0}(z)|^2}\,.
\end{equation} 
\etheo

\dem Let $r$ and $r'$ be the smallest positive integers such that
$\wh{\alpha_0}^{r}\in G$ and $\wh{\beta_0}^{r'}\in G$ respectively, so
that, as seen at the end of the proof of Theorem
\ref{theo:appliarithdim3} (for $K$ as in Theorem
\ref{theo:relcomplcountBianchi}, but this is clear for $K=\QQ$, since
then $m_{\Ga_K}(\alpha_0)=m_G(\alpha_0)=1$),
\begin{equation}\label{eq:relatindex}
[\Ga_{\alpha_0}:G_{\alpha_0}]=
\frac{m_{\Ga_K}(\alpha_0)}{m_G(\alpha_0)}\,r\,,
\end{equation}
and similarly for $\beta_0$.

To prove the counting claims in the above two theorems, we apply
Corollary \ref{coro:complexitycount} with $n=2$ or $n=3$, $\Ga=G$
(which is arithmetic), $\ga_0=\wh{\alpha_0}^{r}$ and
$\ga_1=\wh{\beta_0}^{r'}$. 

Using in the series of equalities below respectively

$\bullet$~ the bijection $\beta\mapsto [\wh\beta]$ from
$G_{\alpha_0}\bs G\cdot\beta_0$ to $G_{\alpha_0}\bs G/G_{\beta_0}$
(with inverse $[\ga] \mapsto \ga\beta_0$), Equation
\eqref{eq:complexities} and the definition of the complexity
$h(\alpha_0)$,

$\bullet$~ the fact that $[G_{\alpha_0}:
\wh{\alpha_0}^{r\ZZ}] = m_G(\alpha_0)$ and similarly for $\beta_0$, as
in the beginning of the proof of Corollary \ref{coro:complexitycount},

$\bullet$~ Corollary \ref{coro:complexitycount}, noting that since
$\PSLZ$ acts freely on $T^1\HH^2_\RR$, all multiplicities $m_{e,\ga}$
are equal to $1$, if $K=\QQ$, and otherwise, as seen in the proof of
Theorem \ref{theo:appliarithdim3}, the multiplicities $m_{e,\ga}$
different from $1$ contribute only in a negligible way to the sums,

\smallskip
\noindent
we have, as $s\ra +\infty$,
\begin{align*}
  &\card\{\beta\in G_{\alpha_0}\bs G\cdot\beta_0\;:\;
  h_{\alpha_0}(\beta)\le s\}\\=\;& \card\big\{[\ga]\in G_{\alpha_0}\bs
  G/G_{\beta_0}\;:\; h_{\wh{\alpha_0}}(\ga\wh{\beta_0}\ga^{-1})\le
  \ln\big(2\,h(\alpha_0)\,s+\bigO(1)\big)\big\}\\=\;&
  \frac{1+\smallo(1)}{m_G(\alpha_0)\,m_G(\beta_0)}\;\card\big\{[\ga]\in
  \wh{\alpha_0}^{r\ZZ}\bs G/\wh{\beta_0}^{r'\ZZ}\;:\;
  h_{\wh{\alpha_0}}(\ga\wh{\beta_0}\ga^{-1})\le
  \ln\big(2\,h(\alpha_0)\,s+\bigO(1)\big)\big\}\\=\;&
\frac{\Vol(\SSS^{n-2})^2\;r\,\ell(\wh{\alpha_0})\,r'\,\ell(\wh{\beta_0})}
{m_G(\alpha_0)\,m_G(\beta_0)\,2^{n-1}\,(n-1)\,
\Vol(\SSS^{n-1})\,\Vol(G\bs\HH^n_\RR)}\;
\big(2\,h(\alpha_0)\,s\big)^{n-1}(1+\bigO(s^{-\kappa}))\,.
\end{align*}

When $K=\QQ$ and $n=2$, this proves the counting claim in Theorem
\ref{theo:relcomplcount}, by using Equations \eqref{eq:volumehdrmodG},
\eqref{eq:translengthregulator} and \eqref{eq:relatindex}, and noting
that $m_{\Ga_K}(\alpha_0)=m_{\Ga_K}(\beta_0)=1$.

When $K$ is quadratic imaginary and $n=3$, this proves the counting
claim in Theorem \ref{theo:relcomplcountBianchi}, by using Equations
\eqref{eq:covolBianchi}, \eqref{eq:formultranslength} and
\eqref{eq:relatindex}.

\bigskip To prove the equidistribution claims in the above two
theorems, we need to use a stronger geometric equidistribution theorem
(also proven in \cite{ParPau14}) than the one stated in Section
\ref{sect:geometry}.

For any $\beta\in G\cdot\beta_0$, let $v_{\alpha_0}(\beta)$ be the
initial tangent vector of the common perpendicular from the geodesic
line $D^-=\axis(\wh{\alpha_0})$ with endpoints $\alpha_0$ and
$\alpha_0^\sigma$ to that with endpoints $\beta$ and $\beta^\sigma$,
if it exists (and the summations below are on $\beta\in G\cdot\beta_0$
such that it does exist). Let $D^+=\axis(\wh{\beta_0})$ and $G_{D^+}$
its stabiliser in $G$, so that as already seen $\Vol(G_{D^+}\bs
D^+)=\frac{r'\,\ell(\wh{\beta_0})}{\iota_{G}(\beta_0)}$. By replacing
the equidistribution of initial points by the one of the initial
tangent vectors (see \cite[Coro.~20]{ParPau14}), and by arguments
similar to the ones leading to Equations \eqref{eq:equidquadirratQ}
and \eqref{eq:equidquadirratK}, we have, by the expression of the
skinning measure and Bowen-Margulis measure recalled in Section
\ref{sect:geometry}, with $n=2$ if $K=\QQ$ and $n=3$ otherwise, as
$t\ra+\infty$,
$$
\frac{2^{n-1}\,(n-1)\,m_G(\beta_0)\,\Vol(\SSS^{n-1})\;\Vol(G\bs \hnr)}
{\iota_{G}(\beta_0)\,\Vol(\SSS^{n-2})\,\Vol(G_{D^+}\bs D^+)\,e^{(n+1)t}}\;
\;\sum_{\beta\in G\cdot\beta_0,\;h_{\wh{\alpha_0}}(\wh\beta)\le t} 
\Delta_{v_{\alpha_0}(\beta)}
\;\weakstar\;\Vol_{\normalout D^-}\,.
$$
Hence, by Equation \eqref{eq:complexities}, as $s\ra+\infty$,
\begin{equation}\label{eq:equidistribnormalout}
\frac{(n-1)\,m_G(\beta_0)\,\Vol(\SSS^{n-1})\;\Vol(G\bs \hnr)}
{\Vol(\SSS^{n-2})\,r'\,\ell(\wh{\beta_0})\,
h(\alpha_0)^{n-1}\;s^{n-1}}\;
\;\sum_{\beta\in G\cdot\beta_0,\;h_{\alpha_0}(\beta)\le s} 
\Delta_{v_{\alpha_0}(\beta)}
\;\weakstar\;\Vol_{\normalout D^-}\,.
\end{equation}

We give in the following two lemmas the main computations used to
deduce from this the equidistribution claims in the above two theorems.

\blemm \label{lem:calcboutscomonperp} If the geodesic lines
$]\alpha_0,\alpha_0^\sigma[$ and $]\beta,\beta^\sigma[$ in $\hnr$ are
disjoint, then the endpoints of the geodesic line containing their
common perpendicular are the points $x^\pm_{\alpha_0}(\beta)$ given by
Equation \eqref{eq:defixpmalphbet}, which are hence well defined.
\elemm

\dem Let $]x^-,x^+[$ be this geodesic line. By Lemma
\ref{lem:complexitycrossratio} (1), and since two nonintersecting
geodesic lines have one and only one geodesic line orthogonal to both
of them, the points $x^-$ and $x^+$ are the two solutions of the pair
of equations $[x^-,x^+,\alpha_0,\alpha_0^\sigma] =-1
=[x^-,x^+,\beta,\beta^\sigma]$. An easy computation shows that this
system is equivalent to
$$
\begin{cases}
  x^-+x^+={\displaystyle 2\,\frac{\n(\alpha_0)-\n(\beta)}
    {\tr\,\alpha_0-\tr\,\beta}}\\
  \phantom{W}\\
  x^-x^+={\displaystyle\frac{\tr\,\beta\,\n(\alpha_0)-\tr\,
\alpha_0\,\n(\beta)}{\tr\,\alpha_0-\tr\,\beta}}\,.
\end{cases}
$$
Solving a quadratic equation gives the result.  
\cqfd

\blemm \label{lem:calcmeaspushskinn} For $n=2,3, 5$, for all distinct
$x,y$ in $\partial_\infty\HH^n_\RR-\{\infty\}$, the pushforward of the
Riemannian measure of the unit normal bundle of the geodesic line
$L=\;]x,y[$ in $\HH^n_\RR$ by the positive endpoint map $v\mapsto z=v_+$
from $\normalout L$ to $\partial_\infty\HH^n_\RR-\{x,y\}$ is
$$
\frac{\|x-y\|^{n-1}}{\|(z-x)(z-y)\|^{n-1}}\;
d\Leb_{\RR^{n-1}}(z)\,.
$$
\elemm

By symmetry, the same result holds for the pushforward by the negative
endpoint map $v\mapsto z=v_-$. When $n=2,3$, $x=\alpha_0$ and $y=
\alpha_0^\sigma$, this measure is $\frac{2^{n-1}\,d\Leb_{\wh K} (z)}
{h(\alpha_0)^{n-1}\; |Q_{\alpha_0}(z)|^{n-1}}$, with the previous notation.

\medskip \dem Let $\wh{K}=\RR$ if $n=2$, $\wh{K}=\CC$ if $n=3$ and
$\wh{K}=\HH$ if $n=5$ (we denote $|z|=\redn(z)^{1/2}$ if $z\in\HH$),
so that $\partial_\infty\HH^n_\RR-\{\infty\} =\wh K$.

If $L_\infty$ is the geodesic line with points at infinity $0\in
\wh{K}$ and $\infty$, then for all $v\in \normalout L_\infty$, we have
$d\Vol_{\normalout L_\infty}(v)= \frac{dt}{t} d\Vol_{\SSS^{n-2}}
(\sigma)$, where $t>0$ is the last coordinate of $\pi(v)$ and
$\sigma\in\SSS^{n-2}$ is the parameter of the Euclidean horizontal
vector $v$ (note that $d\Vol_{\SSS^{0}}$ is the counting measure on
$\{-1,1\}$ if $n=2$). If $z=v_+\in\wh{K}- \{0\}$ is the positive point
at infinity of $v$, then $z=t\,\sigma$ and
$$
d\Vol_{\normalout L_\infty}(v)=\frac{d\Leb_{\wh K}(z)}{|z|^{n-1}}\;.
$$

We may assume that $x>y$ if $n=2$. The homography $A:z\mapsto
(z-x)(z-y)^{-1}$ maps $x$ to $0$ and $y$ to $\infty$. Since it may be
written $z\mapsto(\lambda z-\lambda x)(\lambda z-\lambda y)^{-1}$ with
$\lambda$ in (the center of) $\wh K$ such that $\lambda^2=(x-y)^{-1}$ if
$\wh K=\RR,\CC$ and $\lambda=\frac{1}{\redn(y-x)}$ if $\wh K=\HH$ (see
Equation \eqref{eq:detdieud} to check the Dieudonné determinant is
$1$), it is the extension at infinity of an isometry of $\HH^n_\RR$
(by Poincaré's extension), sending $\Vol_{\normalout L}$ to
$\Vol_{\normalout L_\infty}$. Since $|\dot A\,(z)|=
\frac{|x-y|\;\,}{|z-y|^2}$, we have
$$
(A^{-1})_*\Big(\frac{d\Leb_{\wh K}(z)}{|z|^{n-1}}\Big)
=\frac{|x-y|^{n-1}}{|(z-x)(z-y)|^{n-1}}\;d\Leb_{\wh K}(z)\;.
$$
The result follows. 
\cqfd

\medskip By the continuity of the pushforward maps of measures, by
summing the pushforward of Equation \eqref{eq:equidistribnormalout}
both by the positive and by the negative endpoint map, by the comments
following Lemma \ref{lem:calcmeaspushskinn}, and since
$\{x^-_{\alpha_0} (\beta), x^+_{\alpha_0}(\beta)\}= \{v_{\alpha_0}
(\beta)_-, v_{\alpha_0}(\beta)_+\}$ by Lemma
\ref{lem:calcboutscomonperp}, we have, as $s\to+\infty$,
$$
\frac{(n-1)\,m_G(\beta_0)\,\Vol(\SSS^{n-1})\;\Vol(G\bs \hnr)}
{2^{n-1}\,\Vol(\SSS^{n-2})\,r'\,\ell(\wh{\beta_0})\;s^{n-1}}\, 
\sum_{\beta\in  G\cdot\beta_0,\;h_{\alpha_0}(\beta)\le s}
\Delta_{x^-_{\alpha_0}(\beta)} + \Delta_{x^+_{\alpha_0}(\beta)}
\;\weakstar\;
\frac {2\,d\Leb_{\wh K}(z)}{|Q_{\alpha_0}(z)|^{n-1}}\,.
$$ 
The equidistribution claims in Theorems \ref{theo:relcomplcount} and
\ref{theo:relcomplcountBianchi} follow, as in the end of the proof of
their counting claims.  
\cqfd
%Note that in the proof of the counting result in Theorem
%\ref{theo:relcomplcount} we used the fact that no nonidentity
%orientation-preserving isometry of $\hnr$ fixes a geodesic line
%pointwise. In higher dimension, this is no longer necessarily the case
%but it has no effect on the count as we count modulo $\Ga_{\alpha_0}$
%and the possible rotated copies are counted separately.

\medskip \rem If $\wh K=\HH$, the crossratio of a quadruple of
pairwise distinct points $(a,b,c,d)$ in $\wh K$ may be defined as
$$
[a,b,c,d]=(c-b)^{-1}(c-a)(d-a)^{-1}(d-b)\;.
$$
Though this crossratio does not extend continuously to the quadruples
of pairwise distinct points in the one-point compactification
$\HH\cup\{0\}$, and it is only invariant under $\SLH$ up to conjugation,
it also caracterises the orthogonality of the geodesic lines in $\hcr$
by the same formula (two geodesic lines $]a_-,a_+[$ and $]b_-,b_+[$ in
$\hcr$ are orthogonal if and only if $[a_-,a_+,b_-,b_+] = -1$). Furthermore,
its absolute crossratio $|c-b|^{-1}|c-a||d-a|^{-1}|d-b|$ (we denote
$|z|=\redn(z)^{1/2}$ if $z\in\HH$) does extend continuously (as
already mentionned) and is invariant under $\SLH$ (see for instance
\cite{DelRet13} for all this).

Let $A$ be a definite quaternion algebra over $\QQ$, let $\OOO$ be a
maximal order in $A$, let $G$ be a finite index subgroup of $\PSLO$,
and let $\alpha_0,\beta_0\in\HH$ be two loxodromic quadratic
irrationals over $\OOO$ (see the definitions above Theorem
\ref{theo:appliarithdim5}). We may also define the {\it relative
  height} of an element $\beta$ of $\PSLO\cdot \beta_0$ with respect
to $\alpha_0$ by the same formula \eqref{eq:relhei}. It is still
related to absolute crossratios by Equation
\eqref{eq:relheiwithcrossratio}, and still satisfies the properties of
Lemma \ref{lem:proprirelatheight}. We do have a counting asymptotic of
the elements of $\beta\in G_{\alpha_0}\bs G\cdot \beta_0$ wih relative
height at most $s$ with respect to $\alpha_0$, as $s\ra+\infty$,
analogous to Theorem \ref{theo:relcomplcount} and
\ref{theo:relcomplcountBianchi}, and a related equidistribution
result, but now, due to the higher difficulty of solving quadratic
equations in the noncommutative $\HH$, the points
$x^\pm_{\alpha_0}(\beta)$ are much less explicit: they can only be
defined as the unique two solutions $x^-,x^+$ of the pair of equations
$[x^-,x^+,\alpha_0,\alpha_0^\sigma] =-1
=[x^-,x^+,\beta,\beta^\sigma]$. We leave the precise statement to the
interested reader.

\subsection{Counting crossratios and application to 
generalised Schottky-Klein functions} 
\label{sect:countingcrossratios}

The aim of this subsection is to generalise the work of Pollicott
\cite{Pollicott11}, using the asymptotic of crossratios to prove the
convergence of Schottky-Klein functions, to a much larger class of
Kleinian groups.

Let $\wt M$ be a complete simply connected Riemannian manifold with
(dimension at least $2$ and) pinched negative sectional curvature
$-b^2\le K\le -1$. Let $\Ga$ be a torsion free (only to simplify the
statements) nonelementary (that is, not virtually nilpotent) discrete
group of isometries of $\wt M$, with critical exponent $\delta_\Ga$
and limit set $\Lambda\Ga$ (defined as in Section
\ref{sect:geometry}). Let $D^-$ and $D^+$ be two geodesic lines in
$\wt M$ with points at infinity in the {\it domain of discontinuity}
$\Omega\Ga=\partial_\infty\wt M-\Lambda\Ga$ of $\Ga$.

The Bowen-Margulis measure $m_{\rm BM}$ and skinning measures
$\sigma_{D^\pm}$ are defined as in Section \ref{sect:geometry} where
$\wt M=\HH^n_\RR$. We refer to \cite{DalOtaPei00} for finiteness
criteria of $m_{\rm BM}$ (for instance satisfied if $\Ga$ is
convex-cocompact, or when $\Ga$ is geometrically finite in constant
curvature (in which case $\delta_\Ga$ is the Hausdorff dimension of
$\Lambda\Ga$), but there are many more examples), and
\cite{Babillot02b} for mixing criteria of $m_{\rm BM}$ under the
geodesic flow (for instance satisfied, when $m_{\rm BM}$ is finite, if
$\wt M$ has constant curvature, but there are many more examples).

The following result generalises \cite[Theo.~1.4]{Pollicott11} (whose
proof used symbolic dynamics and transfer operator techniques), which
was only stated for $\wt M=\HH^3_\RR$, $\Ga$ a Schottky group and
$D^-=D^+$, with a nonexplicit multiplicative constant.

\btheo\label{theo:genepol}
 If $m_{\rm BM}$ is finite and mixing, then, as $s\ra+\infty$,
$$
\card\{\ga\in\Ga\;:\; d(D^-,\ga\,D^+)\leq s\}\sim
\frac{\|\sigma_{D^-}\|\,\|\sigma_{D^+}\|}
{\delta_\Ga\,\|m_{\rm BM}\|}\;e^{\delta_\Ga\,s}\,.
$$
\etheo

\dem Since the endpoints of $D^\pm$ are in the domain of discontinuity
of $\ga$, the support of the measure $\wt \sigma_{D^\pm}$ is compact,
hence the skinning measure $\sigma_{D^\pm}$ is finite. The result then
follows from \cite[Coro.~20]{ParPau14}, since the stabiliser of
$D^\pm$ is trivial. 
\cqfd

\medskip Until the end of Subsection \ref{sect:countingcrossratios}, we
assume that $\wt M=\HH^3_\RR$. The next result on asymptotic
properties of crossratios is deduced from Theorem \ref{theo:genepol}
as Theorem 1.3 is deduced from Theorem 1.4 in \cite{Pollicott11}.

\bcoro
If $m_{\rm BM}$ is finite, then for all $z,\xi$ in $\Omega\Ga$,
as $\epsilon\ra 0$,
 $$
 \card\{\ga\in\Ga\;:\; 
\big|\,[z,\xi,\ga z,\ga \xi]-1\big|\geq \epsilon\}\sim
 \frac{4^{\delta_\Ga}\|\sigma_{D^-}\|\,\|\sigma_{D^+}\|}
 {\delta_\Ga\,\|m_{\rm BM}\|}\;\epsilon^{-\delta_\Ga}\,.
$$
\ecoro

The next result is proven as \cite[Coro.~1.2]{Pollicott11}, and
defines a much larger class of holomorphic maps in two variables than
the Schottky-Klein functions (see loc.~cit.~for motivations and
references). We denote by $\Ga_0$ a choice of a representative modulo
inverse of the nontrivial elements of $\Ga$.

\bcoro If $m_{\rm BM}$ is finite and $\delta_\Ga<1$, then the function
$$
(z,\xi)\mapsto (z-\xi) \;\prod_{g\in\Ga_0}\,
\frac{(gz-z)(g\xi- z)}{(g\xi-\xi)(g\xi- z)}\,,
$$
where the terms in the product are ordered by the distance of their
absolute value to $1$, converges uniformly on compact subsets of
$\Omega\Ga\times\Omega\Ga$.  
\ecoro

\section{Counting and equidistribution around binary 
quadratic,  Hermitian and Hamiltonian forms}
\label{sect:errorterms}

Let $Q(x,y)=a\,x^2+b\,xy+c\,y^2$ be a binary quadratic form, which is
primitive integral (its coefficients $a,b,c$ are relatively prime
elements in $\ZZ$). Let $\discr_Q=b^2-4ac$ be the
discriminant of $Q$.  The group $\SLZ$
acts on the set of integral binary quadratic forms by precomposition,
and we denote by $\operatorname{SO}_Q(\ZZ)$ the stabiliser of $Q$.
When $Q$ is positive definite, the following asymptotic is originally due
to Gauss with a weaker error bound:
$$
\card\{(u,v)\in\ZZ\times\ZZ:(u,v)=1,\; Q(u,v)\le s\}=
\frac{12}{\pi\sqrt{-\discr_Q}}s+\bigO(s^{131/416})\,,
$$
the current error bound is proved in \cite{Huxley03}. This
asymptotic, with a less explicit error bound, can be given a geometric
proof using $\PSLZ$-orbits of horoballs and points in $\hdr$ in the
same way as Theorem \ref{theo:posdefHermerror} in Subsection
\ref{sect:posdef} below.

When $Q$ is indefinite, the group of automorphs $\operatorname{SO}_Q
(\ZZ)$ is infinite, and it is appropriate to count
$\operatorname{SO}_Q (\ZZ)$-orbits of relatively prime representations
of integers by $Q$.  Let $$R_Q=\ln\frac{t_Q+u_Q\sqrt{\discr_Q}}{2}$$ 
 be the regulator of $Q$, 
where
$(t_Q,u_Q)$ is the fundamental solution of the Pell-Fermat equation
$t^2-\discr_Q\,u^2=4$. The asymptotic behavior of this
counting function is (see for example \cite[p. 164]{Cohn62})
\begin{equation}\label{eq:quadrformasymptotic}
\card\big( \operatorname{SO}_Q(\ZZ)\backslash 
\{(u,v)\in\ZZ\times\ZZ:(u,v)=1,\; |Q(u,v)|\le s\}\big)=
\frac{12\,R_Q}{\pi^2\sqrt{\discr_Q}}s+\bigO(\sqrt s)\,.
\end{equation}
A geometric proof of the asymptotic \eqref{eq:quadrformasymptotic}
with a less explicit error bound can be obtained from the proof of
\cite[Th\'eo. 4.5]{ParPau12JMD}, using Corollary
\ref{coro:constcurvcount} instead of \cite[Cor. 3.9] {ParPau12JMD}.

The results in Subsection \ref{subsect:errortermintrepbybinary},
generalising the asymptotic \eqref{eq:quadrformasymptotic} to binary
Hermitian and Hamiltonian forms, are the versions with error terms and
with related equidistribution results of respectively Theorem 4,
Theorem 1, Corollary 3 in \cite{ParPau11BLMS}, and Theorem 13,
Corollary 18, Theorem 1, Corollary 2 in \cite{ParPau13ANT}. Similarly,
the error term in Corollary 6 of \cite{ParPau11BLMS} is
$\operatorname{O}(s^{2-\kappa})$ for some $\kappa>0$.

Given an arithmetic group $\Ga$ of isometries of a totally geodesic
subspace $X_\Ga$ of $\HH^n_\RR$, we will denote by $\covol(\Ga)$ the
volume of the Riemannian orbifold $\Ga\bs X_\Ga$.

\subsection{Indefinite binary Hermitian and Hamiltonian 
forms}
\label{subsect:errortermintrepbybinary}

Let $K$ and $\mmm$ be as in the beginning of Subsection
\ref{subsect:mertensbianchi}.  Fix a binary Hermitian form
$f:\CC^2\ra\RR$ with
\begin{equation}\label{eq:form}
f(u,v)=a\,|u|^2+2\,\Re(b\,\ov u\, v) +c\,|v|^2
\end{equation}
which is integral over $\OOO_K$ (its coefficients $a=a(f),b=b(f),
c=c(f)$ satisfy $a,c\in\ZZ$ and $b\in\OOO_K$).  The group $\SLC$ acts
on the right on the binary Hermitian forms by precomposition, and the
action of the Bianchi subgroup $\Ga_K=\SLOK$ preserves the integral
ones.  The group $\autom$ of automorphs of $f$ consists of the
elements $g\in\Ga_K$ stabilising $f$, that is, such that $f\circ g=f$.

In this subsection, we assume that $f$ is indefinite (its discriminant
$\discr(f)=|b|^2-ac$ is positive), a feature which is preserved by the
action of $\SLC$.

Let $G$ be a finite index subgroup of $\Ga_K$. Let $\Ga_{K,\,x,\,y}$
and $G_{x,\,y}$ be the stabilisers of $(x,y)\in K\times K$ in $\Ga_K$ and
$G$ respectively. Let $\iota_G=1$ if $-\id\in G$ and $\iota_G=2$
otherwise. For every $x,y$ in $\OOO_K$ not both zero, and for every
$s>0$, let
\begin{multline*}
\psi_{f,\,G,\,x,\,y}(s)=\\ \card\;\;_{\mbox{$\autom\cap G$}}\bs
\big\{(u,v)\in G(x,y)\;:\; \n(\OOO_Kx+\OOO_Ky)^{-1}|f(u,v)|\leq
s\big\}\;.
\end{multline*}

\btheo\label{theo:mainversionGerror} There exists $\kappa>0$ such
that, as $s$ tends to $+\infty$,
$$
\psi_{f,\,G,\,x,\,y}(s)\sim
\frac{\pi\;\iota_G\;[\Ga_{K,\,x,\,y}:G_{x,\,y}]\;\covol(\autom\cap G)}
{2\;|D_K|\;\zeta_K(2)\;\discr(f)\;[\Ga_K:G]}\;\;\;
s^2 (1+\operatorname{O}(s^{-\kappa}))\;.
$$
If $y\neq 0$, then, as $s\ra+\infty$,
$$
\frac{[\Ga_K:G]\,|D_K|^{\frac{3}2}\,\zeta_K(2)\,\discr(f)\,\n(y)}
{\pi^2\,\iota_G\,\covol(\autom\cap G)\,s^2}
\;\sum_{f'\in\, G\cdot f,\;0<|f'(x,\,y)|\,\leq\, s}
\Delta_{\frac{a(f')\,x\,+\,b(f')\,y}{f'(x,\,y)}}
\;\weakstar\; \Leb_{\CC}\,.
$$
If $\Lambda$ is the lattice of $\lambda\in\OOO_K$ such that
$\pm\Big(\begin{array}{cc} 1&\lambda \\0&1\end{array} \Big)\in G$,
identifying $\OOO_K$ with the upper triangular unipotent sugbroup of
$\Ga_K$, we have, as $s\ra+\infty$,
\begin{multline*}
\card\{f'\in\, \Lambda\bs G\cdot f\;:\;0<|a(f')|\leq s\}\\
= \frac{\pi^2\,\iota_G\,[\OOO_K:\Lambda]\,\covol(\autom\cap G)}
{2\,[\Ga_K:G]\,|D_K|^{\frac{1}2}\,\zeta_K(2)\,\discr(f)}
\;\;
s^2 (1+\operatorname{O}(s^{-\kappa}))\;.
\end{multline*}
\etheo

For smooth functions $\psi$ with compact support on $\CC$, there is an
error term in the above equidistribution claim, of the form
$\bigO(s^{2-\kappa} \|\psi\|_\ell)$ where $\kappa>0$ and
$\|\psi\|_\ell$ is the Sobolev norm of $\psi$ for some $\ell\in\NN$.

\medskip
\dem The first claim is proven in the same way as
\cite[Thm. 4]{ParPau11BLMS}, using Corollary \ref{coro:constcurvcount}
instead of \cite[Cor. 4.9 ]{ParPau12JMD}.

Let us prove the second claim. Let $n=3$ and $\wh K=\CC$. We denote by
$\ov{H}$ the image in $\PSL_2(\wh K)$ of any subgroup $H$ of
$\SL_2(\wh K)$. As in the proof of Theorem \ref{theo:appliHerm}, we
consider the isometric action of $\PSL_2(\wh K)$ on (the upper
halfspace model of) $\hnr$.  As explained in \cite[\S 2]{ParPau11BLMS}
(where the convention for binary Hermitian form was to replace $b$ by
$\overline{b}$), when $a(f)\neq 0$, the totally geodesic hyperplane
$$
\C(f)=\{(z,t)\in\hnr\;:\;f(z,1)+|a(f)|\,t^2=0\}
$$
is preserved by $\autom$, its boundary at infinity is the
$(n-2)$-sphere of center $-\frac{b(f)}{a(f)}$ and radius
$\frac{\sqrt{\discr(f)}} {|a(f)|}$, the hyperbolic orbifold
$\autom\bs\C(f)$ has finite volume (and so does $(\autom\cap
G)\bs\C(f)$~), and, for every $g\in\SL_2(\wh K)$,
$$
\C(f\circ g) = g^{-1}\, \C(f)\,.
$$

Let $\rho=x\,y^{-1}\in K\cup\{\infty\}$, $\ga_\rho=
\begin{pmatrix}\rho& -1\\1& \;\,0\end{pmatrix}\in\SL_2(K)$ if $y\neq
0$ and $\ga_\rho=\id$ otherwise.  We define $R_G(f)=2$
%if $f$ is $G$-reciprocal, that is,
if there exists an element $g$ in $G$ such that $f\circ g=-f$, and
$R_G(f)=1$ otherwise.

To prove the equidistribution claim of Theorem
\ref{theo:mainversionGerror}, we apply Equation
\eqref{eq:distribhorob} with $n=3$, $\Ga$ the arithmetic group
$\ga_\rho^{-1} \,\overline{G}\,\ga_\rho$, $D^-$ the horoball
$\H_\infty=\{(z,t)\in\hnr\;:\;t\geq 1\}$ (noting that $\infty$ is
indeed a parabolic fixed point of $\Ga$, since $\rho=\ga_\rho\infty$
is a parabolic fixed point of $\overline{G}$) and $D^+$ the totally
geodesic hyperplane $\C(f\circ\ga_\rho)=\ga_\rho^{-1}\,\C(f)$. The
stabiliser $\Ga_{D^+}$ of $D^+$ in $\Ga$ indeed has finite covolume in
$D^+$, since (by conjugation for the first equality and by
\cite[Eq.~(8)]{ParPau11BLMS} for the second one)
\begin{equation}\label{eq:volDpmodgaDp}
\Vol(\Ga_{D^+}\bs D^+)=
\Vol\big(\stab_{\,\overline{G}\,}(\C(f))\bs \C(f)\big)=
\frac{1}{R_G(f)}\;\covol(\autom \cap G)\,.
\end{equation}

For every $\ga\in\Ga$, the origin of the common perpendicular from
$D^-$ to $\ga D^+=\C(f\circ\ga_\rho\circ\ga)$ (when it exists, that is
when $a(f\circ\ga_\rho\circ\ga)\neq 0$ and the radius of the circle at
infinity of $\C(f\circ\ga_\rho\circ\ga)$ is strictly less than $1$),
is $\big(- \frac{b(f\circ\ga_\rho\circ\ga)} {a(f\circ\ga_\rho\circ
  \ga)}, 1\big)$, and its hyperbolic length is $\big|\ln\big(\frac
{\sqrt{\discr(f)}} {|a(f\circ\ga_\rho\circ\ga)|} \big)\big|=
\ln\big(\frac {|a(f\circ\ga_\rho\circ\ga)|} {\sqrt{\discr(f)}} \big)$
(by invariance of the discriminant under $\SL_2(\wh K)$). Since $D^+$
has codimension $1$ and $\PSL_2(\wh K)$ preserves the orientation of
$\hnr$, the pointwise stabiliser of $D^+$ in $\Ga$ is trivial. Using a
similar comment about the multiplicities as in the proof of Theorem
\ref{theo:appliarithdim3}, we hence have by taking
$t=\ln\frac{s}{\sqrt{\discr(f)}}$ in Equation \eqref{eq:distribhorob},
as $s\ra+\infty$,
\begin{multline*}
\frac{(n-1)\,\Vol(\SSS^{n-1})\,\Vol(\Ga\bs\HH^n_\RR)}
{\Vol(\SSS^{n-2})\,\Vol(\Ga_{D^+}\bs D^+)\,
\big(\frac{s}{\sqrt{\discr(f)}}\big)^{n-1}}
\;\sum_{\ga\in\Ga/\Ga_{D^+},\;0<|a(f\circ\ga_\rho\circ\ga)|\leq s}
\Delta_{\big(-\frac{b(f\circ\ga_\rho\circ\ga)}
{a(f\circ\ga_\rho\circ\ga)},\,1\big)}\\\;\weakstar \;\Vol_{\partial D^-}\,.
\end{multline*}
Since the map from $G\cdot f$ to $\Ga/\Ga_{D^+}$ defined by $f\circ
\ga\mapsto \ga_\rho^{-1}\ga\ga_\rho\Ga_{D^+}$ is a $R_G(f)$-to-$1$
map, and since the pushforward of measures by the continuous map
$(z,1)\mapsto -\,z$ from $\partial D^-$ to $\wh K$ sends
$\Vol_{\partial D^-}$ to $\Leb_{\wh K}$, we have
\begin{equation}\label{eq:quasifinherm}
\frac{(n-1)\,\Vol(\SSS^{n-1})\,\Vol(\Ga\bs\hnr)\,\discr(f)^{n-1}}
{\Vol(\SSS^{n-2})\,\Vol(\Ga_{D^+}\bs D^+)\,R_G(f)\,s^{n-1}}
\;\sum_{f'\in G\cdot f,\;0<|a(f'\circ\ga_\rho)|\leq s}
\Delta_{\frac{b(f'\circ\ga_\rho)}{a(f'\circ\ga_\rho)}}
\;\weakstar \;\Leb_{\wh K}\,.
\end{equation}
Using Humbert's formula \eqref{eq:covolBianchi}, we have
\begin{equation}\label{eq:modifcovolBianchi}
\Vol(\Ga\bs\htr)=\Vol(\,\overline{G}\,\bs\htr)=
[\,\overline{\Ga_K}:\overline{G}\,] 
\Vol(\,\overline{\Ga_K}\,\bs\htr) =\frac{[\Ga_K:G]}
{\iota_G}\,\frac{|D_K|^{\frac 32}\zeta_K (2)}{4\pi^2}\,.
\end{equation}
By Equations \eqref{eq:quasifinherm},
\eqref{eq:volDpmodgaDp} and \eqref{eq:modifcovolBianchi}, we have
$$
\frac{[\Ga_K:G]\,|D_K|^{\frac{3}2}\,\zeta_K(2)\,\discr(f)}
{\iota_G\,\pi^2\,\covol(\autom\cap G)\,s^2}
\;\sum_{f'\in G\cdot f,\;0<|a(f'\circ\ga_{\rho})|\,\leq\, s}
\Delta_{\frac{b(f'\circ\ga_{\rho})}{a(f'\circ\ga_{\rho})}}
\;\weakstar\; \Leb_{\CC}\,.
$$

When $y=0$, this proves Theorem \ref{theo:mainHermerrorintro} in the
Introduction. By considering the measures induced on the compact
quotient $\Lambda\bs \CC$, whose volume is
$\frac{|D_K|}{2}[\OOO_K:\Lambda]$, and by applying Theorem
\ref{theo:mainHermerrorintro} with error term to the constant function
$1$, we get the last claim of Theorem \ref{theo:mainversionGerror}.

When $y\neq 0$, replacing $s$ by $\frac{s}{|y|^2}$, using the
pushforward of measures on $\CC$ by $z\mapsto -\,\overline{y^{-1}\,
z}$, which sends $\Leb_\CC$ to $|y|^2\Leb_\CC$, and since
$a(f'\circ\ga_\rho)=f'\circ\ga_\rho(1,0)= (y\,\ov y)^{-1} \,f'(x,y)$
and $b(f'\circ\ga_\rho)= -\big(a(f')\,\overline{\rho} +\overline{b(f')}
\,\big)$ by an easy computation, this gives the equidistribution claim
in Theorem \ref{theo:mainversionGerror}.  
\cqfd

\medskip 
For every $s>0$, we consider the integer (depending only on
the ideal class of $\mmm$)
$$
\psi_{f,\,\mmm}(s)=\card\;\;_{\mbox{$\autom$}}\bs
\big\{(u,v)\in\mmm\times\mmm\;:\;\n(\mmm)^{-1}|f(u,v)|\leq s,
\;\;\;u\OOO_K+v\OOO_K=\mmm\big\}\;,
$$
which is the number of nonequivalent $\mmm$-primitive representations
by $f$ of rational integers with absolute value at most $s\,\n(\mmm)$. 

By taking $G=\Ga_K$ and $x,y\in K$ such that $\mmm=\OOO_Kx+\OOO_Ky$
(which exists, see for instance \cite[\S 7]{ElsGruMen98}) in the first
claim of Theorem \ref{theo:mainversionGerror}, we have the following
result.

\btheo \label{theo:mainHermerror}
There exists $\kappa>0$ such that, as $s$ tends to $+\infty$, 
$$
\psi_{f,\,\mmm}(s)\;\sim\; \frac{\pi\;\covol(\autom)}
{2\;|D_K|\;\zeta_K(2)\;\discr(f)}\;\;\;
s^2 (1+\operatorname{O}(s^{-\kappa}))\;.\;\;\;\Box
$$
\etheo

Using the results of \cite{MacRei91} (or other ways of obtaining
formulas for $\covol(\autom)$), one gets very explicit versions of
Theorem \ref{theo:mainHermerror} in special cases. A constant
$\const(f)\in\{1,2,3,6\}$ is defined as follows.  If $\discr(f)\equiv
0\mod 4$, let $\const(f)=2$. If the coefficients $a$ and $c$ of the
form $f$ as in Equation \eqref{eq:form} are both even, let
$\const(f)=3$ if $\discr(f)\equiv 1\mod 4$, and let $\const(f)$ be the
remainder modulo $8$ of $\discr(f)$ if $\discr(f)\equiv 2 \mod 4$.  In
all other cases, let $\const(f)=1$. The following result on integral
binary Hermitian forms $f$ over $\QQ(i)$ follows from the first claim
of Theorem \ref{theo:mainHermerror} using $K=\QQ(i)$ as in the proof
of \cite[Coro.~3]{ParPau11BLMS}.

\bcoro \label{coro:gaussianerror} 
There exists $\kappa>0$ such that, as $s$ tends to $+\infty$,
\begin{multline*}
\card\;\;_{\mbox{$\operatorname{SU}_f(\ZZ[i])$}}\bs
\big\{(u,v)\in\ZZ[i]^2\;:\;u\,\ZZ[i]+v\,\ZZ[i]\in\ZZ[i],
\;\;|f(u,v)|\leq s\big\}\\=
\frac{\pi^2}{8\;\const(f)\;\zeta_{\QQ(i)}(2)}\,
\prod_{p|\discr(f)} \big(1+\bigg(\frac{-1}p\bigg)p^{-1}\big)
\;\;s^2 (1+\operatorname{O}(s^{-\kappa}))\;,
\end{multline*}
where $p$ ranges over the odd positive rational primes and
$\big(\frac{-1}p\big)$ is a Legendre symbol. \cqfd
\ecoro

\medskip
Let $\HH, A,\OOO$ and the associated notation, be as in the
beginning of Section \ref{subsect:mertenshamilton}. Let
$f:\HH\times\HH\ra \RR$ be a binary Hamiltonian form, with
\begin{equation}\label{eq:formANT}
f(u,v)=a\,\redn(u)+ \redtr(\ov u\, b \,v) +c\,\redn(v)\;,
\end{equation}
which is integral over $\OOO$ (its coefficients $a=a(f),b=b(f),
c=c(f)$ satisfy $a,c\in\ZZ$ and $b\in\OOO$).  In this subsection, we
assume that $f$ is indefinite (its discriminant $\discr(f)=\n(b)-ac$
is positive).  We denote by $\Ga_\OOO=\SLO$ the Hamilton-Bianchi group
of invertible $2\times 2$ matrices with coefficients in $\OOO$ (see
Subsection \ref{subsect:mertenshamilton} or \cite[\S 3]{ParPau13ANT}
for definitions).  The group $\SLH$ acts on the right on the binary
Hermitian forms by precomposition and the action of $\Ga_\OOO$
preserves the integral ones.  The group $\automH$ of automorphs of $f$
consists of the elements $g\in\Ga_\OOO$ stabilising $f$, that is, such
that $f\circ g=f$.

Let $G$ be a finite index subgroup of $\Ga_\OOO$. For all $x,y$ in
$A$ not both zero, and for every $s>0$, let
$$
\psi_{f,\,G,\,x,\,y}(s)=\card\;\;_{\mbox{$\automH\cap G$}}\bs
\big\{(u,v)\in G(x,y)\;:\;
\n(\OOO x+\OOO y)^{-1}|f(u,v)|\leq s\big\}\;.
$$
Let $\Ga_{\OOO,\,x,\,y}$ and $G_{x,y}$ be the stabilisers of $(x,y)\in
A\times A$ for the left linear actions of $\Ga_\OOO$ and $G$
respectively. Let $K_{x,\,y}$ be the left fractional ideal $\OOO$ if
$xy=0$ and $\OOO x\cap\OOO y$ otherwise, and $\OOO_r(K_{x,\,y})$ its
right order.  Let $\iota_G=1$ if $-\id\in G$, and $\iota_G=2$
otherwise.

\btheo\label{theo:mainversionGANTerror}
There exists $\kappa>0$ such that, as $s$ tends to $+\infty$, 
\begin{multline*}
\psi_{f,\,G,\,x,\,y}(s)= \\
\frac{540\;\iota_G\;[\Ga_{\OOO,\,x,\,y}:G_{x,y}]\;\covol(\automH\cap G)}
{\pi^2\;\zeta(3)\;|\OOO_r(K_{x,\,y})^\times|\;\discr(f)^2
\;[\Ga_\OOO:G]\prod_{p|D_A}(p^3-1)(1-p^{-1})} \;\;\;
s^4(1+\operatorname{O}(s^{-\kappa}))\;,
\end{multline*}
with $p$ ranging over positive rational primes dividing $D_A$. As
$s\to+\infty$, we have
$$
\frac{\zeta(3)\,[\Ga_\OOO:G]\,\discr(f)^4\,\prod_{p|D_A}(p^3-1)(p-1)}
{2160\,\iota_G\,\covol(\automH\cap G)\,s^4}\;
\sum_{f'\in G\cdot f\;:\; 0<|a(f')|\le s}\Delta_{\frac{b(f')}{a(f')}}
\;\weakstar\;\Leb_\HH
$$
and if $y\neq 0$,
$$
\frac{\zeta(3)\,[\Ga_\OOO:G]\,\discr(f)^4\,\prod_{p|D_A}(p^3-1)(p-1)}
{2160\,\iota_G\,\covol(\automH\cap G)\,s^4}
\;\sum_{f'\in\, G\cdot f,\;0<|f'(x,\,y)|\,\leq\, s}
\Delta_{\frac{a(f')\,x\,+\,b(f')\,y}{f'(x,\,y)}}
\;\weakstar\; \Leb_{\HH}\,.
$$
If $\Lambda$ is the lattice of $\lambda\in\OOO$ such that
$\pm\Big(\begin{array}{cc} 1&\lambda \\0&1\end{array} \Big)\in G$,
identifying $\OOO$ with the upper triangular unipotent sugbroup of
$\Ga_\OOO$, we have, as $s\ra+\infty$,
\begin{multline*}
\card\{f'\in\, \Lambda\bs G\cdot f\;:\;0<|a(f')|\leq s\}\\
= \frac{540\,\iota_G\,D_A\,\covol(\automH\cap G)}
{\zeta(3)\,[\Ga_\OOO:G]\,\discr(f)^4\,\prod_{p|D_A}(p^3-1)(p-1)}
\;\;
s^4 (1+\operatorname{O}(s^{-\kappa}))\;.
\end{multline*}
\etheo

For smooth functions $\psi$ with compact support on $\HH$, there is an
error term in the above equidistribution claim, of the form
$\bigO(s^{4-\kappa} \|\psi\|_\ell)$ where $\kappa>0$ and
$\|\psi\|_\ell$ is the Sobolev norm of $\psi$ for some $\ell\in\NN$.

\medskip \dem The proof is completely analogous to that of Theorem
\ref{theo:mainHermerror}, mostly replacing $n=3$ by $n=5$, $\wh K=\CC$
by $\wh K=\HH$, $K$ by $A$, $\OOO_K$ by $\OOO$, $\n$ by $\redn$, and
references to \cite{ParPau11BLMS} to references to \cite{ParPau13ANT},
so that Equation \eqref{eq:quasifinherm} is still valid, and then one
replaces Humbert's formula \eqref{eq:covolBianchi} by Emery's formula
\eqref{eq:covolHamilton}, and the formula $\Vol(\OOO_K\bs \CC)=
\frac{\sqrt{|D_K|}}{4}$ by $\Vol(\OOO\bs \HH)=\frac{D_A}{4}$.  \cqfd

\bigskip Given two left fractional ideals $\mmm, \mmm'$ of $\OOO$ and
$s\geq 0$, let
\begin{multline*}
\psi_{f,\,\mmm,\,\mmm'}(s)\\=\card\;\;_{\mbox{$\automH$}}\bs
\big\{(u,v)\in \mmm\times\mmm\;:\;
\frac{|f(u,v)|}{\n(\mmm)}\leq s\;,\;\;\OOO u+\OOO v=\mmm\;,\;\;
[K_{u,\,v}]=[\mmm'] \ \big\}\;.
\end{multline*}
The following result follows from the first claim of Theorem
\ref{theo:mainversionGANTerror} as in the proof of
\cite[Coro.~18]{ParPau13ANT}.

\bcoro\label{coro:doublefracidealerror}
There exists $\kappa>0$ such that, as $s$ tends to $+\infty$,
$$
\psi_{f,\,\mmm,\,\mmm'}(s)\sim \;\;
\frac{540\;\covol(\automH)}
{\pi^2\;\zeta(3)\;|\OOO_r(\mmm')^\times|\;\discr(f)^2\;
\prod_{p|D_A}(p^3-1)(1-p^{-1})}\;\;\;
s^4(1+\operatorname{O}(s^{-\kappa}))\;. \;\;\;\Box
$$
\ecoro

For every $s\geq 0$, we consider the integer
$$
\psi_{f,\,\mmm}(s)=\card\;\;_{\mbox{$\automH$}}\bs
\big\{(u,v)\in\mmm\times\mmm\;:\;\redn(\mmm)^{-1}|f(u,v)|\leq s,
\;\;\;\OOO u+\OOO v=\mmm\big\}\;,
$$
which is the number of nonequivalent $\mmm$-primitive representations
by $f$ of rational integers with absolute value at most
$s\,\redn(\mmm)$. The next result then follows from Corollary
\ref{coro:doublefracidealerror} as in the proof given in Corollary 19
of Theorem 1 in \cite{ParPau13ANT}.

\bcoro\label{coro:mainintroANTerror} There exists $\kappa>0$ such
that, as $s$ tends to $+\infty$, 
$$
\psi_{f,\,\mmm}(s)= \frac{45\;D_A\;\covol(\automH)}
{2\,\pi^2\;\zeta(3)\;\discr(f)^2\;\prod_{p|D_A}(p^3-1)}\;\;
s^4(1+\operatorname{O}(s^{-\kappa}))\;.\;\;\;\Box
$$
\ecoro

\medskip
The group of automorphs of the binary Hamiltonian form $(u,v)\mapsto
\redtr(\ov u\,v)$ (which is the standard real scalar product on $\HH$)
is
$$
\SpO=\Big\{g\in\SLO\;:\;
^t\overline{g}\,\Big(\begin{array}{cc}0& 1\\1&
  0\end{array}\Big)\,g=\Big(\begin{array}{cc}0& 1\\1&
  0\end{array}\Big)\Big\}\;,
$$
which is an arithmetic lattice in the symplectic group
$\operatorname{Sp}_1(\HH)$ over Hamilton's quaternion algebra.  The
following result follows from Corollary \ref{coro:mainintroANTerror}
as in the proof of \cite[Coro.~2]{ParPau13ANT}.

\bcoro\label{coro:introspunoooerror}
There exists $\kappa>0$ such that, as $s$ tends to $+\infty$,
\begin{multline*}
\card\;\;_{\mbox{$\SpO$}}\bs
\big\{(u,v)\in \OOO\times \OOO\;:\;
|\redtr(\ov u\,v)|\leq s,\;\OOO u+ \OOO v=\OOO\;\big\}
= \\\frac{D_A}
{48\,\zeta(3)}\prod_{p|D_A}\frac{p^2+1}{p^2+p+1}\;\;
s^4(1+\operatorname{O}(s^{-\kappa}))\;.\;\;\;\Box
\end{multline*}
\ecoro

\subsection{Positive definite binary Hermitian and 
Hamiltonian forms}
\label{sect:posdef}

Let $K$ and $\mmm$ be as in the beginning of Subsection
\ref{subsect:mertensbianchi}. Let $f$ be an integral binary Hermitian
form over $\OOO_K$ as in the beginning of Subsection
\ref{subsect:errortermintrepbybinary}.  In this subsection, we assume
that $f$ is positive definite (that is $\discr(f)<0$ and $a(f)>0$), a
feature which is preserved by the action of the Bianchi group
$\Ga_K=\SLOK$. Note that the group $\autom$ of automorphs of $f$ is
then finite.

For every $s>0$, we consider the following integer (its finiteness is
part of the following proof), depending only on the ideal class of
$\mmm$,
$$
\psi^+_{f,\,\mmm}(s)=\card\;
\big\{(u,v)\in\mmm\times\mmm\;:\;\n(\mmm)^{-1}f(u,v)\leq s,
\;\OOO_K u+\OOO_K v=\mmm\big\}\;,
$$
which is the number of nonequivalent $\mmm$-primitive representations
by $f$ of the rational integers at most $s\,\n(\mmm)$.

\btheo \label{theo:posdefHermerror}
There exists $\kappa>0$ such that, as $s$ tends to $+\infty$, we have
$$
\psi^+_{f,\,\mmm}(s)\;\sim\; \frac{\pi^2}
{|D_K|\;\zeta_K(2)\;|\discr(f)|}\;\;\;
s^2 (1+\operatorname{O}(s^{-\kappa}))\;.
$$
If $G$ is a finite index subgroup of $\Ga_K$, with $\iota_G=1$ if
$-\id\in G$ and $\iota_G=2$ otherwise, then, as $s\to+\infty$,
$$
\frac{|\autom\cap G|\,[\Ga_K:G]\,|D_K|^{\frac 32}\;\zeta_K(2)\;
|\discr(f)|}{2\,\pi^2\,\iota_G\;s^2}\;
\sum_{f'\in G\cdot f,\; a(f')\le s}\Delta_{\frac{b(f')}{a(f')}}
\;\weakstar\;\Leb_\CC\,.
$$
If $\Lambda$ is the lattice of $\lambda\in\OOO_K$ such that
$\pm\Big(\begin{array}{cc} 1&\lambda \\0&1\end{array} \Big)\in G$,
identifying $\OOO_K$ with the upper triangular unipotent sugbroup of
$\Ga_K$, we have, as $s\ra+\infty$,
\begin{multline*}
\card\{f'\in\, \Lambda\bs G\cdot f\;:\;a(f')\leq s\}\\
= \frac{\pi^2\,\iota_G\,[\OOO_K:\Lambda]}
{|\autom\cap G|\,[\Ga_K:G]\,|D_K|\,\zeta_K(2)\,|\discr(f)|}
\;\;
s^2 (1+\operatorname{O}(s^{-\kappa}))\;.
\end{multline*}
\etheo

For smooth functions $\psi$ with compact support on $\CC$, there is an
error term in the above equidistribution claim, of the form
$\bigO(s^{2-\kappa} \|\psi\|_\ell)$ where $\kappa>0$ and
$\|\psi\|_\ell$ is the Sobolev norm of $\psi$ for some $\ell\in\NN$.

\medskip \dem Let $n=3$ and $\wh K=\CC$. 

We start the proof by describing the space of positive definite binary
Hermitian forms in hyperbolic geometry terms. Let $\Q^+$ be the cone
of positive definite binary Hermitian forms. The multiplicative group
$\RR_+= \; ]0,+\infty[$ acts on $\Q^+$ by multiplication and we define
$\ov\Q^+ =\Q^+/\RR_+$. The map $\Phi:\ov\Q^+\to\hnr=\wh K\times\RR_+$
induced by $f\mapsto \big(-\frac{b(f)}{a(f)}, \frac{\sqrt{-\discr(f)}}
{a(f)}\big)$ is a homeomorphism, which is (anti-)equivariant in the
sense that $\Phi(f\circ g)=g^{-1}\Phi(f)$ for all $g\in\SL_2(\wh K)$
and $f\in\Q^+$ (see for example \cite[Prop. 22]{ParPau13ANT} for a
proof for positive definite binary Hamiltonian forms, the above claims
being obtained by embedding as usual $\CC$ in Hamilton's quaternion
algebra $\HH$).

Now, let $x,y$ be not both zero in $\OOO_K$, with $x=1$ if $y=0$, and
let $I_{x,\,y}=\OOO_Kx+\OOO_Ky$. We will more generally study the following
counting function of the representations, in a given orbit of
$G$, of integers by $f$, defined, for all $s\geq 0$, by
$$
\psi^+_{f,\,G,\,x,\,y}(s)= \card\;
\big\{(u,v)\in G(x,y)\;:\; \n(I_{x,\,y})^{-1}\,f(u,v)\leq s\big\}\;.
$$

Let $\rho$, $\ga_\rho$ be as in the proof of Theorem
\ref{theo:mainversionGerror}.  We will apply Corollary
\ref{coro:constcurvcount} with $\Ga$ the arithmetic group
$\ga_\rho^{-1}\,\overline{G}\,\ga_\rho$, $D^-$ the horoball
$\H_\infty=\{(z,t)\in\hnr\;:\;t\geq 1\}$ (which is centered at a
parabolic fixed point of $\Ga$, and whose pointwise stabiliser is
trivial), and $D^+$ the singleton consisting of
$\Phi(f\circ\ga_\rho)$.

For any $g\in\Ga_K$, the distance of $\Phi(f\circ g)$ from $\H_\infty$
is $\ln\frac{a(f\circ g)} {\sqrt{-\discr(f))}}$, when this number is
positive, which is the case except for finitely many right classes of
$g$ under the stabiliser of $\H_\infty$.  Note that $\ga\in \Ga$
belongs to the stabiliser $\Ga_{D^+}$ of $D^+$ in $\Ga$ if and only if
$\ga^{-1}$ stabilises the positive homothety class of $f\circ
\ga_\rho$, that is, since $\SL_2(\wh K)$ preserves the discriminant,
if and only if $\ga^{-1}$ fixes $f\circ\ga_\rho$. Thus $\Ga_{D^+}=
\ga_\rho^{-1}( \autom\cap G)\ga_\rho$, which is also the pointwise
stabiliser of $D^+$, whose order is hence $|\autom\cap G|$.

Let $G_{\H_\rho}$ and $G_{x,\,y}$ be the stabilisers in $G$ of the
horoball $\H_\rho= \ga_\rho\H_\infty$ in $\hnr$ and of the pair
$(x,y)$ in $\wh K\times \wh K$.  As previously, we denote by $\ov{H}$
the image in $\PSL_2(\wh K)$ of any subgroup $H$ of $\SL_2(\wh K)$. We
use the (surprising !) convention that $\n(y)=1$ if $y=0$ to avoid
considering cases. The counting result with error term of Corollary
\ref{coro:constcurvcount}, and the value of $C(D^-,D^+)$ in case (4)
just above it, using cosets and multiplicities arguments already seen,
shows that
\begin{align}
  \psi^+_{f,\,G,\,x,\,y}(s)&= \card\; \big\{[\ga]\in G/G_{x,\,y}\;:\;
  \n(I_{x,\,y})^{-1}\,f(\ga(x,y))\leq s\big\}
  \nonumber\\
  &= \frac{2}{\iota_G}\,\card\; \big\{[\ga]\in
  \ov{G}/\,\overline{G_{x,\,y}}\;:\;
  \,f\circ\ga\circ\ga_\rho(1,0)\leq \frac{s\,\n(I_{x,\,y})}{\n(y)}\big\}
  \nonumber\\
  &= \frac{2}{\iota_G}\,[\,\overline{G_{\H_\rho}}:\overline{G_{x,\,y}}\,]\,
  \card\; \big\{[\ga]\in \ov{G}/\,\overline{G_{\H_\rho}}\;:\;
  \,a(f\circ\ga\circ\ga_\rho)\leq \frac{s\,\n(I_{x,\,y})}{\n(y)}\big\}
  +\bigO(1) \nonumber \\&=
  \frac{2}{\iota_G}\,[\,\overline{G_{\H_\rho}}:\overline{G_{x,\,y}}\,]\,
  |\autom\cap G|\, \nonumber\\&\quad\quad\card\; \big\{[\ga]\in
  \Ga_{D^+}\bs\Ga/\,\Ga_{D^-} \;:\;
  \,a(f\circ\ga_\rho\circ\ga)\leq
  \frac{s\,\n(I_{x,\,y})}{\n(y)}\big\}+\bigO(1)
  \nonumber\\&
 = \frac{2}{\iota_G}\,[\,\overline{G_{\H_\rho}}:\overline{G_{x,\,y}}\,]\,
  |\autom\cap G|\, \N_{D^+,\,D^-}\big(\ln
  \frac{s\,\n(I_{x,\,y})}{\n(y)\,\sqrt{-\discr(f)}}\big)+\bigO(1)
  \nonumber\\&
 = \frac{2\,[\,\overline{G_{\H_\rho}}:\overline{G_{x,\,y}}\,]\,
\Vol(\Ga_{D^-}\bs D^-)}{\iota_G\,\Vol(\Ga\bs \hnr)}\Big(
  \frac{s\,\n(I_{x,\,y})}{\n(y)\,\sqrt{-\discr(f)}}\Big)^{n-1}
(1+\bigO(s^{-\kappa}))\,.\label{eq:loncalcdefpos}
\end{align}
As seen in the proof of Theorem \ref{theo:appliHerm}, we have
\begin{align}
\Vol(\Ga_{D^-}\bs D^-)&=\Vol(\,\overline{G_{\H_\rho}}\,\bs \H_\rho)=
[\,\overline{(\Ga_K)_{\H_\rho}}:\overline{G_{\H_\rho}}\,]
\,\Vol(\,\overline{(\Ga_K)_{\H_\rho}}\,\bs \H_\rho)\nonumber\\ &=
[\,\overline{(\Ga_K)_{\H_\rho}}:\overline{G_{\H_\rho}}\,]\;
\frac{\sqrt{|D_K|}}{2\,{\omega_K}}\frac{\n(y)^2}{\n(I_{x,\,y})^2}\,.
\label{eq:calccovolhororho}
\end{align}
We have
\begin{align}
[\,\overline{(\Ga_K)_{\H_\rho}}:\overline{G_{\H_\rho}}\,]
[\,\overline{G_{\H_\rho}}:\overline{G_{x,\,y}}\,]&=
[\,\overline{(\Ga_K)_{\H_\rho}}:\overline{\Ga_{K,\,x,\,y}}\,]
[\,\overline{\Ga_{K,\,x,\,y}}:\overline{G_{x,\,y}}\,]\nonumber\\ &
=\frac{\omega_K}{2}\,[\Ga_{K,\,x,\,y}:G_{x,\,y}]\,.
\label{eq:calcindiceherm}
\end{align}
Using Equation \eqref{eq:modifcovolBianchi}, Equation
\eqref{eq:loncalcdefpos} then gives the following formula, interesting
in itself,
$$
\psi^+_{f,\,G,\,x,\,y}(s)=
\frac{2\,\pi^2\,[\Ga_{K,\,x,\,y}:G_{x,\,y}]}
{|D_K|\,\zeta_K(2)\,|\discr(f)|\,[\Ga_{K}:G]}\;s^2(1+\bigO(s^{-\kappa}))\;.
$$

Since any nonzero ideal $\mmm$ in $\OOO_K$ may be written $\mmm=
\OOO_Kx_\mmm+\OOO_ky_\mmm$ for some $x_\mmm,y_\mmm$ not both zero in
$\OOO_K$ and $\psi^+_{f,\,\mmm}=
\psi^+_{f,\,\Ga_K,\,x_\mmm,\,y_\mmm}$, the first claim of Theorem
\ref{theo:posdefHermerror} follows.

\medskip Let us now prove the equidistribution result in Theorem
\ref{theo:posdefHermerror}. The main additional remark is that for all
$f'\in G\cdot f$, the initial point of the common perpendicular from
$\H_\infty$ to $\Phi(f'\circ\ga_\rho)$ is $\big(-\frac{b(f'\circ\ga_\rho)}
{a(f'\circ\ga_\rho)},\,1\big)$. Taking $t=\ln\frac{s}
{\sqrt{-\discr(f)}}$ in Equation \eqref{eq:distribhorob} in
Corollary \ref{coro:constcurvcount}, gives, as $s\ra+\infty$,
$$
\frac{|\autom\cap G|\,(n-1)\,\Vol(\Ga\bs \hnr)}
{\big(\frac{s}{\sqrt{-\discr(f)}}\big)^{n-1}}
\;\sum_{\ga\in\Ga/\Ga_{D^+},\;a(f\circ\ga_\rho\circ\ga)\leq s}
\Delta_{\big(-\frac{b(f\circ\ga_\rho\circ\ga)}
{a(f\circ\ga_\rho\circ\ga)},\,1\big)}\;\weakstar \;\Vol_{\partial D^-}\,.
$$ Hence, since the map from $G\cdot f$ to $\Ga/\Ga_{D^+}$ defined by
$f\circ \ga\mapsto \ga_\rho^{-1}\ga\ga_\rho\Ga_{D^+}$ is a bijection,
using the pushforward of measures by the map $(z,1)\mapsto -\,z$ from
$\partial D^-$ to $\wh K$, we have
\begin{multline}\label{eq:equidposdef}
|\autom\cap G|\,(n-1)\,\Vol(\Ga\bs \hnr)
\,|\discr(f)|^{\frac{n-1}{2}}\;s^{-(n-1)}
\;\sum_{f'\in G\cdot f,\;a(f'\circ\ga_\rho)\leq s}
\Delta_{\frac{b(f'\circ\ga_\rho)}{a(f'\circ\ga_\rho)}}\\
\;\weakstar \;\Leb_{\wh K}\,.
\end{multline}

Taking $y=0$ (then $\ga_\rho=\id$) and using Equation
\eqref{eq:modifcovolBianchi}, this gives the equidistribution claim
in Theorem \ref{theo:posdefHermerror}.  The last claim of Theorem
\ref{theo:mainversionGerror} follows from it as in the proof of
Theorem \ref{theo:mainversionGerror}.

The case $y\neq 0$, using Equation \eqref{eq:modifcovolBianchi}, gives
the following equidistribution result, interesting in itself, as in the
end of the proof of Theorem \ref{theo:mainversionGerror}:
\begin{multline*}
\frac{|\autom\cap G|\,[\Ga_K:G]\,|D_K|^{\frac 32}\;\zeta_K(2)\;
|\discr(f)|\,\n(y)}{2\,\pi^2\,\iota_G\;s^2}
\;\sum_{f'\in\, G\cdot f,\;f'(x,\,y)\,\leq\, s}
\Delta_{\frac{a(f')\,x\,+\,b(f')\,y}{f'(x,\,y)}}\\
\;\weakstar\; \Leb_{\CC}\,.\;\;\;\Box
\end{multline*}

\medskip Let $\HH, A,\OOO,\mmm$ and the associated notation be as in
the beginning of Section \ref{subsect:mertenshamilton}. Let $f$ be an
integral binary Hamiltonian form over $\OOO$, as in Subsection
\ref{subsect:errortermintrepbybinary}. In this subsection, we assume
that $f$ is positive definite (that is $\discr(f)<0$ and $a(f)>0$), a
feature which is preserved by the action (on the right) of the
Hamilton-Bianchi group $\Ga_\OOO=\SLO$ (by precomposition). Note that
the group $\automH$ of automorphs of $f$ is then finite.

For every $s>0$, we consider the integer
$$
\psi^+_{f,\,\mmm}(s)=
\card\;\big\{(u,v)\in\mmm\times\mmm\;:\;\redn(\mmm)^{-1}f(u,v)\leq
s, \;\;\;\OOO u+\OOO v=\mmm\big\}\;,
$$
which is the number of nonequivalent $\mmm$-primitive representations
by $f$ of the rational integers  at most $s\redn(\mmm)$.

\btheo \label{theo:posdefHamerror}
There exists $\kappa>0$ such that, as $s$ tends to $+\infty$,
$$
\psi^+_{f,\,\mmm}(s)\;\sim\; \frac{60\,D_A}
{\zeta(3)\,|\discr(f)|^2\,\prod_{p|D_A}(p^3-1)}\;\;\;
s^4 (1+\operatorname{O}(s^{-\kappa}))\;.
$$
If $G$ is a finite index subgroup of $\Ga_\OOO$, with $\iota_G=1$ if
$-\id\in G$ and $\iota_G=2$ otherwise, then, as $s\to+\infty$,
\begin{multline*}
\frac{|\automH\cap G|\,[\Ga_\OOO:G]\,\zeta(3)\;
|\discr(f)|^2\,\prod_{p|D_A}(p^3-1)(p-1)}{2880\,\iota_G\;s^4}\;
\sum_{f'\in G\cdot f,\; a(f')\le s}\Delta_{\frac{b(f')}{a(f')}}\\
\;\weakstar\;\Leb_\HH\,.
\end{multline*}
If $\Lambda$ is the lattice of $\lambda\in\OOO$ such that
$\pm\Big(\begin{array}{cc} 1&\lambda \\0&1\end{array} \Big)\in G$,
identifying $\OOO$ with the upper triangular unipotent sugbroup of
$\Ga_\OOO$, we have, as $s\ra+\infty$,
\begin{multline*}
\card\{f'\in\, \Lambda\bs G\cdot f\;:\;a(f')\leq s\}\\
= \frac{720\,\iota_G\,[\OOO:\Lambda]\,D_A}
{|\automH\cap G|\,[\Ga_\OOO:G]\,\zeta(3)\;
|\discr(f)|^2\,\prod_{p|D_A}(p^3-1)(p-1)}
\;\;
s^4 (1+\operatorname{O}(s^{-\kappa}))\;.
\end{multline*}
\etheo

\dem The proof is completely analogous to that of Theorem
\ref{theo:posdefHermerror}, mostly replacing $n=3$ by $n=5$, $\wh
K=\CC$ by $\wh K=\HH$, $K$ by $A$, $\OOO_K$ by $\OOO$, $\n$ by
$\redn$.  Given $x,y$ not both zero in $\OOO$, with $x=1$ if $y=0$,
let $I_{x,\,y}$ and $K_{x,\,y}$ be as in the beginning of the proof of
Theorem \ref{theo:appliHam}. We also introduce the following counting
function of the representations of integers by $f$, in a given orbit
of $G$:
$$
\psi^+_{f,\,G,\,x,\,y}(s)= \card\;
\big\{(u,v)\in G(x,y)\;:\; \redn(I_{x,\,y})^{-1}\,f(u,v)\leq s\big\}\;.
$$ 
Keeping the relevant notation of the proof Theorem
\ref{theo:posdefHermerror}, Equation \eqref{eq:loncalcdefpos} is still
valid. By Equation \eqref{eq:xycuspvol} where $\tau=1$, we have
\begin{align*}
\Vol(\Ga_{D^-}\bs D^-)&=\Vol(\,\overline{G_{\H_\rho}}\,\bs \H_\rho)=
[\,\overline{(\Ga_\OOO)_{\H_\rho}}:\overline{G_{\H_\rho}}\,]
\,\Vol(\,\overline{(\Ga_\OOO)_{\H_\rho}}\,\bs \H_\rho)\\ &=
[\,\overline{(\Ga_\OOO)_{\H_\rho}}:\overline{G_{\H_\rho}}\,]\;
\frac{D_A\,\redn(y)^4}{16\,|\OOO_r(K_{x,\,y})^\times|\,  
[\,\overline{(\Ga_\OOO)_{\H_\rho}}:\overline{(\Ga_\OOO)_{x,\,y}}\,]\,
\redn(I_{x,\,y})^4}\,.
\end{align*}
By Emery's formula \eqref{eq:covolHamilton}, we have
\begin{align*}
\Vol(\Ga\bs\hnr)&=\Vol(\,\overline{G}\,\bs\hcr)=
[\,\overline{\Ga_\OOO}:\overline{G}\,]
\Vol(\,\overline{\Ga_\OOO}\,\bs\hcr) \\ &
=\frac{[\Ga_\OOO:G]}{\iota_G}\,
\frac{\zeta(3)\prod_{p|D_A}(p^3-1)(p-1)}{11520}\,.
\end{align*}
Since $\frac{[\,\overline{G_{\H_\rho}}:\overline{G_{x,\,y}}\,]\,
  [\,\overline{(\Ga_\OOO)_{\H_\rho}}:\overline{G_{\H_\rho}}\,]}
{[\,\overline{(\Ga_\OOO)_{\H_\rho}}:\overline{(\Ga_\OOO)_{x,\,y}}\,]}
=[\,\overline{(\Ga_\OOO)_{x,\,y}}:\overline{G_{x,\,y}}\,]=
[(\Ga_\OOO)_{x,\,y}:G_{x,\,y}]$, we hence have the following counting
asymptotic, interesting in itself:
\begin{equation}\label{eq:countdefposgeneHam}
\psi^+_{f,\,G,\,x,\,y}(s)=
\frac{1440\,D_A\,\iota_G\,[(\Ga_\OOO)_{x,\,y}:G_{x,\,y}]}
{\zeta(3)\,[\Ga_\OOO:G]\,|\OOO_r(K_{x,\,y})^\times|\,|\discr(f)|^2\,
\prod_{p|D_A}(p^3-1)(p-1)}\;\;s^4(1+\bigO(s^{-\kappa}))\,.
\end{equation}
Given two left fractional ideals $\mmm, \mmm'$ of $\OOO$ and
$s\geq 0$, let 
$$
\psi^+_{f,\,\mmm,\,\mmm'}(s)=
\card\;\big\{(u,v)\in\mmm\times\mmm\;:\;\n(\mmm)^{-1}f(u,v)\leq
s, \;I_{u,\,v}=\mmm, \, [K_{u,\,v}]=[\mmm']\big\}\;.
$$
As explained in the proof of Theorem \ref{theo:appliHam}, taking
$(x,y)\in\OOO\times\OOO$ such that $[I_{x,\,y}]=[\mmm]$ and
$[K_{x,\,y}]=[\mmm']$, we have
$\psi^+_{f,\,\mmm,\,\mmm'}=\psi^+_{f,\,\Ga_\OOO,\,x,\,y}$. Hence
\begin{align*}
\psi^+_{f,\,\mmm,\,\mmm'}(s)=
\frac{1440\,D_A}
{\zeta(3)\,|\OOO_r(\mmm')^\times|\,|\discr(f)|^2\,
\prod_{p|D_A}(p^3-1)(p-1)}\;\;s^4(1+\bigO(s^{-\kappa}))\,.
\end{align*}
Now, the first claim of Theorem \ref{theo:posdefHamerror} follows from
Equation \eqref{eq:deuring}, since
$$
\psi^+_{f,\,\mmm}=\sum_{[\mmm']\;\in\;_\OOO\!\I}\;\psi^+_{f,\,\mmm,\,\mmm'}\;,
$$

The last two assertions of Theorem \ref{theo:posdefHamerror} are
proven in the same way as the last two assertions of Theorem
\ref{theo:posdefHermerror}, since Equation \eqref{eq:equidposdef} is
still valid, and if $y\neq 0$, we furthermore have
\begin{multline*}
\frac{|\automH\cap G|\,[\Ga_\OOO:G]\,\zeta(3)\,
|\discr(f)|^2\,\n(y)^{2}\,\prod_{p|D_A}(p^3-1)(p-1)}
{2880\,\iota_G\;s^4}\;\times\\
\;\sum_{f'\in\, G\cdot f,\;f'(x,\,y)\,\leq\, s}
\Delta_{\frac{a(f')\,x\,+\,b(f')\,y}{f'(x,\,y)}}
\;\weakstar\; \Leb_{\HH}\,.\;\;\;\Box
\end{multline*}

\subsection{Representation of algebraic integers by 
integral binary quadratic forms}
\label{subsec:normform}

In this final subsection, we study counting and equidistribution
problems of the representations of algebraic integers by quadratic norm
forms (for related results, see for instance \cite{Nagell62},
\cite[Thm. R]{Odoni78}, \cite[Chap.~VI]{Lang94}, as well as \cite[\S
3.1]{GorPau12} for an ergodic approach).

Let $K$ and its associated notation be as in the beginning of Section
\ref{sect:irrationals}. Let $\alpha_0\in\wh K$ be a fixed quadratic
irrational over $\OOO_K$. We also denote by $\n$ the relative norm in
$K(\alpha_0)$ over $K$ (which is consistent with the notation of the
beginning of Section \ref{sect:irrationals}). We consider the
(relative) norm form, seen as a map from $K\times K$ to $K$, defined
by
$$
N_{\alpha_0}:(u,v)\mapsto \n(u-\alpha_0 v)=
u^2-\tr\,\alpha_0\;uv+\n(\alpha_0)\,v^2\,,
$$
which is the homogeneous form of the minimal polynomial of $\alpha_0$
over $K$. Its values on $\OOO_K\times\OOO_K$ belong to $\OOO_K$ if
$\alpha_0$ is an algebraic integer, and to $a\OOO_K$ for some $a\in K$ in
general. We will study the representations of elements of $K$ by this
norm form $N_{\alpha_0}$, in orbits of (finite index subgroups of) the
modular groups $\Ga_K=\SL_2(\OOO_K)$ (a minor change of notation from
the beginning of Section \ref{sect:irrationals}).

We denote by $\ov{H}$ the image in $\PSL_2(\wh K)$ of any subgroup $H$
of $\SL_2(\wh K)$. We denote by $H_P$ the stabiliser, for any group
action of a group $H$ on a set $X$, of a point $P$ or a subset $P$ in
$X$. The stabiliser $\stab_{\Ga_K}N_{\alpha_0}$ of the norm form
$N_{\alpha_0}$, for the action (on the right) of $\Ga_K$ by
precomposition by the linear action, is exactly
$(\Ga_K)_{\{\alpha_0, \,\alpha_0^\sigma\}}$.

As a warm-up, when $K=\QQ$, we give an equidistribution result of the
fractions of the representations, in an aformentioned orbit, of usual
rationals by this norm from (which, in this case, is an indefinite
quadratic form, hence the classical counting results of the beginning
of Section \ref{sect:errorterms} apply).

\btheo\label{theo:posdefherm} For every finite index subgroup $G$ of
$\Ga_\QQ=\SLZ$, with $\iota_G=1$ if $-\id\in G$ and $\iota_G=2$
otherwise, we have the following convergence of
measures on $\RR-\{\alpha_0,\alpha_0^\sigma\}$ as $s\to+\infty$:
$$
\frac{\pi^2\,\iota_G\,[\Ga_\QQ:G]}{12\,[(\Ga_\QQ)_{\infty}:G_\infty]\;s}\;
\sum_{(u,v)\in G(0,1),\;|\n(u-\alpha_0 v)\,|\leq s} \Delta_{\frac uv}
\;\weakstar\;\frac{d\Leb_\RR(t)}{|Q_{\alpha_0}(t)|}\,.
$$
\etheo

For smooth functions $\psi$ with compact support on
$\RR-\{\alpha_0,\alpha_0,\}$, there is an error term in the above
equidistribution claims, of the form $\bigO(s^{1-\kappa}
\|\psi\|_\ell)$ where $\kappa>0$ and $\|\psi\|_\ell$ is the Sobolev
norm of $\psi$ for some $\ell\in\NN$.

For instance, taking $G=\SLZ$ gives Theorem \ref{theo:normformintro}
in the introduction, and taking for $G$ the Hecke congruence subgroup
modulo $k\in\NN-\{0\}$, we have, by the index computation of
for instance \cite[p.~24]{Shimura71},
$$
\frac{\pi^2\,k\prod_{\ppp|k}\big(1+\frac{1}{\ppp}\big)}{12\,s}\;
\sum_{(u,v)\in \ZZ^2,\;(u,v)=1,\; v\equiv 0\;[k],\;|\n(u-\alpha_0 v)\,|\leq s} 
\Delta_{\frac uv} \;\weakstar\;\frac{d\Leb_\RR(t)}{|Q_{\alpha_0}(t)|}\,.
$$

\medskip \dem We take $n=2$ and $\wh K=\RR$ in this proof. We will
apply Equation \eqref{eq:distribtotgeod} with $\Ga$ the arithmetic
group $\overline{G}$, $D^-$ the geodesic line with points at infinity
$\alpha_0,\alpha_0^\sigma$ (on which the group $\wh{\alpha_0}^\ZZ\cap
\overline{G}$ acts with compact quotient) and $D^+$ the horoball
$\H_\infty=\{(z,t)\in\hnr=\wh K\times\RR_+\;:\;t\geq 1\}$ (which is
centered at a parabolic fixed point of $\overline{G}$).

For every $\ga=\begin{pmatrix} a & b \\ c & d\end{pmatrix}$ in $\SLOK$,
the distance $\ell_{\ga}$ between $D^-$ and $\ga D^+$, when they
are disjoint (which is the case except for finitely many double
classes of $\ga$ in $\Ga_{D^-}\bs\Ga/\Ga_{D^+}$), is equal to the
distance between $D^+=\H_\infty$ and the geodesic line $\ga^{-1}
D^-$ with points at infinity $\ga^{-1}\alpha_0, \ga^{-1}
\alpha_0^\sigma$.  Hence by \cite[Lem. 4.2]{ParPau12JMD}
$$
\ell_{\ga}=
\Big|\ln\frac{|\ga^{-1}\alpha_0-\ga^{-1}\alpha_0^\sigma|}{2}\,\Big|
=\ln h(\ga^{-1}\alpha_0)=
\ln \big(h(\alpha_0)\, |\n(u-v\alpha_0)\,|\big)\;.
$$
Let $v_{\ga}$ be the initial tangent vector of the common
perpendicular from $D^-$ to $\ga D^+$, when they are disjoint, and
note that its positive point at infinity $(v_{\ga})_+$ is the
point at infinity of the horoball $\ga D^+$, which is
$\ga\infty=\frac{a}{c}$.  

By coset and multiplicities argument already seen, by the version for
initial tangent vectors of Equation \eqref{eq:distribtotgeod} seen in
Subsection \ref{subsec:relcompquadirrat}, we have
\begin{equation}\label{eq:geodhorodimde}
\frac{\Vol(\SSS^{n-1})\,\Vol(\Ga\bs\hnr)}
{\Vol(\SSS^{n-2})\,\Vol(\Ga_{D^+}\bs D^+)\;e^{t}}
\sum_{\ga\in\Ga/\Ga_{D^+}}\Delta_{v_{\ga}}
\;\weakstar\;\Vol_{\normalout D^-}\,.
\end{equation}
We have 
$$
\Vol(\Ga\bs\hnr)=[\,\overline{\Ga_\QQ}:\overline{G}\,]\,
\Vol(\,\overline{\Ga_\QQ}\,\bs\hdr)=
\frac{[\Ga_\QQ:G]}{\iota_G}\frac{\pi}{3}\,,
$$
and 
$$
\Vol(\Ga_{D^+}\bs D^+)=[\,\overline{(\Ga_\QQ)_\infty}:\overline{G_\infty}\,]\,
\Vol(\,\overline{(\Ga_\QQ)_\infty}\,\bs D^+)=
\frac{[(\Ga_\QQ)_\infty:G_\infty]}{\iota_G}\,.
$$
Taking $t=\ln(h(\alpha_0)s)$ in Equation \eqref{eq:geodhorodimde},
since the image of $\Vol_{\normalout D^-}$ by the pushforward of
measures by the positive endpoint map is $\frac{2^{n-1}\,d\Leb_{\wh K}
  (z)} {h(\alpha_0)^{n-1}\; |Q_{\alpha_0}(z)|^{n-1}}$ by (the comment
following) Lemma \ref{lem:calcmeaspushskinn}, and since the canonical
map $G/G_{1,0}\ra \ov{G}/\overline{G_\infty}=\Ga/\Ga_{D_{D^+}}$ is
$\frac{2}{\iota_G}$-to-$1$, the result follows.  
\cqfd

\bigskip We now assume that $K$ is an imaginary quadratic number
field. Let $\mmm$ be a (nonzero) fractional ideal of $\OOO_K$. Let
$\alpha_0\in\wh K$ be a quadratic irrational over $K$.  We consider
the counting function $\ov\psi_{\mmm}$ of nonequivalent
$\mmm$-primitive representations, by the norm form $N_{\alpha_0}$, of
elements of $K$ (necessary integers in $K$ if $\alpha_0$ is an
algebraic integer), defined on $[0,+\infty[$ by
$$
s\mapsto\card\;\;_{\stab_{\Ga_K}\!N_{\alpha_0}}\bs
\big\{(u,v)\in\mmm\times\mmm\;:\;\n(\mmm)^{-1}|\n(u-v\alpha_0)|\le s,
\ \OOO_K u+\OOO_K v=\mmm\big\}\,.
$$
Since for any $b\in\OOO_K$, we have $\n(b\mmm)=|b|^2\n(\mmm)$, the
counting function $\ov\psi_{\mmm}$ depends only on the ideal class of
$\mmm$.

\btheo \label{theo:laquillebordel}
There exists $\kappa>0$ such that, as $s$ tends to $+\infty$,
$$
\ov\psi_\mmm(s)= \frac{2\,\pi^2\,h(\alpha_0)^2\,
\big|\log\big|\frac{\operatorname{tr}\,\wh{\alpha_0}+
    \sqrt{(\operatorname{tr}\,\wh{\alpha_0})^2-4}}{2}\big|\;\big|}
{m_{\Ga_K}(\alpha_0)\,\iota_{\Ga_K}(\alpha_0)\,|D_K|\,\zeta_K (2)}
\;s^2+\operatorname{O}(s^{2-\kappa})\,.
$$
Furthermore, for every finite index subgroup $G$ of $\Ga_K=\SLOK$,
with $\iota_G=1$ if $-\id\in G$ and $\iota_G=2$ otherwise, we have the
following convergence of measures on
$\CC-\{\alpha_0,\alpha_0^\sigma\}$ as $s\to+\infty$:
$$
\frac{[\Ga_K:G]\,|D_K|\,\zeta_K(2)\,\iota_G\,\omega_K}
{8\,\pi^2\,[(\Ga_K)_\infty:G_\infty]\,\,s^2}
\sum_{(u,v)\in G(0,1),\;|\n(u-\alpha_0 v)\,|\leq s} \Delta_{\frac uv}
\;\weakstar\;\frac{d\Leb_\CC(z)}{|Q_{\alpha_0}(z)|^2}\;.
$$
\etheo

For smooth functions $\psi$ with compact support on $\CC$, there is an
error term in the above equidistribution claim, of the form
$\bigO(s^{2-\kappa} \|\psi\|_\ell)$ where $\kappa>0$ and
$\|\psi\|_\ell$ is the Sobolev norm of $\psi$ for some $\ell\in\NN$.

We leave to the reader a version of this equidistribution claim where
$(0,1)$ is replaced by any fixed pair $(x,y)$ of  elements
of $\OOO_K$ which are not both zero (replacing $D^+=\H_\infty$ by $D^+=\ga_\rho\H_\infty$ in
its proof, with $\ga_\rho$ as in the proof of Theorem
\ref{theo:mainversionGerror}).

For instance, taking $G=\Ga_K$, we have
$$
\frac{|D_K|\,\zeta_K(2)\,\omega_K}
{8\,\pi^2\,\,s^2}
\sum_{(u,v)\in \OOO_K\times\OOO_K,\;\OOO_K u+\OOO_Kv=\OOO_K,\;
|\n(u-\alpha_0 v)\,|\leq s} \Delta_{\frac uv}
\;\weakstar\;\frac{d\Leb_\CC(z)}{|Q_{\alpha_0}(z)|^2}\;.
$$

\medskip \dem We will prove a stronger counting claim.  Let $x,y$ in
$\OOO_K$ be elements which are not both zero, with $x=1$ if $y=0$ and the same strange
convention that $\n(y)=1$ if $y=0$. Let $I_{x,\,y}=x\OOO_K+y\OOO_K$,
and for every $s>0$, let
$$
\ov\psi_{G,\,x,\,y}(s)= \card\;\;\stab_{G}N_{\alpha_0}\bs
\big\{(u,v)\in G(x,y)\;:\; \n(I_{x,\,y})^{-1}|\n(u-v\alpha_0)|\leq
s\big\}\;.
$$
Let $\rho=x\,y^{-1}\in K\cup\{\infty\}$, $\ga_\rho=
\begin{pmatrix}\rho& -1\\1& \;\,0\end{pmatrix}\in\SL_2(K)$ if $y\neq
0$ and $\ga_\rho=\id$ otherwise. Let $\H_\infty= \{(z,t)\in\htr
\;:\;t\geq 1\}$ and $\H_\rho=\ga_\rho\H_\infty$.  We will apply the
counting claim of Corollary \ref{coro:constcurvcount} with $n=3$,
$\Ga$ the arithmetic group $\overline{G}$, $D^-$ the geodesic line
with points at infinity $\alpha_0,\alpha_0^\sigma$ (whose stabiliser
is exactly $\stab_{\Ga_K}N_{\alpha_0}\cap G$, and acts with compact
quotient on $D^-$) and $D^+$ the horoball $\H_\rho$ (which is centered
at a parabolic fixed point of $\ov{G}\,$).

For every $\ga\in\SLOK$, the distance $\ell_{\ga}$ between $D^-$ and
$\ga D^+=\ga\ga_\rho\H_\infty $, when they are disjoint, is equal to
the distance between $\H_\infty$ and the geodesic line with points at
infinity $\ga_{\rho}^{-1}\ga^{-1}\alpha_0, \ga_{\rho}^{-1}\ga^{-1}
\alpha_0^\sigma$.  Hence again by \cite[Lem. 4.2]{ParPau12JMD}, if
$\ga\ga_{\rho}(1,0)=(u',v')$, then $\ell_{\ga}= \ln \big(h(\alpha_0)\,
|\n(u'-v'\alpha_0)\,|\big)$.  If $y\neq 0$, we have
$\ga\ga_{\rho}(1,0)=\frac{1}{y}\ga(x,y)$, hence
$$
\ell_{\ga}= \ln \big(\frac{h(\alpha_0)}{\n(y)}\,
|\n(u-v\alpha_0)\,|\big)
$$ 
if $\ga(x,y)=(u,v)$.  By the above convention, this is also true if
$y=0$.

As in the proof of Theorem \ref{theo:appliHerm}, we then have
\begin{align}
\ov\psi_{G,\,x,\,y}(s)& =\card\;
\big\{[\ga]\in G_{\{\alpha_0,\,\alpha_0^\sigma\}}\bs G/G_{(x,\,y)}\;:\;
\ell_\ga\leq \ln\frac{h(\alpha_0)\,\n(I_{x,\,y})\,s}{\n(y)}\big\}
+\bigO(1)
\nonumber\\ &=\frac{2}{\iota_G}\;\card\;
\big\{[\ga]\in \overline{G_{\{\alpha_0,\,\alpha_0^\sigma\}}}\,\bs
 \overline{G}/\, \overline{G_{(x,\,y)}}\;:\;
\ell_\ga\leq \ln\frac{h(\alpha_0)\,\n(I_{x,\,y})\,s}{\n(y)}\big\}
+\bigO(1)
\nonumber\\ &=\frac{2}{\iota_G}\;
[\,\overline{G_{\H_\rho}}:\overline{G_{(x,\,y)}}\,]\;
\N_{D^-,\,D^+}\big(\ln\frac{h(\alpha_0)\,\n(I_{x,\,y})\,s}{\n(y)}\big)
+\bigO(1)
\nonumber\\ &=
\frac{2\,[\,\overline{G_{\H_\rho}}:\overline{G_{(x,\,y)}}\,]}
{\iota_G}\;\frac{\Vol(\SSS^1)\,\Vol(\Ga_{D^-}\bs D^-)\,
\Vol(\Ga_{D^+}\bs D^+)}{m_G(\alpha_0)\,\Vol(\SSS^2)\,
\Vol(\Ga\bs\htr)}\;\times\nonumber\\ 
&\qquad\qquad\big(\frac{h(\alpha_0)\,\n(I_{x,\,y})\,s}
{\n(y)}\big)^2(1+\bigO(s^{-\kappa}))\;.
\label{eq:quadrbianchi}
\end{align}
By Equation \eqref{eq:calcvollonggeod}, we have $\Vol(\Ga_{D^-}\bs
D^-)= \frac{m_G(\alpha_0)\,[(\,\overline{\Ga_K})_{\alpha_0}:\,
  \overline{G}_{\alpha_0}]} {m_{\Ga_K}(\alpha_0)\,\iota_G(\alpha_0)}\,
\ell(\wh{\alpha_0})$. By Equation \eqref{eq:calccovolhororho}, we have
$\Vol(\Ga_{D^+}\bs D^+)= [\,\overline{(\Ga_K)_{\H_\rho}}:
\overline{G_{\H_\rho}}\,]\;
\frac{\sqrt{|D_K|}\,\n(y)^2}{2\,{\omega_K}\,\n(I_{x,\,y})^2}$. By
Equation \eqref{eq:modifcovolBianchi}, we have
$\Vol(\Ga\bs\htr)=\frac{[\Ga_K:G]\,|D_K|^{\frac 32}\,\zeta_K (2)}
{4\,\pi^2\,\iota_G}$. We have $[(\,\overline{\Ga_K})_{\alpha_0}:\,
\overline{G}_{\alpha_0}]=\frac{[(\Ga_K)_{\alpha_0}:G_{\alpha_0}]}
{\iota_G}$. By Equation \eqref{eq:calcindiceherm}, we have
$[\,\overline{(\Ga_K)_{\H_\rho}}:\overline{G_{\H_\rho}}\,]
[\,\overline{G_{\H_\rho}}:\overline{G_{(x,\,y)}}\,]
=\frac{\omega_K}{2}\,[(\Ga_K)_{(x,\,y)}:G_{(x,\,y)}]$.  Therefore
Equation \eqref{eq:quadrbianchi} gives the following result,
interesting in itself,
\begin{multline*}
\ov\psi_{G,\,x,\,y}(s) =
\frac{2\,\pi^2\,[(\Ga_K)_{(x,\,y)}:G_{(x,\,y)}]\,
[(\Ga_K)_{\alpha_0}:G_{\alpha_0}]\,h(\alpha_0)^2\,
\big|\log\big|\frac{\operatorname{tr}\,\wh{\alpha_0}+
    \sqrt{(\operatorname{tr}\,\wh{\alpha_0})^2-4}}{2}\big|\;\big|}
{m_{\Ga_K}(\alpha_0)\,\iota_G(\alpha_0)\,[\Ga_K:G]\,
|D_K|\,\zeta_K (2)}\,\times\\\;s^2(1+\bigO(s^{-\kappa}))\,.
\end{multline*}
Since any nonzero ideal $\mmm$ in $\OOO_K$ may be written $\mmm=
x_\mmm\OOO_K+y_\mmm\OOO_K$ for some $x_\mmm,y_\mmm$ not both zero in
$\OOO_K$ and $\ov\psi_{\mmm}= \ov\psi_{\Ga_K,\,x_\mmm,\,y_\mmm}$, the
first claim of Theorem \ref{theo:laquillebordel} follows.

\medskip
The proof of the equidistribution claim is similar to the one in
Theorem \ref{theo:posdefherm} (which was written with greater care
than necessary for this purpose).
\cqfd

{\small \bibliography{../biblio} }
%{\small \bibliography{../viitteet} }

\bigskip
{\small
\noindent \begin{tabular}{l} 
Department of Mathematics and Statistics, P.O. Box 35\\ 
40014 University of Jyv\"askyl\"a, FINLAND.\\
{\it e-mail: jouni.t.parkkonen@jyu.fi}
\end{tabular}
\medskip

\noindent \begin{tabular}{l}
D\'epartement de math\'ematique, UMR 8628 CNRS, B\^at.~425\\
Universit\'e Paris-Sud,
91405 ORSAY Cedex, FRANCE\\
{\it e-mail: frederic.paulin@math.u-psud.fr}
\end{tabular}
}

\end{document}